\documentclass[12pt, reqno]{amsart}

\usepackage{mathtools}
\mathtoolsset{showonlyrefs}
\numberwithin{equation}{section}

\usepackage{graphicx}

\usepackage{mathtools}
\mathtoolsset{showonlyrefs}

\usepackage{amssymb}
\usepackage{amsthm}
\usepackage{amscd}
\usepackage{amsmath}
\usepackage{epic,eepic}
\usepackage{mathrsfs}
\usepackage{latexsym}
\usepackage{pifont}
\usepackage[all]{xy}
\usepackage {cancel}

\theoremstyle{plain}
\newtheorem{lemma}{Lemma}[section]
\newtheorem{prop}[lemma]{Proposition}
\newtheorem{thm}[lemma]{Theorem}
\newtheorem{cor}[lemma]{Corollary}
\newtheorem{intthm}{Theorem}

\theoremstyle{definition}

\newtheorem{rem}[lemma]{Remark}
\newtheorem{defi}[lemma]{Definition}
\newtheorem{exa}[lemma]{Example}

\newtheorem*{problem}{Problem}

\newcommand{\bde}{\begin{defi}}
\newcommand{\ede}{\end{defi}\vspace{1mm}}
\newcommand{\ble}{\begin{lemma}}
\newcommand{\ele}{\end{lemma}}
\newcommand{\bpr}{\begin{prop}}
\newcommand{\epr}{\end{prop}}
\newcommand{\bt}{\begin{thm}}
\newcommand{\et}{\end{thm}}
\newcommand{\bco}{\begin{cor}}
\newcommand{\eco}{\end{cor}}
\newcommand{\bre}{\begin{rem}}
\newcommand{\ere}{\end{rem}}
\newcommand{\bex}{\begin{exa}}
\newcommand{\eex}{\end{exa}}
\newcommand{\bpf}{\begin{proof}}
\newcommand{\epf}{\end{proof}}

\newcommand{\mcB}{\mathcal{B}}

\newcommand{\mcD}{\mathcal{D}}
\newcommand{\mcE}{\mathcal{E}}
\newcommand{\mcF}{\mathcal{F}}
\newcommand{\mcG}{\mathcal{G}}
\newcommand{\mcH}{\mathcal{H}}
\newcommand{\mcI}{\mathcal{I}}

\newcommand{\mcK}{\mathcal{K}}
\newcommand{\mcL}{\mathcal{L}}
\newcommand{\mcM}{\mathcal{M}}
\newcommand{\mcN}{\mathcal{N}}
\newcommand{\mcO}{\mathcal{O}}
\newcommand{\mcP}{\mathcal{P}}
\newcommand{\mcQ}{\mathcal{Q}}

\newcommand{\mcS}{\mathcal{S}}
\newcommand{\mcT}{\mathcal{T}}
\newcommand{\mcU}{\mathcal{U}}
\newcommand{\mcV}{\mathcal{V}}

\newcommand{\mcY}{\mathcal{Y}}

\newcommand{\mbA}{\mathbb{A}}

\newcommand{\mbC}{\mathbb{C}}

\newcommand{\mbF}{\mathbb{F}}
\newcommand{\mbG}{\mathbb{G}}
\newcommand{\mbH}{\mathbb{H}}

\newcommand{\mbN}{\mathbb{N}}

\newcommand{\mbP}{\mathbb{P}}

\newcommand{\mbR}{\mathbb{R}}

\newcommand{\mbV}{\mathbb{V}}

\newcommand{\mbZ}{\mathbb{Z}}

\newcommand{\mfS}{\mathfrak{S}}

\newcommand{\mfb}{\mathfrak{b}}
\newcommand{\mfc}{\mathfrak{c}}

\newcommand{\mfg}{\mathfrak{g}}

\newcommand{\mfl}{\mathfrak{l}}
\newcommand{\mfm}{\mathfrak{m}}

\newcommand{\mfo}{\mathfrak{o}}

\newcommand{\mfs}{\mathfrak{s}}
\newcommand{\mft}{\mathfrak{t}}

\newcommand{\msD}{\mathscr{D}}
\newcommand{\msE}{\mathscr{E}}
\newcommand{\msF}{\mathscr{F}}
\newcommand{\msG}{\mathscr{G}}

\newcommand{\msO}{\mathscr{O}}

\newcommand{\msV}{\mathscr{V}}

\newcommand{\msX}{\mathscr{X}}

\newcommand{\SSP}{\vspace{3mm}}
\newcommand{\LSP}{\vspace{5mm}}
\newcommand{\mr}{\mathrm}
\newcommand{\N}{N}
\newcommand{\M}{m}
\newcommand{\DMO}{\nabla}
\newcommand{\STR}{\phi}
\newcommand{\IN}{a}
\newcommand{\EX}{d}
\newcommand{\Diag}{\rotatebox[origin=c]{45}{$\Leftarrow$}}
\newcommand{\Diagg}{\rotatebox[origin=c]{45}{$\Rightarrow$}}

\newcommand{\ma}{}

\begin{document}

\title[The irreducibility of the moduli space of curves with dormant $\mathrm{PGL}_2^{(N)}$-oper]{
The irreducibility of the moduli space  of \\ pointed stable curves with dormant $\mathrm{PGL}_2^{(N)}$-oper}
\author{Yasuhiro Wakabayashi}
\date{\today}
\markboth{}{}
\maketitle
\footnotetext{Y. Wakabayashi: 
Graduate School of Information Science and Technology, The University of Osaka, Suita, Osaka 565-0871, Japan;}
\footnotetext{e-mail: {\tt wakabayashi@ist.osaka-u.ac.jp};}
\begin{abstract}
A $\mathrm{PGL}_n^{(N)}$-oper is a specific type of flat $\mathrm{PGL}_n$-bundle on an algebraic curve in prime characteristic $p$ enhanced by an action of the sheaf of differential operators of level $N-1$. In this paper, we introduce and study a higher-level generalization of the Hitchin-Mochizuki morphism on the moduli space of $\mathrm{PGL}_n^{(N)}$-opers, defined via the characteristic polynomials of their $p^N$-curvatures.
As an application, we prove the irreducibility of the moduli space classifying  pointed stable curves equipped with dormant $\mathrm{PGL}_2^{(N)}$-opers, i.e., $\mathrm{PGL}_2^{(N)}$-opers with vanishing $p^N$-curvature.
\end{abstract}
\tableofcontents

\section{Introduction} \label{S2}

A {\it $G$-oper} for a reductive group $G$ is a specific type of    flat $G$-bundle on  an algebraic curve.
$G$-opers were  originally   introduced 
   in the study of the geometric Langlands correspondence for constructing Hecke eigensheaves  on the moduli space of bundles via quantization of  Hitchin's integrable system (cf. \cite{BeDr1}).

 For instance, 
   $\mr{GL}_n$-opers  (for $n \geq 1$) may be described as   
flat vector bundles of rank $n$ equipped with complete flags, corresponding locally  to  
  scalar differential equations of the form  
   \begin{align} \label{EQ330}
  Dy = 0, \ \ \  D = \frac{d^n}{dx^n} + a_1 \frac{d^{n-1}}{dx^{n-1}} + \cdots + a_{n-1} \frac{d}{dx} + a_n,
   \end{align}
where $x$ is 
  a local coordinate on the underlying curve  and $a_1, \cdots, a_n$ are coefficient functions.

In another important example,  $\mr{PGL}_2$-opers on a complex projective smooth curve $X$
have been studied  from the perspective of  Teichm\"{u}ller theory.
Indeed, if $X^\mr{an}$ denotes the associated Riemann surface, then a  $\mr{PGL}_2$-oper determines (up to conjugation)
a  representation  of 
 the fundamental group $\pi_1 (X^\mr{an}) \rightarrow \mr{PGL}_2 (\mbC)$, and
 such representations encode 
refinements of the complex structure of $X^\mr{an}$, known as 
 {\it projective structures}.
   This means that each $\mr{PGL}_2$-oper  corresponds to an atlas of coordinate charts  into the  Riemann sphere  $\hat{\mbC} := \mbC \sqcup \{ \infty \}$, with translation functions given by  M\"{o}bius transformations.
  A canonical example is the {\it uniformizing projective structure}, which arises from 
   the Fuchsian representation associated to the uniformization of $X^\mr{an}$.

 $G$-opers  on a hyperbolic curve  {\it in prime characteristic $p > 0$} have been studied in the context of   a characteristic-$p$ analogue  of the geometric Langlands correspondence (cf. ~\cite{BeTr}),  as well as various other topics, including 
      $p$-adic Teichm\"{u}ller theory (cf., e.g.,  ~\cite{Moc1},  \cite{Moc2},  ~\cite{JRXY}, ~\cite{JoPa}, \  ~\cite{LasPa0}, and ~\cite{LiOs}).
 A key common feature in  these developments is the notion of  $p$-curvature of a connection, which serves as an  
 invariant measuring  the obstruction to the compatibility of $p$-power structures that arise in certain spaces of infinitesimal symmetries.
 This invariant also plays a central role in 
the Grothendieck-Katz conjecture, which provides  a conjectural  criterion for  
 the algebraicity of solutions to linear differential equations 
  (cf. ~\cite{NKa3}, ~\cite{And}).
  
A $G$-oper is said to be  {\it dormant} if its $p$-curvature  vanishes.
For example, a $\mr{GL}_n$-oper is dormant if and only if the corresponding differential equation admits  a full set of solutions.
The moduli  theory  of dormant $G$-opers for general $G$  has been  developed in, e.g.,  ~\cite{Wak2}, ~\cite{Wak3}, ~\cite{Wak5},
\  ~\cite{Wak6},  ~\cite{Wak7},  ~\cite{Wak10}, and ~\cite{Wak20}.
In these works, we introduced a generalization of  dormant $\mr{PGL}_n$-opers
 using the sheaf of differential operators of finite level, as  developed  by  P. Berthelot and C. Montagnon (cf. ~\cite{Ber1}, ~\cite{Ber2}, ~\cite{Mon}).
 These generalized objects are referred to as 
{\it dormant $\mr{PGL}_n^{(\N)}$-opers}
  for $\N \in \mbZ_{> 0}$ (cf.  ~\cite{Wak6}, ~\cite{Wak20}).
  In ~\cite[Theorem A]{Hos} (or ~\cite[Theorem A, (i)]{Wak25}), it was  shown that dormant $\mr{PGL}_n^{(\N)}$-opers correspond bijectively  to {\it $F^\N$-projective structures} (cf. ~\cite[Definition 3.1]{Hos}, ~\cite[Definition 1.2.1]{Wak25}).
Within the framework of  Teichmuller theory in characteristic $p$, 
 these may be regarded as analogues of ``nice" projective structures on hyperbolic Riemann surfaces, such as 
 those  with real monodromy, including the uniformizing projective structures\footnote{
 The analogy between real-monodromy opers and dormant opers appears  in several aspects of Teichm\"{u}ler theory. 
 One such aspect  is that both  exhibit an Eichler-Shimura-type decomposition in the cohomology of symmetric products of $\mr{PGL}_2$-opers. In the case of real-monodromy opers, this phenomenon arises via an extension of  Faltings' result on the {\it discreteness} of 
 the space of real-monodromy $\mr{PGL}_2$-opers (cf. ~\cite{Fal}, ~\cite{Wak25}). In the other, the decomposition is induced by the {\it generic \'{e}taleness} of the space of dormant $\mr{PGL}_2$-opers (cf. ~\cite{Wak24}).
Another point of analogy lies in the $2$d  TQFT mentioned in this Introduction. 
Dormant $\mr{PGL}_2$-opers on a $3$-pointed projective line are completely determined by elements of $\mbF_p$, called ``radii" (i.e., eigenvalues of the residue matrices at the marked points), which satisfy a certain ``triangle inequality" (cf. ~\cite[Introduction, Theorem 1.3]{Moc2}). On the other hand, real-monodromy $\mr{PGL}_2$-opers on a pair of pants are characterized by a triangle inequality for the hyperbolic lengths of the boundary components (cf. ~\cite{Gol}).
See also ~\cite[Remark A.11]{BeTr} for a discussion of this analogy from a viewpoint related to 
 an analytic version of the geometric Langlands correspondence.
 }.

In the case of elliptic curves, {\it dormant generic Miura $\mr{PGL}_2^{(\N)}$-opers} 
and {\it $F^\N$-affine structures},  which are equivalent to each other (cf. ~\cite[Theorem B, Definition 2.1]{Hos2},  ~\cite[Theorem A, (ii), Definition 1.2.1]{Wak25}), play the role of their counterparts. 
These are further described by  {\it Igusa-structures of level $p^\N$}, which are, by definition,  generators of the iterated Verschiebung  (cf. ~\cite[Definition 12.3.1]{KaMa}). 
Such a description allows us  to interpret  moduli spaces of dormant $\mr{PGL}_2^{(\N)}$-opers   as {\it higher-genus generalizations of  the classical  Igusa curves} (cf. Section \ref{S300}).

\vspace{-10mm}

\hspace{-2mm} 
 \includegraphics[width=20cm,bb=0 0 1112 380,clip]{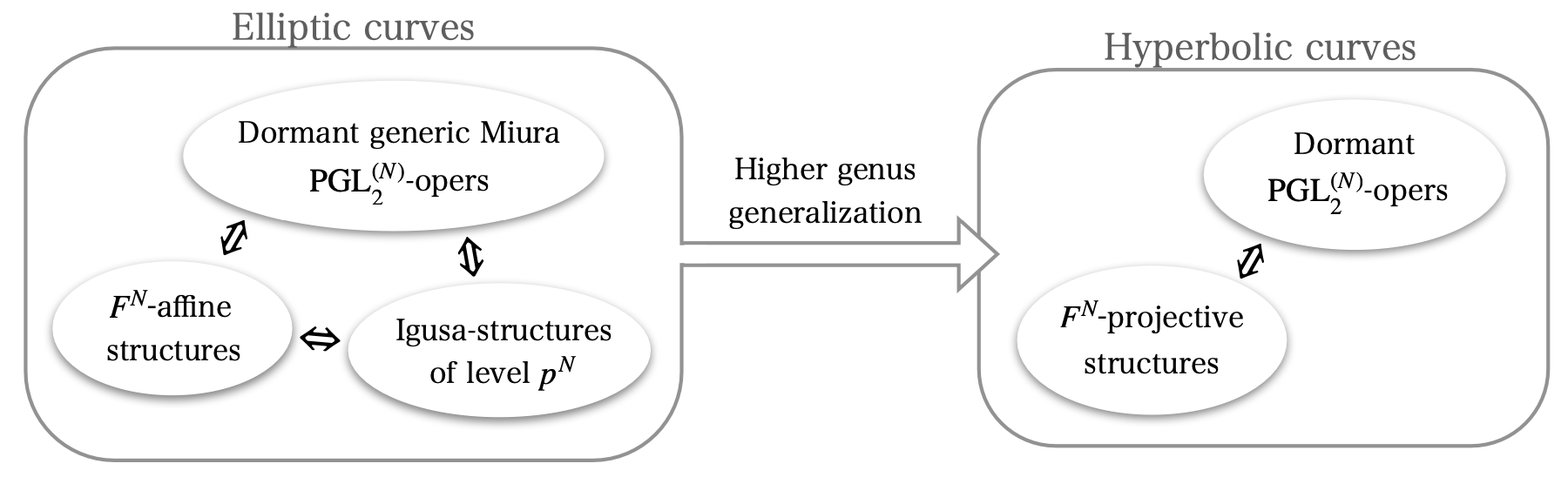}

\vspace{6mm}

Now, we fix a positive integer $n$ satisfying $1 < n < p$, 
an algebraically closed  field $k$ of characteristic $p$, and a pair of nonnegative integers $(g, r)$ satisfying  $2g-2+r > 0$.
Let $\overline{\mcM}_{g, r}$  denote the moduli stack of $r$-pointed stable curves of genus $g$ over $k$.
Given a positive integer $\N$ and an $r$-tuple $\rho := (\rho_i)_{i=1}^r$ of elements of  
a certain finite set $\Xi_{n, \N}$ (cf. \eqref{Er3}),
we define the moduli category 
\begin{align}
\mcO p_{\N, \rho, g, r} \ \left(\text{resp.,} \ 
\mcO p^{^\mr{Zzz...}}_{\N, \rho, g, r}; \text{resp.,} \ \mcO p^{^\mr{nilp}}_{\N, \rho, g, r} \right)
\end{align}
 that  classifies  pairs $(\msX, \msE^\spadesuit)$
 consisting of a pointed curve $\msX$ in $\overline{\mcM}_{g, r}$ and  a $\mr{PGL}_2^{(\N)}$-oper 
(resp., a dormant $\mr{PGL}_2^{(\N)}$-oper; resp., a $p^\N$-nilpotent $\mr{PGL}_2^{(\N)}$-oper) $\msE^\spadesuit$ on $\msX$ of radii $\rho$.
In particular, we have a sequence of inclusions
$\mcO p^{^\mr{Zzz...}}_{\N, \rho, g, r} \subseteq \mcO p^{^\mr{nilp}}_{\N, \rho, g, r} \subseteq \mcO p_{\N, \rho, g, r}$.

According to   ~\cite[Theorem B]{Wak20}, $\mcO p^{^\mr{Zzz...}}_{\N, \rho, g, r}$ can be represented by  a (possibly empty) proper  Deligne-Mumford stack   over $k$, and the projection $\mcO p_{\N, \rho, g, r}^{^\mr{Zzz...}} \rightarrow \overline{\mcM}_{g, r}$ given by $(\msX, \msE^\spadesuit) \mapsto \msX$ is finite.
Moreover, if either ``$n=2$" or ``$\N =1$ and $2n < p$" is fulfilled, then this projection is generically \'{e}tale (cf. ~\cite[Theorem G]{Wak5}, ~\cite[Theorem C, (i)]{Wak20}).
These geometric properties play crucial roles  in establishing a factorization rule of 
 the total number of dormant $\mr{PGL}_n^{(\N)}$-opers on a general curve according to  various gluing  procedures for underlying curves.
This factorization  satisfies a $2$d TQFT, i.e., a $2$-dimensional topological quantum field theory (cf.  ~\cite[Theorem C, (ii)]{Wak20}).
In the special  case  $\N =1$, this $2$d TQFT is closely  related to the Gromov-Witten theory of Grassmannian varieties and to  the conformal field theory  for affine Lie algebras (cf. ~\cite{Wak5}).

The generic \'{e}taleness for $n=2$ also gives rise to  a certain arithmetic lifting  of each dormant $\mr{PGL}_2^{(\N)}$-oper to characteristic $p^\N$, called the {\it canonical diagonal lifting} (cf.  ~\cite[Theorem D]{Wak20}).
Thus, 
a detailed investigation of the  moduli stack $\mcO p^{^\mr{Zzz...}}_{\N, \rho, g, r}$ serves as a bridge between the theory of opers (as well as of linear differential equations) in arithmetic settings and  various other types  of enumerative geometry, thereby offering a unified framework for their mutual development.

\hspace{0mm} 
 \includegraphics[width=19cm,bb=0 0 1012 200,clip]{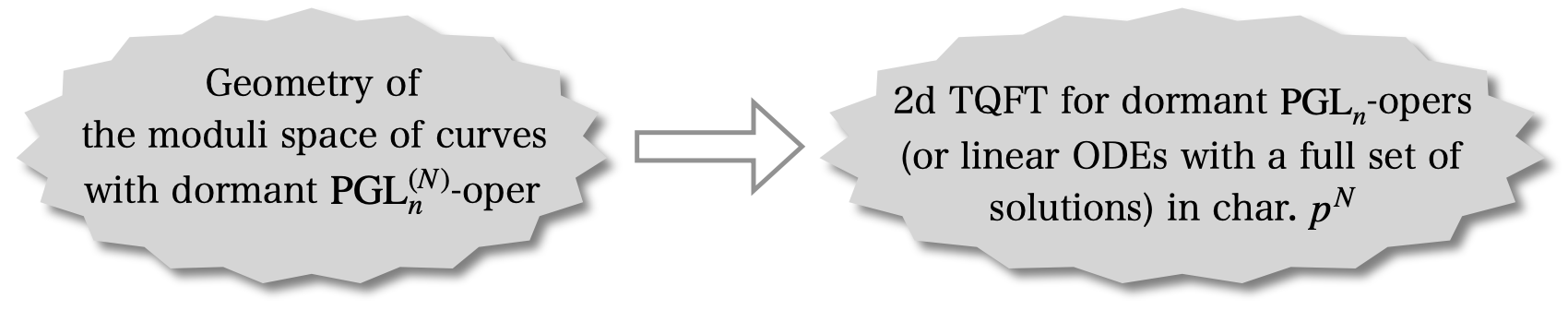}

\vspace{3mm}







In this paper, we focus on and prove certain fundamental properties of the stack $\mcO p^{^\mr{Zzz...}}_{\N, \rho, g, r}$, regarded  as a characteristic-$p$  analogue of the moduli space of Riemann surfaces (equipped with uniformizing projective structures), as well as a higher-genus generalization of Igusa curves. 
Our central result is the {\it irreducibility of $\mcO p^{^\mr{Zzz...}}_{\N, \rho, g, r}$}.
The irreducibility of 
the moduli space of smooth projective curves  over an arbitrary field
  was  established in ~\cite{DeMu} via   the contractibility of  Teichm\"{u}ller spaces.
  Furthermore,
in the context of $p$-adic Teichm\"{u}ller theory,
the same property for   $\mcO p^{^\mr{Zzz...}}_{\N, \rho, g, r}$  with  $\N =1$ has already been shown in ~\cite[Chapter II, Theorem 2.8]{Moc2} 
 (cf. ~\cite[Theorem 12.6.1, Corollary 12.6.2]{KaMa} for the case  of Igusa curves).


A key tool in our approach is a higher-level generalization of the so-called Hitchin-Mochizuki morphism, previously studied  in ~\cite{BeTr}, ~\cite{LasPa0}, ~\cite{JoPa}.
Using a certain relative affine space $B_{\mr{univ}, 0}^\flat$ over  $\overline{\mcM}_{g, r}$ (cf. \eqref{Ey34}),
we construct a morphism
\begin{align}
\mr{HM}^{(\N)}_{\widetilde{\rho}} : \mcO p_{\N +1, \widetilde{\rho}, g, r} \rightarrow B^\flat_{\mr{univ}, 0},
\end{align}
for each $\widetilde{\rho} \in \Xi_{n, \N+1}^{r\ma} \left( = \Xi_{n, \N+1} \times \cdots \times \Xi_{n, \N+1}\right)$ lifting  $\rho$ (cf. \eqref{Er247}), 
defined  by taking the characteristic polynomials of the $p^\N$-curvatures of $\mr{PGL}_{n}^{(\N +1)}$-opers.
We refer  to  this morphism as the {\bf level-$\N$ Hitchin-Mochizuki morphism} (cf. Definition \ref{Defgg9}).

Under certain assumptions,
we prove that 
this morphism satisfies some natural and desirable properties; see Theorem \ref{Cor20} for the precise statements.
Since $\mcO p^{^\mr{nilp}}_{\N +1, \widetilde{\rho}, g, r}$ can be obtained as the inverse image 
 of the zero section $0_{\mr{univ}} : \overline{\mcM}_{g, r} \rightarrow B^\flat_{\mr{univ}, 0}$,
the level-reduction morphism   $\Lambda_{\widetilde{\rho}, g, r} : \mcO p^{^\mr{nilp}}_{\N +1, \widetilde{\rho}, g, r} \rightarrow \mcO p^{^\mr{Zzz...}}_{\N, \rho, g, r}$
turns out to be  (finite, flat, and) surjective (cf. Theorem \ref{Th665}).
By combining this fact  
and the standard degeneration technique for  underlying curves,
 we establish some  desired assertions   (including the connectedness  of $\mcO p^{^\mr{nilp}}_{\N, \rho, g, r}$)  via induction on the dimension.
Our main results are summarized   as follows.

\SSP
\begin{intthm}[cf. Theorem \ref{Cor11}]  \label{ThA}
Suppose that $n=2$, and 
let $\rho$ be an element of $\Xi_{2, \N}^{r\ma}$, where  we set $\rho := \emptyset$ when  $r =0$.
Then, the following assertions hold:
\begin{itemize}
\item[(i)]
If $\mcO p_{\N, \rho, g, r}^{^\mr{nilp}} \neq \emptyset$,
then any two irreducible components of $\mcO p_{\N, \rho, g, r}^{^\mr{nilp}}$ intersect.
In particular,  $\mcO p_{\N, \rho, g, r}^{^\mr{nilp}}$ is connected (when nonempty).
\item[(ii)]
If $\mcO p^{^\mr{Zzz...}}_{\N, \rho, g, r} \neq \emptyset$,
then $\mcO p^{^\mr{Zzz...}}_{\N, \rho, g, r}$ is irreducible.
 \end{itemize}
 \end{intthm}
\SSP

In the above theorem, we assume the nonemptiness of  the relevant moduli spaces. 
For results concerning this property,  we refer to Proposition \ref{Prop77}, Corollary \ref{Cor15}, and  ~\cite[Corollary 6.28, (ii), Proposition 10.18]{Wak20}.

In Section \ref{S8},
we present several  direct applications of the irreducibility  of $\mcO p^{^\mr{Zzz...}}_{\N, \rho, g, r}$.
These include
applications   to 
the generalized  Verschiebung map on the moduli space of stable vector bundles  (cf. Proposition \ref{Prop491}, Corollary \ref{Cor4}), and to  certain arithmetic deformations 
of  dormant $\mr{PGL}_2^{(\N)}$-opers (cf. Propositions \ref{Proptt6}, \ref{Proptt5}).
Regarding the former topic, we have also established  an explicit  and combinatorial procedure for computing   the generic  degree of the generalized  Verschiebung map, based on the results of   this paper  (cf. ~\cite{KoWa}).

Moreover,  motivated by the present work, we intend to pursue further research from a coding-theoretic perspective
on the tower of stacks 
\begin{align} 
\cdots \xrightarrow{} \mcO p^{^\mr{Zzz...}}_{N, \rho_\N, g, r} \xrightarrow{} \cdots   \xrightarrow{}   \mcO p^{^\mr{Zzz...}}_{3, \rho_2, g, r}\xrightarrow{}  \mcO p^{^\mr{Zzz...}}_{2, \rho_2, g, r} \xrightarrow{}  \mcO p^{^\mr{Zzz...}}_{1, \rho_1, g, r}.
\end{align}
(cf. \eqref{Eq12})
formed by varying the level $\N$ (cf. ~\cite{Wak80}).

\LSP
\subsection*{Notation and Conventions} \label{SS0}

Throughout  this  paper,
all schemes (and Deligne-Mumford stacks) are assumed  to be  locally noetherian.

We fix 
an odd prime $p$,  an algebraically closed  field $k$ of characteristic $p$, and 
 a positive integer  $n$  with $1 < n < p$.

Let  $\mr{GL}_n$ (resp., $\mr{PGL}_n$) denote  the general (resp.,  projective) linear group of $k^{\oplus n}$. 
The algebraic group $\mr{PGL}_n$ contains 
the Borel subgroup $B$  consisting of  the images of upper triangular invertible  $n \times n$ matrices under  the natural quotient $\mr{GL}_n \twoheadrightarrow \mr{PGL}_n$.
The Lie algebra of $\mr{PGL}_n$ (resp., $B$) is denoted by 
$\mfg := \mr{Lie} (\mr{PGL}_n)$ (resp.,  $\mfb := \mr{Lie}(B)$).

Let $S$ be a scheme (or more generally, a Deligne-Mumford stack).
We write   $\mcS ch_{/S}$ for  the category of $S$-schemes.
For a vector bundle (i.e., a locally free coherent sheaf)  $\mcF$ on $S$, we denote by $\mbV (\mcF)$ the relative affine scheme over $S$  associated to $\mcF$, i.e., the spectrum $\mcS pec (\mr{Sym}_{\mcO_S} (\mcV^\vee))$ of the symmetric algebra $\mr{Sym}_{\mcO_S}(\mcF^\vee)$ of  its dual $\mcF^\vee$ over $\mcO_S$.

Also, denote by $\mcE nd (\mcF)\left( :=\mcE nd_{\mcO_S} (\mcF)\right)$ the sheaf of  $\mcO_S$-linear endomorphisms of $\mcF$,
and by $\overline{\mcE nd} (\mcF)$ the cokernel of the diagonal embedding $\mcO_S \hookrightarrow \mcE nd (\mcF)$.
In particular,  if $\mcF$ has rank $n$, then  $\mcE nd (\mcF)$ (resp., $\overline{\mcE nd} (\mcF)$) can be identified with the adjoint vector bundle associated to the $\mr{GL}_n$-bundle (resp., $\mr{PGL}_n$-bundle) determined by  $\mcF$.

\vspace{10mm}
\section{Preliminaries on $\mcD$-modules of finite  level} \label{S1}
\LSP

In this section, we review basic materials  on  $\mcD$-modules of finite level in positive characteristic, including their logarithmic generalizations.
For detailed accounts, we refer  the reader to 
 ~\cite{Ber1}, ~\cite{Ber2}, ~\cite{Mon},  and ~\cite{Wak20}.

\LSP
\subsection{$\mcD^{(\M)}$-modules in the logarithmic setting} \label{SS11}

Let $(g, r)$ be a pair of nonnegative integers with $2g-2+r > 0$ and 
$S$  a scheme (or more generally, a Deligne-Mumford stack) over $k$.
Also, let   $\msX := (f : X \rightarrow S, \{ \sigma_i \}_{i=1}^r)$  be  an $r$-pointed stable curve of genus $g$ over $S$, where
$f: X\rightarrow S$ denotes the underlying prestable curve and 
each  $\sigma_i$ ($i=1, \cdots, r$) specifies  the $i$-th marked point $S \rightarrow X$.

As shown in  ~\cite[Theorem 2.6]{FKa}, 
 both $S$ and $X$  carry natural log structures.
 We denote the resulting log structures by $S^\mr{log}$ and $X^\mr{log}$, respectively.
 The morphism $f$ extends to a morphism of log schemes  $f^\mr{log} : X^\mr{log} \rightarrow S^\mr{log}$, under which $X^\mr{log}$ becomes  a log curve over $S^\mr{log}$ in the sense of ~\cite[Definition 2.1]{Wak20}.
 We denote the sheaf of logarithmic $1$-forms on $X^\mr{log}$ over $S^\mr{log}$ by  $\Omega_{X^\mr{log}/S^\mr{log}}$  and its dual  by 
 $\mcT_{X^\mr{log}/S^\mr{log}}$.

Next, let 
 $\M$ be  a nonnegative integer.
Let  $F_S$ (resp., $F_X$) denote  the absolute Frobenius endomorphism of $S$ (resp., $X$).
We define 
the {\bf $(\M +1)$-st Frobenius twist} of $X$ over $S$ to be 
the base-change  $f^{(\M +1)} : X^{(\M +1)}\left(:=S \times_{F_S^{\M +1}, S} X\right) \rightarrow S$ of the $S$-scheme $f: X \rightarrow S$ along  the $(\M +1)$-st iterate $F_S^{\M +1}$ of $F_S$.
The associated  morphism $F_{X/S}^{(\M +1)}  \left(:= (f, F_X^{\M +1}) \right) : X \rightarrow X^{(\M +1)}$ is called the {\bf $(\M +1)$-st relative Frobenius morphism} of $X$ over $S$.
The $i$-th marked point $\sigma_i$ ($i=1, \cdots, r$) determines a marked point $\sigma^{(\M +1)} : S \rightarrow X^{(\M +1)}$ of $X^{(\M +1)}$, and the resulting collection
\begin{align}
\msX^{(\M +1)} := (f^{(\M +1)}: X^{(\M +1)}\rightarrow S, \{ \sigma_i^{(\M +1)}\}_{i=1}^r)
\end{align}
 forms an $r$-pointed stable curve.

To simplify the notation, we write $F_{X/S} := F_{X/S}^{(1)}$, and write $\Omega^{(m+1)}:= \Omega_{X^{(m+1)\mr{log}}/S^\mr{log}}$,
$\mcT^{(m+1)} := \mcT_{X^{(m+1)\mr{log}}/S^\mr{log}}$.
Also, we occasionally write $X^{(0)}$, $f^{(0)}$, $\Omega^{(0)}$, and $\mcT^{(0)}$ instead of $X$, $f$, $\Omega$, and $\mcT$, respectively, for notational convenience.
 
For each $\ell \in \mbZ_{\geq 0}$,
we denote by
$\mcD^{(m)}_{X^\mr{log}/S^\mr{log}, \leq \ell}$ 
the sheaf of logarithmic differential operators of level $m$ and order $\leq \ell$, 
as introduced in ~\cite[D\'{e}finition 2.3.1]{Mon}.
We   set
\begin{align}
\mcD_{X^\mr{log}/S^\mr{log}}^{(m)}
 := \bigcup_{\ell \in \mbZ_{\geq 0}} \mcD^{(m)}_{X^\mr{log}/S^\mr{log}, \leq \ell},
\end{align}
which is equipped with  a structure of (non-commutative) algebras.
When there is no fear of confusion, we 
abbreviate
 $\mcD_{\leq \ell}^{(\M)} := \mcD^{(m)}_{X^\mr{log}/S^\mr{log}, \leq \ell}$ and $\mcD_{\leq \infty}^{(\M)} \  \left(\text{or simply} \  \mcD_{}^{(\M)} \right) := \mcD_{X^\mr{log}/S^\mr{log}}^{(m)}$. 

Let $x$ be  a (locally defined) logarithmic coordinate on  $X$ (relative to $S$).
Associated to such a coordinate  $x$, one has  
 a  local $\mcO_X$-basis $\{ \partial^{\langle  j \rangle}\}_{j \leq \ell}$ of ${^L}\mcD_{\leq \ell}^{(\M)}$ (cf. ~\cite[Section 1.2.3]{Mon}).
For  each  $j \in \mbZ_{\geq 0}$,  the pair  of nonnegative integers $(q_j, r_j)$ is  uniquely determined by the condition $j = p^\M \cdot q_j + r_j$ with  $0 \leq r_j < p^\M$.
Then,
the identity 
\begin{align} \label{e364}
\partial^{\langle j' \rangle} \cdot \partial^{\langle j'' \rangle} = \sum_{j = \mr{max} \{ j', j'' \} }^{j' + j''} \frac{j!}{(j' + j'' - j)! \cdot (j - j')! \cdot  (j-j'')!} \cdot \frac{q_{j'}! \cdot q_{j''}!}{q_j !} \cdot \partial^{\langle j \rangle}
\end{align}
holds for any nonnegative  integers $j'$ and $j''$ (cf. ~\cite[Lemme 2.3.4]{Mon}).
In particular, 
the sections  $ \partial^{\langle  j \rangle}$ satisfy the commutativity relation  
\begin{align} \label{eQ10}
\partial^{\langle j' \rangle} \cdot \partial^{\langle j'' \rangle} = \partial^{\langle j'' \rangle} \cdot \partial^{\langle j' \rangle}.
\end{align}

Given each  $\ell \in \mbZ_{\geq 0} \sqcup \{ \infty \}$,
we shall write  ${^L}\mcD_{\leq \ell}^{(m)}$ (resp., ${^R}\mcD_{\leq \ell}^{(m)}$) for the sheaf $\mcD_{\leq \ell}^{(m)}$ endowed with a structure of $\mcO_X$-module arising from left (resp., right) multiplication by sections of $\mcD_{\leq 0}^{(m)} \left(= \mcO_X \right)$.
A {\bf (left) $\mcD^{(\M)}$-module structure} on an $\mcO_X$-module $\mcF$ is defined as  a left $\mcD^{(\M)}$-action $\nabla : {^L}\mcD^{(\M)} \rightarrow \mcE nd_{\mcO_S}(\mcF)$ on $\mcF$ extending its $\mcO_X$-module structure.
When it is necessary to emphasize the level ``$\M$", we write $\nabla^{(\M)}$ instead of $\nabla$.
Note that giving such a structure for $\M = 0$ amounts to  giving  an {\bf $S^\mr{log}$-connection} in the sense of ~\cite[Definition 4.1]{Wak5}, i.e., an $f^{-1}(\mcO_S)$-linear morphism $\mcF \rightarrow \Omega \otimes \mcF$ satisfying the usual Leibnitz rule.
An $\mcO_X$-module equipped with a $\mcD^{(\M)}$-module structure (resp., an $S^\mr{log}$-connection) is called a {\bf (left) $\mcD^{(\M)}$-module} (resp., a {\bf flat module} on $X^\mr{log}/S^\mr{log}$).

For an $\mcO_X$-module $\mcF$, we define  the tensor product $\mcD_{\leq \ell}^{(m)} \otimes \mcF := {^R}\mcD_{\leq \ell}^{(m)} \otimes \mcF$ and equip it 
with the $\mcO_X$-module structure on $\mcD_{\leq \ell}^{(m)}$  given by left multiplication.

Given a $\mcD^{(\M)}$-module $(\mcF, \nabla)$, we define 
\begin{align} \label{eQ1}
\mcS ol (\nabla^{(m)})
\end{align}
to be  the subsheaf of $\mcF$ on which $\mcD_+^{(\M)}$ acts as zero, where $\mcD_+^{(\M)}$ denotes the kernel of the canonical projection $\mcD^{(\M)} \twoheadrightarrow \mcO_X$.
The sheaf $\mcS ol (\nabla^{(m)})$ can be regarded as an $\mcO_{X^{(\M +1)}}$-module via the underlying homeomorphism of $F_{X/S}^{(\M +1)}$.

Let $(\mcF, \DMO^{(\M)})$ be a $\mcD^{(\M)}$-module. 
For an  integer $a$ with $0\leq a \leq m$, 
 the composite
\begin{align} \label{eQ81}
\nabla^{(m)\Rightarrow (a)} : {^L}\mcD^{(a)} \rightarrow \mcE nd_{\mcO_S} (\mcF)
\end{align}
of $\nabla^{(m)}$ and the natural morphism ${^L}\mcD^{(a)} \rightarrow {^L}\mcD^{(m)}$ defines a $\mcD^{(a)}$-module structure on $\mcF$.
In particular,   we obtain a $\mcD^{(a)}$-module $(\mcF, \nabla^{(m) \Rightarrow (a)})$.
For an  integer $a$ with $0\leq a \leq m+1$,  we  
set 
\begin{align} \label{dE48}
\mcF^{[a]} := \begin{cases} \mcF & \text{if $a = 0$};  \\
\mcS ol (\DMO^{(\M) \Rightarrow (a-1)}) & \text{if $a > 0$},
\end{cases}
\end{align}
which  is an $\mcO_{X^{(a)}}$-module.

In what follows, we  define an $S^\mr{log}$-connection on the $\mcO_{X^{(a)}}$-module $\mcF^{[a]}$ for $a = 0, \cdots, \M$ (cf. ~\cite[Section 2.6]{Wak20}).
We begin by setting  $\nabla^{[0]}$ (or $(\nabla^{(\M)})^{[0]}$) $:= \nabla^{(\M) \Rightarrow (0)}$.
Now, let us fix  an integer  $a \in \{1, \cdots, \M\}$.
Since  
$\mcD_{\leq p^a -1}^{(a)} = \mr{Im}\left(\mcD^{(a-1)}_{\leq p^{a}} \rightarrow \mcD^{(a)}_{\leq p^{a}}\right)$ and 
 $\mcD_{\leq p^a }^{(a)}/\mcD_{\leq p^a -1}^{(a)} = \mcT^{\otimes p^a}$, we obtain 
an exact sequence
\begin{align} \label{dE45}
\mcD^{(a-1)}_{\leq p^{a}}\cap \mcD^{(a-1)}_{+} \rightarrow
 \mcD^{(a)}_{\leq p^{a}}\cap \mcD^{(a)}_{+}  \xrightarrow{\delta} \mcT^{\otimes p^{a}} \left(= F_{X/S}^{(a)*}(\mcT^{(a)}) \right) \rightarrow 0,
\end{align}
where the first arrow is the morphism obtained from  the natural  morphism $\mcD^{(a-1)} \rightarrow \mcD^{(a)}$.
Let us take a local section $\partial$ of $\mcT^{(a)}$.
There exists locally  a section $\widetilde{\partial}$ in $\delta^{-1} ((F_{X/S}^{(a)})^{-1}(\partial)) \left(\subseteq \mcD^{(a)}_{\leq p^{a}}\cap \mcD^{(a)}_{+} \right)$.
The exactness of \eqref{dE45} and the definition of $\mcF^{[a]}$ imply  that 
 the $f^{-1}(\mcO_S)$-linear endomorphism $\DMO_\partial^{[a]} := \DMO^{(\M)\Rightarrow {(a)}} (\widetilde{\partial})$ of $\mcF^{[a]}$ does not depend on the choice of  $\widetilde{\partial}$ (i.e., depends only on $\partial$).
 Hence, we obtain a well-defined  morphism 
 \begin{align} \label{dE60}
 \nabla^{[a]}  \left( \text{or} \ (\DMO^{(\M)})^{[a]} \right) : \mcF^{[a]} \rightarrow \Omega^{(a)} \otimes \mcF^{[a]}
 \end{align}
 determined  by assigning $\partial \mapsto \DMO_\partial^{[a]}$.
 This morphism satisfies the Leibniz rule and hence defines 
  an $S^\mr{log}$-connection.
We thus obtain  a flat module
\begin{align} \label{J16}
(\mcF^{[a]}, \DMO^{[a]})
\end{align}
on the log curve $X^{(a)\mr{log}}/S^\mr{log}$.
 The collection $\{ \mcF^{[a]} \}_{0 \leq a \leq \M +1}$ forms  a decreasing filtration on $\mcF$ such that $\mr{Ker}(\nabla^{[a]}) = \mcF^{[a+1]}$ for every $a = 0,1, \cdots, \M$.

\LSP
\subsection{The relative $p$-Hitchin base} \label{SS5r1}

We use  the notation $\mfc$ to denote the GIT quotient $\mfg /\!/\mr{PGL}_n$ of $\mfg$ by the adjoint $\mr{PGL}_n$-action.
The $k$-scheme  $\mfc$ admits a $\mbG_m$-action coming from the homotheties on $\mfg$.
If $\mcL$ is a line bundle on either $X$ or $X^{(\M+1)}$, then there exists a natural identification between  the affine space associated to $\bigoplus_{\ell =2}^{n} \mcL^{\otimes \ell}$ and $\mcL^\times \times^{\mbG_m} \mfc$, where $\mcL^\times$ denotes the $\mbG_m$-bundle corresponding to $\mcL$.
That is to say,  there exists a canonical identification
\begin{align} \label{Ey5}
\mcL^\times \times^{\mbG_m} \mfc \cong  \mbV (\bigoplus_{\ell =2}^n \mcL^{\otimes \ell}).
\end{align}

Denote by 
\begin{align} \label{Er24}
B_\msX \ \left(\text{resp.,} \ B^\flat_\msX \right)
\end{align}
the set-valued
contravariant functor on $\mcS ch_{/S}$ which assigns, to any $S$-scheme $s : S' \rightarrow S$,
 the set of morphisms $S' \times_S X \rightarrow (\Omega^{\otimes p^{\M +1}})^\times \times^{\mbG_m} \mfc$ over $X$
 (resp., $S' \times_S X^{(\M +1)} \rightarrow (\Omega^{(m+1)})^\times \times^{\mbG_m} \mfc$ over $X^{(\M +1)}$).
This functor 
 can be represented by a smooth $S$-scheme.
Indeed, the identification \eqref{Ey5} yields  an isomorphism of $S$-schemes
\begin{align} \label{Er99}
B_\msX \xrightarrow{\sim} \mbV (f_* (\bigoplus_{\ell=2}^{n} \Omega^{\otimes p^{\M +1} \cdot \ell}))  \ \left(\text{resp.,} \  B^{\flat}_\msX  \xrightarrow{\sim} \mbV (f_*^{(\M +1)} (\bigoplus_{\ell=2}^{n} (\Omega^{(\M +1)})^{\otimes \ell}))   \right),
\end{align}
which  is functorial with respect to $S$.
In particular, the relative dimension of $B^\flat_\msX$ is computed by the following chain of equalities:
\begin{align} \label{Eu3}
\mr{dim} (B^\flat_\msX/S) &= 
\mr{rk} \left( f_*^{(\M +1)} (\bigoplus_{\ell=2}^{n} (\Omega^{(\M +1)})^{\otimes \ell})\right) \\
& = \sum_{\ell =2}^n \mr{rk} \left(f_{*}^{(\M +1)} ((\Omega^{(\M +1)})^{\otimes \ell}) \right) \notag \\
& =  \sum_{\ell =2}^n \left(\ell \cdot (2g-2+r) + 1-g \right) \\
& =  (n^2 -1)(g-1) + \frac{1}{2}(n^2 +n -2)r.\notag
\end{align}

Since $F^{(m+1)*}_{X/S}(\Omega^{(m+1)}) \cong \Omega^{\otimes p^{\M +1}}$, the pull-back via $F_{X/S}^{(m+1)}$ gives rise to  
an $S$-morphism 
\begin{align} \label{Er91}
\iota_\msX : B^\flat_\msX \rightarrow B_\msX,
\end{align}
and it makes the following square diagram commute:
\begin{align} \label{Er208}
\vcenter{\xymatrix@C=46pt@R=36pt{
 B_\msX^\flat \ar[r]^-{\eqref{Er99}}_-{\sim} \ar[d]_-{\iota_\msX} & \mbV (f_*^{(\M +1)} (\bigoplus_{\ell=2}^{n} (\Omega^{(\M +1)})^{\otimes \ell})) \ar[d] \\
 B_\msX \ar[r]_-{\eqref{Er99}}^-{\sim} & \mbV (f_* (\bigoplus_{\ell=2}^{n} \Omega^{\otimes p^{\M +1} \cdot \ell})),
 }}
\end{align}
where the right-hand vertical arrow arises from the direct sum $\bigoplus_{\ell =2}^n \omega_\ell$ of the morphisms
 \begin{align}
 \omega_\ell : (\Omega^{(m+1)})^{\otimes \ell} \rightarrow \left( F_{X/S*}^{(m+1)} (F_{X/S}^{(m+1)*} ((\Omega^{(m+1)})^{\otimes \ell})) \cong \right)
F_{X/S*}^{(m+1)} (\Omega^{\otimes p^{m+1}\cdot \ell})
\end{align}
induced by 
the adjunction relations ``$F_{X/S}^{(m+1)*} (-) \dashv F_{X/S*}^{(m+1)} (-)$".

\SSP
\bpr \label{Prop8}
\begin{itemize}
\item[(i)]
For each $\ell =2, \cdots, n$,
the direct image 
\begin{align} \label{Ey10}
f^{(\M +1)}_* (\omega_\ell) : f^{(\M +1)}_* ((\Omega^{(m+1)})^{\otimes \ell}) \rightarrow f_* (\Omega^{\otimes p^{m+1}\cdot \ell})
\end{align}
of $\omega_\ell$ via $f^{(\M +1)}$
 is injective, and its cokernel is flat over $S$.
\item[(ii)]
The morphism $\iota_\msX$ is a closed immersion.
\end{itemize}
\epr
\begin{proof}
Assertion (i) follows from  an argument similar to the proof of  ~\cite[Proposition 3.23]{Wak5}.
Assertion (ii) follows from assertion (i) together with the commutativity of \eqref{Er208}, since the right-hand vertical arrow is obtained by applying   the functor $\mcS pec (\mr{Sym}_{\mcO_S} (-)^\vee)$ to the morphism  $\bigoplus_{\ell =2}^n f^{(\M +1)}_* (\omega_\ell)$.
\end{proof}
\SSP

Next,  we shall write
\begin{align} \label{Et61}
B_{\msX, 0}^\flat
\end{align}
for the subfunctor of $B_\msX^\flat$ that classifies, for each $S$-scheme $s : S' \rightarrow S$,  morphisms $S' \times_S X^{(\M +1)} \rightarrow (\Omega^{(\M+1)})^\times \times^{\mbG_m} \mfc$ 
inducing the zero section of $\mfc$ when restricted to every marked point $\sigma_i^{(\M +1)}$ of $\msX^{(\M +1)}$ under the identification $\sigma^{(\M +1)*}_i ((\Omega^{(\M +1)})^\times) \cong \mbG_m \times_k S$ arising from the residue isomorphism.
This functor can be represented by a closed subscheme of $B_{\msX}^\flat$, and (the resp'd portion of) the isomorphism \eqref{Er99}  restricts to
an isomorphism
\begin{align} \label{Et111}
B_{\msX, 0}^\flat \xrightarrow{\sim} \mbV (f^{(\M +1)}_* (\bigoplus_{\ell =2}^n (\Omega^{(\M +1)})^{\otimes \ell}(-D^{(\M +1)}))),
\end{align}
where $D^{(\M +1)}$ denotes the relative effective divisor on $X^{(\M +1)}$ defined as the union of $\sigma^{(\M +1)}_i$'s ($i=1, \cdots, r$).
Similarly to  \eqref{Eu3},
the relative dimension of  $B_{\msX, 0}^\flat/S$ is  given  by 
\begin{align} \label{Eu4}
\mr{dim}(B_{\msX, 0}^\flat/S) = (n^2 -1) (g-1) + \frac{1}{2} n(n-1)r.
\end{align}

\LSP
\subsection{$p^{m+1}$-curvature of a $\mcD^{(m)}$-module} \label{SS51}

Note that the image of the natural morphism $\mcD_{\leq p^{m+1}}^{(m)} \rightarrow \mcD_{\leq p^{m+1}}^{(m+1)}$ coincides with $\mcD_{\leq p^{m+1}-1}^{(m+1)}$.
If $\varpi : \mcD_{\leq p^{m+1}}^{(m)} \rightarrow \mcD_{p^{m+1}-1}^{(m+1)}$ denotes the resulting surjection, then the composite
\begin{align} \label{Eq1}
\mcD_{\leq p^{m+1}-1}^{(m)} \xrightarrow{\mr{inclusion}} \mcD_{\leq p^{m+1}}^{(m)} \xrightarrow{\varpi} \mcD_{\leq p^{m+1}-1}^{(m+1)}
\end{align}
is an isomorphism.
The composite $\varpi' : \mcD_{\leq p^{m+1}}^{(m)} \rightarrow \mcD_{\leq p^{m+1}-1}^{(m)}$ of $\varpi$ and the inverse of  \eqref{Eq1} defines a split surjection of the short exact sequence
\begin{align} \label{Eq2}
0 \rightarrow \mcD_{p^{m+1}-1}^{(m)} \rightarrow \mcD_{\leq p^{m+1}}^{(m)} \rightarrow \left(\mcD_{\leq p^{m+1}}^{(m)}/\mcD_{\leq p^{m+1}-1}^{(m)}= \right) \mcT^{\otimes p^{m+1}} \rightarrow 0.
\end{align}
Thus, we obtain an  $\mcO_X$-linear morphism
\begin{align} \label{Eq40}
\psi_{X^\mr{log}/S^\mr{log}} :
 \mcT^{\otimes p^{m+1}}  \hookrightarrow \mcD_{\leq p^{m+1}}^{(m)} \xrightarrow{\mr{inclusion}} \mcD^{(m)},
\end{align}
where the first arrow denotes the split injection of  \eqref{Eq2} corresponding to $\varpi'$.
This morphism coincides, via the adjunction relation $F^{(m+1)*}_{X/S} (-) \dashv F_{X/S*}^{(m+1)} (-)$, with the $p^{m+1}$-curvature map $\mcT^{(m+1)} \rightarrow F_{X/S*}^{(m+1)} (\mcD^{(m)})$ discussed  in ~\cite[Definition 3.10]{Ohk}.

Now, let us consider  a rank $n$ vector bundle
 $\mcF$ on $X$  equipped with 
a $\mcD^{(\M)}$-module structure $\nabla$.
Recall that the composite
\begin{align}
\psi (\nabla) : 
\mcT^{\otimes p^{m+1}}
 \xrightarrow{\psi_{X^\mr{log}/S^\mr{log}}} \mcD^{(m)} \xrightarrow{\nabla} \mcE nd_{\mcO_S} (\mcF)
\end{align}
is called the {\bf $p^{m+1}$-curvature (map)} of $\nabla$.
As proved in ~\cite[Lemma 4.15]{Ohk}, the image of $\psi (\nabla)$ actually  lies in $\mcE nd_{} (\mcF) := \mcE nd_{\mcO_X} (\mcF) \left(\subseteq \mcE nd_{\mcO_S} (\mcF) \right)$.
Also, if $a$ is a nonnegative integer with 
$a < m$, then we see that 
 the $p^{a+1}$-curvature of $\nabla^{(m)\Rightarrow (a)}$ vanishes.

We denote by 
\begin{align}
\mcY := [(\Omega^{\otimes p^{\M +1}})^\times \times^{\mbG_m} \mfg) / \mr{PGL}_n]
\end{align}
the quotient stack of $(\Omega^{\otimes p^{\M +1}})^\times \times^{\mbG_m} \mfg$ by the adjoint  $\mr{PGL}_n$-action.
Then,  
the pair $(\mcF, \psi (\nabla))$ specifies, via the quotient $\mr{Lie} (\mr{GL}_n) \twoheadrightarrow \mfg$,  a global section of this  stack over $X$.
The composite of this section with the natural quotient  
$\mcY \rightarrow (\Omega^{\otimes p^{\M +1}})^\times \times^{\mbG_m} \mfc$
 defines  an $S$-rational point 
\begin{align} \label{Er100}
\overline{\psi} (\nabla) \in B_\msX (S)
\end{align}
of 
$B_\msX$.
Under the identification \eqref{Er99},
this $S$-rational point corresponds to the coefficients $a_\ell \in H^0 (X, \Omega^{\otimes p^{m+1}\cdot \ell})$  ($2 \leq \ell \leq n$)
of 
  the characteristic polynomial
\begin{align} \label{EQ400}
\mr{Char} (\breve{\psi} (\nabla)) := \mr{det} (t - \breve{\psi} (\nabla)) = t^n +  a_{2}\cdot t^{n-2} + a_{3} \cdot t^{n-3} + \cdots + a_{n},
\end{align}
where $\breve{\psi} (\nabla)$ denotes the trace-free part of the $p^{\M +1}$-curvature $\psi (\nabla)$, i.e., 
$\breve{\psi} (\nabla) := \psi (\nabla) - \frac{1}{n} \cdot \mr{tr} (\psi (\nabla)) \otimes \mr{id}_\mcF$.

\SSP
\bpr \label{Prop8}
Let us keep the above notation, and we regard  $B_\msX^\flat (S)$ as a subset of $B_\msX (S)$ via $\iota_\msX$.
Then, the $S$-rational point
$\overline{\psi} (\nabla)$ lies in $B^\flat_\msX (S)$, or equivalently,  the element $a_\ell$ lies in the subspace $H^0 (X^{(\M +1)}, (\Omega^{(\M +1)})^{\otimes \ell})$ via the injection \eqref{Ey10}  for every $\ell = 2, \cdots, n$.
\epr
\begin{proof}
Denote by $\mcE nd (\mcF, \nabla)$ (resp., $\mcE nd (\mcF^{[m]}, \nabla^{[m]})$)
the subsheaf of $\mcE nd_{\mcO_S} (\mcF)$ (resp., $\mcE nd_{\mcO_S}(\mcF^{[m]})$)
consisting of endomorphisms preserving the $\mcD^{(\M)}$-module structure $\nabla$ (resp., the $S^\mr{log}$-connection $\nabla^{[m]}$).
The sheaf $\mcE nd (\mcF, \nabla)$ (resp., $\mcE nd (\mcF^{[m]}, \nabla^{[m]})$)
carries  a natural $\mcO_{X^{(m+1)}}$-module structure.
The $p^{m+1}$-curvature $\psi (\nabla) :  \mcT^{\otimes p^{\M +1}} \rightarrow \mcE nd_{\mcO_S} (\mcF)$ of $\nabla$
(resp., $p$-curvature $\psi (\nabla^{[m]}) : (\mcT^{(m+1)})^{\otimes p} \rightarrow \mcE nd_{\mcO_S} (\mcF^{[m]})$ of $\nabla^{[m]}$)
restricts to an $\mcO_{X^{(m+1)}}$-linear morphism
\begin{align}
\psi (\nabla)^\flat : \mcT^{(m+1)} \rightarrow \mcE nd (\mcF, \nabla) \  \left(\text{resp.,} \  \psi (\nabla^{[m]})^\flat  : \mcT^{(m+1)} \rightarrow \mcE nd (\mcF^{[m]}, \nabla^{[m]}) \right)
\end{align}
under the identification $F_{X/S}^{(m+1)*} (\mcT^{(m+1)}) \cong \mcT^{\otimes p^{\M +1}}$ (resp., $F^*_{X^{(\M)}/S} (\mcT^{(\M +1)}) \cong (\mcT^{(\M +1)})^{\otimes p}$).
By the definition of this  morphism,  the  following diagram commutes:
\begin{align} \label{eQ26}
\vcenter{\xymatrix@C=46pt@R=36pt{
& \mcT^{(m+1)} \ar[dl]_-{\psi (\nabla)^\flat} \ar[dr]^-{\psi (\nabla^{[m]})^\flat} &
\\
\mcE nd (\mcF, \nabla)\ar[rr]_-{\mr{res}_{(\mcF, \nabla)}}&& \mcE nd (\mcF^{[m]}, \nabla^{[m]}),
 }}
\end{align}
where the lower horizontal arrow $\mr{res}_{(\mcF, \nabla)}$ is  given by $h \mapsto h |_{\mcF^{[m]}}$
for any $h \in \mcE nd (\mcF, \nabla)$.

We here denote by $U$ the smooth locus  of $X \setminus \bigcup_{i=1}^r \mr{Im} (\sigma_i)$ relative to $S$, which is 
a scheme-theoretic dense open subscheme of $X$.
Then,  $\mr{res}_{(\mcF, \nabla)}$ is an isomorphism over  $U^{(\M +1)}$.
According to   ~\cite[Proposition 3.1]{ChZh0} (or  ~\cite[Proposition 3.2]{LasPa0}, ~\cite[Proposition-Definition 3.24]{Wak5}),
 the characteristic polynomial $\mr{Char} (\breve{\psi} (\nabla^{[m]})|_U)$ of the traceless $p$-curvature  $\breve{\psi} (\nabla^{[m]})|_U$ lies in $\bigoplus_{\ell =2}^n H^0 (U^{(\M +1)}, (\Omega^{(\M +1)})^{\otimes \ell})$.
The commutativity  of \eqref{eQ26}  implies 
that this characteristic polynomial 
coincides with 
 $\mr{Char} (\breve{\psi} (\nabla)|_U)$ via $\mr{res}_{(\mcF, \nabla)}$,
 and hence, we have 
 \begin{align} \label{Re220}
 \mr{Char} (\breve{\psi} (\nabla)|_U) \in \bigoplus_{\ell =2}^n H^0 (U^{(\M +1)}, (\Omega^{(\M +1)})^{\otimes \ell}).
 \end{align}

Furthermore,  by an argument similar to the proof of ~\cite[Proposition 3.3]{Wak5},
 the cokernel of $\omega_\ell$ for every $\ell$ is relatively torsion-free over $S$.
It follows from this fact together with \eqref{Re220} that 
 $\mr{Char} (\breve{\psi} (\nabla)$ lies in $\bigoplus_{\ell =2}^n H^0 (X^{(\M +1)}, (\Omega^{(\M +1)})^{\otimes \ell})$.
 This completes the proof of this assertion.
\end{proof}

\LSP
\subsection{Exponents of $\mcD^{(\M)}$-modules} \label{SS90}

In what follows, let us consider the situation where the underlying curve is a relative formal disk over $S$.
We shall set $U_\oslash := \mr{Spec} (\mcO_S [\![t]\!])$, where $t$ denotes  a formal parameter.
For simplicity, we write $\mcO_\oslash := \mcO_{U_\oslash}$.
We equip $U_\oslash$ with the log structure associated to
the marked point 
 $\sigma_\oslash : S \rightarrow U_\oslash$
 corresponding to the surjection $\mcO_S [\![t]\!] \twoheadrightarrow \mcO_S$ given by $t \mapsto 0$.
More explicitly, this log structure is given by 
 the monoid morphism $\mbN \rightarrow \mcO_{\oslash}$ given by $\ell \mapsto t^\ell$.
Let  $U_\oslash^\mr{log}$ denote the resulting log scheme.

 We then define  a sheaf of noncommutative rings on $U_\oslash$ by 
 $\mcD^{(\M)}_{\oslash}
   := \varprojlim_{\ell \geq 1} \mcD_{U_{\oslash, \ell}^\mr{log}/S}^{(\M)}$,
where $U_{\oslash, \ell}^\mr{log}$ ($\ell \geq  0$) denotes the strict closed subscheme of $U_\oslash^{\mr{log}}$ determined  by the ideal sheaf $(t^\ell) \subseteq \mcO_\oslash$.
This sheaf carries  two $\mcO_\oslash$-module structures  ${^L}\mcD_\oslash^{(\M)}$ and ${^R}\mcD_\oslash^{(\M)}$ arising from those of ${^L}\mcD_{U^\mr{log}_{\oslash, n}/k}^{(\M)}$'s and ${^R}\mcD_{U^\mr{log}_{\oslash, n}/k}^{(\M)}$'s (defined as in Section \ref{SS11}), respectively.
The $\mcO_\oslash$-module ${^L}\mcD_{\oslash}^{(\M)}$ admits a decomposition  as a direct sum $\bigoplus_{j \in \mbZ_{\geq 0}} \mcO_\oslash \cdot \partial^{\langle j \rangle}$, where $\partial^{\langle j \rangle}$'s are  the sections associated to the logarithmic coordinate $x$, as in ~\cite[Section 2.2]{Wak20}.

For each $d \in \mbZ/p^{\M +1}\mbZ$, 
we denote by $\widetilde{d}$ the integer defined as the unique lifting of $d$ via the natural surjection $\mbZ \twoheadrightarrow \mbZ/p^{\M +1} \mbZ$ satisfying $0 \leq \widetilde{d} < p^{\M +1}$.
There exists a unique $\mcD_\oslash^{(\M)}$-module structure
\begin{align}
\nabla_{\oslash, d}^{(\M)}  : {^L}\mcD_\oslash^{(\M)} \rightarrow \mcE nd_{\mcO_S} (\mcO_\oslash)
\end{align}
on $\mcO_\oslash$ determined by the condition that $\nabla_{\oslash, d}^{(\M)} (\partial^{\langle j \rangle}) (t^\ell) = q_j ! \cdot \binom{\ell - \widetilde{d}}{j} \cdot t^\ell$ for every $j, \ell \in \mbZ_{\geq 0}$.
The resulting $\mcD_\oslash^{(\M)}$-module
\begin{align} \label{UU3}
\msO_{\oslash, d}^{(\M)} := (\mcO_\oslash, \nabla_{\oslash, d}^{(\M)})
\end{align}
is isomorphic to the unique extension of $\msO_{\oslash, 0}^{(\M)}$ to $t^{-\widetilde{d}} \cdot \mcO_\oslash \left(\supseteq \mcO_\oslash \right)$.

The following lemma will be applied in the discussion of Section \ref{SS19e9} several times.

\SSP
\ble \label{Lem4490}
Let $d$ be an element of  $\mbZ/p^{\M +1} \mbZ$.
Then,  the composite 
\begin{align}
\mcS ol (\nabla_{\oslash, d}^{(\M)}) \hookrightarrow \mcO_\oslash \twoheadrightarrow \sigma_{\oslash*} (\sigma^*_\oslash (\mcO_\oslash))
\end{align}
is surjective if and only if $d = 0$.
\ele
\begin{proof}
The assertion follows from ~\cite[Proposition 4.4]{Wak20}.
\end{proof}
\SSP

Now, let us return  to the global setting (hence $\mcD^{(\M)} := \mcD_{X^\mr{log}/S^\mr{log}}^{(\M)}$ for $\msX := (X/S, \{ \sigma_i \}_i)$ as before).
Suppose that $S$ is connected, and 
let $(\mcF, \nabla)$ be a 
$\mcD^{(\M)}$-module with vanishing $p^{\M +1}$-curvature such that $\mcF$ is a vector bundle of rank $n$.
Also, let $i \in \{1, \cdots, r \}$.
According to  ~\cite[Proposition-Definition 4.8]{Wak20}, there exists 
a unique multiset
\begin{align} \label{Eq1114}
e_i (\nabla) := [d_1, \cdots, d_n]
\end{align}
with cardinality $n$ 
satisfying the following condition:
for each $k$-rational  point $q$ of $S$,
there exists, after base-change  along $\mr{Spec}(\mcO_{S, q}) \rightarrow S$,  an isomorphism of $\mcD_\oslash^{(\M)}$-modules  $\bigoplus_{i=1}^n \msO^{(\M)}_{\oslash, d_i} \xrightarrow{\sim} (\mcF, \nabla)|_{U_\oslash}$   under 
any identification of $U_\oslash^\mr{log}$ with the formal neighborhood of $\mr{Im} (\sigma_i) \left(\subseteq X \right)$.
We  refer to $e_i (\nabla)$ as the {\bf (local) exponent} of $\nabla$ at $\sigma_i$.

\LSP
\subsection{Monodromy operators} \label{SS1076}

We recall the notion of monodromy operator introduced in ~\cite[Section 4.2]{Wak20}.
Let us  set
$\mcB_S := \bigoplus_{j \in \mbZ_{\geq 0}} \mcO_S \cdot \partial_\mcB^{\langle j \rangle}$,
where $\partial_\mcB^{\langle i \rangle}$'s are   abstract symbols.
We equip $\mcB_S$ with a structure of $\mcO_S$-algebra by declaring  
\begin{align} \label{YY12}
\partial_\mcB^{\langle j' \rangle} \cdot \partial_\mcB^{\langle j''\rangle} = \sum_{j = \mr{max}\{j', j'' \}}^{j' + j''} \frac{j!}{(j' + j''-j)! \cdot (j-j')! \cdot (j-j'')!}\cdot \frac{q_{j'}! \cdot q_{j''}!}{q_j !} \cdot \partial_\mcB^{\langle j\rangle}.
\end{align}
In particular,  $\mcB_S$ is a commutative $\mcO_S$-algebra  generated by the sections $\partial_\mcB^{\langle 1 \rangle}, \partial_\mcB^{\langle p \rangle}, \cdots, \partial_\mcB^{\langle p^{\M}\rangle}$.

Suppose that we are given 
 an $\mcO_S$-module $\mcG$ equipped with 
  an $\mcO_S$-algebra morphism  
$\mu : \mcB_S \rightarrow \mcE nd_{\mcO_S} (\mcG)$.
For each $a\in \mbZ_{\geq 0}$,
we set 
$\mu^{\langle p^a \rangle} := \mu  (\partial_\mcB^{\langle p^a \rangle})\in \mr{End}_{\mcO_S}(\mcG)$.
The morphism $\mu$ is completely  determined by  
the $(\M+1)$-tuple
\begin{align} \label{e187}
\mu^{\langle \bullet \rangle} := (\mu^{\langle 1 \rangle}, \mu^{\langle p \rangle} \cdots,  \mu^{\langle p^{\M} \rangle}) \in \mr{End}_{\mcO_S}(\mcG)^{\oplus (\M +1)}.
\end{align}

Let $U_\oslash$, $\mcD_\oslash^{(\M)}$, and $\sigma_\oslash$ be as before.
It follows from \eqref{e364} and  \eqref{YY12} that the assignment $\partial^{\langle j \rangle} \mapsto \partial_\mcB^{\langle j \rangle}$ ($j \in \mbZ_{\geq 0}$) determines an isomorphism of $\mcO_S$-algebras 
\begin{align} \label{Eq48}
\sigma_\oslash^*({^L}\mcD_\oslash^{(\M)}) \xrightarrow{\sim} \mcB_S.
\end{align}

Now, let us take a $\mcD_\oslash^{(\M)}$-module $(\mcF, \nabla)$.
The  $\mcD_\oslash^{(\M)}$-module structure  $\nabla$ induces  a $\sigma_\oslash^*({^L}\mcD_\oslash^{(\M)})$-action $\sigma_\oslash^*(\nabla)$  on $\sigma_\oslash^*(\mcF)$.
This yields a composite
\begin{align}
\mu (\nabla^{(\M)}) : \mcB_S \xrightarrow{\eqref{Eq48}^{-1}} \sigma_\oslash^*({^L}\mcD_\oslash^{(\M)}) \xrightarrow{\sigma^*_\oslash (\nabla)} \mcE nd_{\mcO_S} (\sigma_\oslash^* (\mcF)).
\end{align}
In particular, we obtain  a tuple 
\begin{align} \label{e18d7}
\mu (\nabla)^{\langle \bullet \rangle} := (\mu (\nabla)^{\langle 1 \rangle}, \mu (\nabla)^{\langle p \rangle} \cdots,  \mu (\nabla)^{\langle p^{\M} \rangle}) \in \mr{End}_{\mcO_S}(\sigma_\oslash^*(\mcF))^{\oplus (\M +1)}.
\end{align} 
Since $\mcB_S$ is commutative,
these endomorphisms 
  commute with one another.
We refer to $\mu (\nabla)$ (and the associated tuple $\mu (\nabla)^{\langle \bullet \rangle}$) as 
the {\bf monodromy operator} of $\nabla$, following the terminology of  ~\cite[Definition 4.2]{Wak20}.
 Moreover, 
 for each $a = 0, \cdots, \M$,
we call  $\mu (\nabla)^{\langle p^a \rangle}$  the {\bf $a$-th monodromy operator} of $\nabla$.
The relationship between the notions of exponent and  monodromy operator is described in ~\cite[Proposition-Definition 4.8]{Wak20}.

Next, we shall write $\mr{End}_{\mcO_S}^0 (\sigma_i^* (\mcF))$ for  the subspace of $\mr{End}_{\mcO_S} (\sigma_i^* (\mcF))$ consisting of endomorphisms with vanishing trace.
We shall 
write
\begin{align}
\breve{\mu} (\nabla)^{\langle \bullet \rangle} := (\breve{\mu} (\nabla)^{\langle 1 \rangle}, \breve{\mu} (\nabla)^{\langle p \rangle}, \cdots, \breve{\mu} (\nabla)^{\langle p^{m} \rangle}) \in (\mr{End}_{\mcO_S}^0 (\sigma_\oslash^* (\mcF)))^{\oplus (m+1)},
\end{align}
where 
$\breve{\mu} (\nabla)^{\langle p^a \rangle} := \mu (\nabla)^{\langle p^a\rangle} - \frac{1}{n} \cdot \mr{tr} (\mu (\nabla)^{\langle p^a\rangle}) \cdot \mr{id}_{\sigma_i^*(\mcF)}$.

Next, 
let us  return  to the global setting.
Let $\msX$ be as before and fix $i \in \{1, \cdots, r \}$.
Also,
 let us take a $\mcD^{(m)}$-module $(\mcF, \nabla)$.
Assume that 
the image $\mr{Im}(\sigma_i)$ is determined by ``$t = 0$" for  some  function $t$  on its  open neighborhood.
This assumption allows us to identify  $U_\oslash$ with the formal neighborhood of $\mr{Im}(\sigma_i)$.
In particular, we obtain $\mcD_{\oslash}^{(\N -1)}$-module $(\mcF |_{U_\oslash}, \nabla |_{U_\oslash})$, which determines  
elements
 \begin{align} \label{Eq58}
 \mu_i (\nabla)^{\langle \bullet \rangle} := 
 \mu (\nabla |_{U_\oslash})^{\langle \bullet \rangle},
\hspace{10mm}
\breve{\mu}_i (\nabla)^{\langle \bullet \rangle} := 
\breve{\mu} (\nabla |_{U_\oslash})^{\langle \bullet \rangle} 
 \end{align}
 of $\mr{End}_{\mcO_S} (\sigma_i^* (\mcF))^{\oplus (m+1)}$
 under the natural  identification $\sigma_i^* (\mcF) = \sigma_\oslash^* (\mcF |_{U_\oslash})$.

\vspace{10mm}
\section{The moduli space of dormant $\mr{PGL}_n^{(\N)}$-opers} \label{S5}
\LSP

This section is devoted to the study of  dormant $\mr{PGL}_n^{(\N)}$-opers and   their moduli space; for a detailed treatment of the subject, we refer the reader to   ~\cite{Wak20}.
We  also  introduce a generalization of the Hitchin-Mochizuki morphism, defined via  higher-level $p$-curvature.
As an application, we establish several geometric properties of this morphism (cf. Theorems \ref{Cor20}, \ref{Th665}).

\LSP
\subsection{Dormant/$p^\N$-nilpotent $\mr{PGL}_n^{(\N)}$-opers} \label{SS100}

Let $(g, r)$ be a pair of nonnegative integers with $2g-2+r > 0$,
$S$   a scheme (or more generally, a Deligne-Mumford stack) over $k$,  and
$\msX := (f: X \rightarrow S, \{ \sigma_i \}_{i=1}^r)$ an $r$-pointed stable curve  of genus $g$ over $S$, as in the previous section.
Also, we set $\N := \M +1$, i.e., $\N$ is a positive integer.

Consider a pair $\msE^\spadesuit := (\mcE_B, \phi)$ consisting of 
a $B$-bundle  $\mcE_B$ on $X$ and an $(\N -1)$-PD stratification on $\mcE := \mcE_B \times^B \mr{PGL}_n$.
 (The precise definition of a  higher-level PD stratification  on a $G$-bundle will not be recalled here, as 
the discussion below will only concern its translation into  a $\mcD^{(\N-1)}$-actions on a vector bundle.
For full details, see   ~\cite[Definition 2.3]{Wak20}.)
Denote by $\nabla_\phi$  the $S^\mr{log}$-connection on $\mcE$ corresponding to the $0$-PD stratification induced by  $\phi$.

\SSP
\bde[cf. ~\cite{Wak20}, Definitions 5.2, 5.3] \label{Def3}
\begin{itemize}
\item[(i)]
We say that $\msE^\spadesuit$ is a {\bf $\mr{PGL}_n^{(\N)}$-oper}  (or a {\bf $\mr{PGL}_n$-oper of level $\N$})
on $\msX$ if 
the pair $(\mcE_B, \nabla_\phi)$ forms a $\mr{PGL}_n$-oper (equivalently,  an $\mfs \mfl_n$-oper) on $\msX$, in the sense of ~\cite[Definition 2.1]{Wak5}.
One can define the notion of an isomorphism between two   $\mr{PGL}_n^{(\N)}$-opers.
\item[(ii)]
 A $\mr{PGL}_n$-oper $\msE^\spadesuit := (\mcE_B, \phi)$ is said to be {\bf dormant} if
  the $p^\N$-curvature of $\phi$ (in the sense of ~\cite[Definition 2.10, (ii)]{Wak20}) vanishes identically.
\end{itemize}
\ede
\SSP

To describe 
 $\mr{PGL}_n^{(\N)}$-opers in terms of $\mcD^{(\N -1)}$-modules, we  recall the notion of an $n^{(\N)}$-theta characteristic introduced in ~\cite{Wak20}, as follows.

\SSP
\bde[cf. ~\cite{Wak20}, Definition 5.12] \label{Def6}
An {\bf $n^{(N)}$-theta characteristic} of $\msX$ is a pair
\begin{align} \label{Eq7}
\vartheta := (\varTheta, \nabla_\vartheta)
\end{align}
consisting of a line bundle $\varTheta$ on $X$ and a $\mcD^{(N -1)}$-module structure $\nabla_\vartheta$ on the line bundle $\mcT^{\otimes \frac{n (n-1)}{2}} \otimes \varTheta^{\otimes n}$.
We say that  an $n^{(N)}$-theta characteristic $\vartheta := (\varTheta, \nabla_\vartheta)$ is   {\bf dormant} if 
$\nabla_\vartheta$ has vanishing $p^\N$-curvature.
  \ede
\SSP

Let us take a dormant $n^{(N)}$-theta characteristic $\vartheta := (\varTheta, \nabla_\vartheta)$ of $\msX$.
(According to ~\cite[Proposition 5.14]{Wak20}, such an object always exists.)
We set
\begin{align}
\mcF_\varTheta := \mcD^{(\N -1)}_{\leq n-1} \otimes \varTheta, \hspace{5mm} \mcF_\varTheta^j := \mcD^{(\N -1)}_{\leq n  -j-1} \otimes \varTheta
 \ (j=0, \cdots, n).
\end{align}
Note   that these vector bundles do not depend on the level $\N$ because of the assumption $n < p$.
Since $\mcF_\varTheta^{j}/\mcF_\varTheta^{j+1}$ ($j=0, \cdots, n-1$)
can be identified with $\mcT^{\otimes (n-j-1)} \otimes \varTheta$, we obtain a composite isomorphism
\begin{align}
\mr{det}(\mcF_\varTheta) \xrightarrow{\sim} \bigotimes_{j=0}^{n-1} \mcF_\varTheta^j/\mcF_\varTheta^{j+1} \xrightarrow{\sim}
\bigotimes_{j=0}^{n-1} (\mcT^{\otimes (n-j-1)} \otimes \varTheta)
\xrightarrow{\sim} \mcT^{\otimes \frac{n (n-1)}{2}} \otimes \varTheta^{\otimes n}.
\end{align}

\SSP
\bde[cf. ~\cite{Wak20}, Definition 5.15] \label{Def8}
\begin{itemize}
\item[(i)]
A {\bf  $(\mr{GL}_n^{(N)}, \vartheta)$-oper} on $\msX$ is a $\mcD^{(N-1)}$-module structure $\nabla^\diamondsuit$ on $\mcF_\varTheta$   such that, for every $j=0, \cdots, n-1$,
the $\mcO_X$-linear morphism $\mcD^{(\N -1)} \otimes \mcF_\varTheta \rightarrow \mcF_\varTheta$ induced by $\nabla^\diamondsuit$ restricts to an isomorphism
\begin{align} \label{EQ200}
\mcD_{\leq n-j-1}^{(\N -1)} \otimes \mcF_\varTheta^{n-1} \xrightarrow{\sim} \mcF_\varTheta^j.
\end{align}
We say that  a  $(\mr{GL}_n^{(\N)}, \vartheta)$-oper $\nabla^\diamondsuit$ is  {\bf dormant} if its $p^\N$-curvature vanishes identically.
\item[(ii)]
Let $\nabla^\diamondsuit_\circ$ and $\nabla^\diamondsuit_\bullet$ be (dormant) $(\mr{GL}_n^{(N)}, \vartheta)$-opers on $\msX$.
An {\bf isomorphism of (dormant) $(\mr{GL}_n^{(\N)}, \vartheta)$-opers} from  $\nabla^\diamondsuit_\circ$ to $\nabla^\diamondsuit_\bullet$
is an isomorphism of $\mcD^{(\N -1)}$-modules $(\mcF_\varTheta, \nabla^\diamondsuit_\circ) \xrightarrow{\sim} (\mcF_\varTheta, \nabla^\diamondsuit_\bullet)$ preserving the filtration $\{ \mcF_\varTheta^j \}_{j=0}^n$.
\end{itemize}
  \ede
\SSP

 Let $\nabla^\diamondsuit$ be a   $(\mr{GL}_n^{(\N)}, \vartheta)$-oper on $\msX$.
 Denote by $\mcE$  the $\mr{PGL}_n$-bundle obtained by projectivizing $\mcF_\varTheta$.
  Then,  $\nabla^\diamondsuit$ induces an $(\N -1)$-PD stratification $\phi_{\nabla^\diamondsuit}$ on $\mcE$.
 The filtration $\{ \mcF_\varTheta^j \}_{j=0}^n$
  defines  a $B$-reduction $\mcE_B$ of $\mcE$, and
  the pair
 \begin{align}
 \nabla^{\diamondsuit \Rightarrow \spadesuit} := (\mcE_B, \phi_{\nabla^\diamondsuit})
 \end{align}
 constitutes  a $\mr{PGL}_n^{(\N)}$-oper, since 
 the morphisms \eqref{EQ200} associated to  the  various $j$'s  are isomorphisms.
 The resulting assignment $\nabla^\diamondsuit \mapsto \nabla^{\diamondsuit \Rightarrow \spadesuit}$ determines a well-defined bijection of sets 
 \begin{align} \label{Eq107v}
\left(\begin{matrix}\text{the set of isomorphism classes of} \\ \text{$(\mr{GL}_n^{(\N)}, \vartheta)$-opers on $\msX$} \end{matrix} \right)
\xrightarrow{\sim}
\left(\begin{matrix}\text{the set of isomorphism classes of} \\ \text{$\mr{PGL}_n$-opers on $\msX$} \end{matrix} \right)
\end{align}
(cf. ~\cite[Theorem 5.18]{Wak20}).

Moreover, since $\nabla_\vartheta$ is dormant, we see that   $\nabla^\diamondsuit$ is dormant if and only if the corresponding $\mr{PGL}_2^{(\N)}$-oper $ \nabla^{\diamondsuit \Rightarrow \spadesuit}$ is dormant.
It follows that 
\eqref{Eq107v} restricts to a bijection
\begin{align} \label{Eq107}
\left(\begin{matrix}\text{the set of isomorphism classes of} \\ \text{dormant $(\mr{GL}_2^{(\N)}, \vartheta)$-opers on $\msX$} \end{matrix} \right)
\xrightarrow{\sim}
\left(\begin{matrix}\text{the set of isomorphism classes of} \\ \text{dormant $\mr{PGL}_n^{(\N)}$-opers on $\msX$} \end{matrix} \right).
\end{align}

\SSP
\bde \label{Def44}
\begin{itemize}
\item[(i)]
A $(\mr{GL}_n^{(\N)}, \vartheta)$-oper  $\nabla^\diamondsuit$  is said to be {\bf $p^\N$-nilpotent} if the equality $\mr{Char} (\breve{\psi}(\nabla^\diamondsuit)) = 0$ holds (cf. \eqref{EQ400}).
\item[(ii)]
For 
  a $\mr{PGL}_n^{(\N)}$-oper $\msE^\spadesuit$,
  we shall set
  \begin{align} \label{Er115}
  \psi\text{-}\mr{Char} (\msE^\spadesuit) := 
  \mr{Char} (\breve{\psi}(\nabla^\diamondsuit)),
  \end{align}
  where $\nabla^\diamondsuit$ denotes 
  the corresponding $(\mr{GL}_2^{(\N)}, \vartheta)$-oper via \eqref{Eq107v}.
  (Note that 
 this definition
   does not depend on the choice of a dormant $n^{(\N)}$-theta characteristic $\vartheta$.)
 We say that $\msE^\spadesuit$  is   {\bf $p^\N$-nilpotent} if
 $\psi\text{-}\mr{Char} (\msE^\spadesuit) = 0$, i.e., 
  the corresponding $(\mr{GL}_2^{(\N)}, \vartheta)$-oper  $\nabla^\diamondsuit$ is $p^\N$-nilpotent.
\end{itemize}
\ede

\LSP
\subsection{Radii of $\mr{PGL}_n^{(\N)}$-opers} \label{SS43}

 Let us consider the direct product  $\mfg^{\N} := \mfg \times \cdots \times \mfg$ of $\N$ copies of $\mfg$.
The adjoint $\mr{PGL}_n$-action on $\mfg$ induces a diagonal  $\mr{PGL}_n$-action on this product.
Accordingly, we define  its  GIT quotient by  
\begin{align} \label{Eq449}
\mfc_{\N} := \mfg^{\N\ma}/\!/\mr{PGL}_n,
\end{align}
which is a $k$-scheme and coincides with $\mfc$ when $\N =1$.
If   $[\mfg^{\N\ma}/\mr{PGL}_n]$ denotes the quotient stack of $\mfg^{\N\ma}$ by the same $\mr{PGL}_n$-action,
then there exists a natural morphism $[\mfg^{\N\ma}/\mr{PGL}_n] \rightarrow \mfc_\N$.

Given  a positive integer $\N'$ with $\N' \leq \N$ and an $r$-tuple $\rho := (\rho_i)_{i=1}^r \in \mfc_\N (S)^{r\ma}$, 
we set 
\begin{align} \label{Ey78}
\rho^{(\N')} := (\rho_i^{(\N')})_{i=1}^r,
\end{align}
 where $\rho_i^{(\N')}$ ($i=1, \cdots, r$) denotes the image of $\rho_i$ via $\mfc_\N \rightarrow \mfc_{\N'}$ induced from the projection onto the first $\N'$ factors $\mfg^{\N} \twoheadrightarrow \mfg^{\N'}$. 

Now, let $\nabla^\diamondsuit$ be a $(\mr{GL}_n^{(\N)}, \vartheta)$-oper on $\msX$.
For each $i=1, \cdots, r$,
the element   $\mu_i (\nabla^\diamondsuit)^{\langle \bullet \rangle}$ determines  an $S$-rational point of  $[\mfg^{\N}/\mr{PGL}_n]$ via $\mr{Lie}(\mr{GL}_n) \twoheadrightarrow \mfg$.
The image of this point under the morphism  $[\mfg^{\N}/\mr{PGL}_n] \rightarrow \mfc_\N$ defines an $S$-rational point
\begin{align} \label{Er6}
\rho_i (\nabla^\diamondsuit) \in \mfc_\N (S).
\end{align}

Moreover,  given a $\mr{PGL}_n$-oper $\msE^\spadesuit$ on $\msX$,
we shall set
\begin{align} \label{Er7}
\rho_i (\msE^\spadesuit) := \rho_i (\nabla^\diamondsuit_{\msE^\spadesuit}) \in \mfc_\N (S),
\end{align}
where $\nabla^\diamondsuit_{\msE^\spadesuit}$ is the $(\mr{GL}_n^{(\N)}, \vartheta)$-oper  corresponding to $\msE^\spadesuit$.
Note that $\rho_i (\msE^\spadesuit)$ does not depend on the choice of $\vartheta$.

\SSP
\bde \label{Def66}
\begin{itemize}
\item[(i)]
We refer to $\rho_i (\nabla^\diamondsuit)$ (resp., $\rho_i (\msE^\spadesuit)$)
as the {\bf radius} of $\nabla^\diamondsuit$ (resp., $\msE^\spadesuit$).
\item[(ii)]
Let 
 $\rho := (\rho_i)_{i=1}^r$ be an element of $\mfc_\N (S)^{r\ma}$.
We say that a $(\mr{GL}_n^{(\N)}, \vartheta)$-oper $\nabla^\diamondsuit$  (resp.,  a $\mr{PGL}_n^{(\N)}$-oper $\msE^\spadesuit$) is {\bf of radii $\rho$} if $\rho_i = \rho_i (\nabla^\diamondsuit)$ (resp., $\rho_i = \rho_i (\msE^\spadesuit)$) for every $i=1, \cdots, r$.
When $r = 0$, any $(\mr{GL}_n^{(\N)}, \vartheta)$-oper (resp.,  $\mr{PGL}_2^{(\N)}$-oper)  is said to be  {\bf of radii $\emptyset$}.
\end{itemize}
\ede
\SSP

For an element $d$ of $\mbZ/p^\N \mbZ$,
let $\widetilde{d}_{[0]}, \cdots, \widetilde{d}_{[\N -1]}$ be the collection of integers uniquely determined by the condition that
$\widetilde{d} = \sum_{a=0}^{\N -1} p^a \cdot \widetilde{d}_{[a]}$ and $0 \leq \widetilde{d}_{[a]} < p$ ($a=0, \cdots, \N -1$).
Given each $a = 0, \cdots, \N -1$, we write $\widetilde{d}_{[0, a]} := \sum_{b=0}^a p^b \cdot \widetilde{d}_{[b]}$, i.e., the remainder obtained by dividing $\widetilde{d}$ by $p^{a+1}$, and write
$d_{[a]}$ (resp., $d_{[0, a]}$) for the image of $\widetilde{d}_{[a]}$ (resp., $\widetilde{d}_{[0, a]}$)
via the natural projection $\mbZ \twoheadrightarrow \mbF_p$
 (resp., $\mbZ \twoheadrightarrow \mbZ/p^{a+1} \mbZ$).

Denote by $\Delta$ the image of the diagonal embedding $\mbZ/p^\N \mbZ \hookrightarrow (\mbZ/p^\N \mbZ)^{n\ma}$, which is a group homomorphism.
In particular, this gives rise the quotient set $(\mbZ/p^\N \mbZ)^{n\ma} /\Delta$.
The set $(\mbZ/p^\N \mbZ)^{n\ma}$ is equipped with the action of the symmetric group  of $n$ letters  $\mfS_n$ by permutation;
this  induces a well-defined $\mfS_n$-action on $(\mbZ/p^\N \mbZ)^{n\ma}/\Delta$.
Hence, we obtain the  quotient sets
\begin{align} \label{Er1}
\mfS_n \backslash (\mbZ/p^\N \mbZ)^{n\ma},
\hspace{10mm}
\mfS_n \backslash (\mbZ/p^\N \mbZ)^{n\ma}/\Delta,
\end{align}
and there is a natural projection map 
$\mfS_n \backslash (\mbZ/p^\N \mbZ)^{n\ma} \twoheadrightarrow 
\mfS_n \backslash (\mbZ/p^\N \mbZ)^{n\ma}/\Delta$.
Each element of $ \mfS_n \backslash (\mbZ/p^\N \mbZ)^{n\ma}$ can be interpreted  as a multiset of elements of $\mbZ/p^\N \mbZ$  with  cardinality  $n$.

Let  $\widetilde{\Xi}_{n, \N}$ denote  the subset of $\mfS_n \backslash (\mbZ/p^\N \mbZ)^{n\ma}$ consisting of multisets $[d_1, \cdots, d_n]$ such that
the induced elements $d_{1 [0]}, \cdots, d_{n [0]}$ of $\mbF_p$ 
are pairwise distinct.
(In particular, $\widetilde{\Xi}_{n, \N}$ can be identified with the  set of $n$-element  subsets of $\mbZ/p^\N \mbZ$ whose images in $\mbF_p$ are distinct.)
Then, we obtain the set 
\begin{align} \label{Er3}
\Xi_{n, \N} \ \left(\subseteq \mfS_n \backslash (\mbZ/p^\N \mbZ)^{n\ma}/\Delta \right)
\end{align}
consisting of all elements  admitting representatives in  $\widetilde{\Xi}_{n, \N}$.
We occasionally regard $\Xi_{n, \N}$ as the scheme defined to be  the disjoint union of copies of $\mr{Spec} (k)$ indexed by the elements in $\Xi_{n, \N}$. 

\SSP
\begin{exa} \label{Exam009}
Let us consider the case  $n=2$.
Denote by 
$(\mbZ/p^\N \mbZ)/\{ \pm 1 \}$ the set of equivalence classes of elements $a \in \mbZ/p^\N \mbZ$, in which $a$ and $-a$ are identified.
Then, the assignment 
$[a, b] \mapsto \overline{\big(\frac{a-b}{2}\big)}$
determines a bijection
\begin{align} \label{Eq311}
\mfS_2 \backslash (\mbZ/p^\N \mbZ)^{2\ma} /\Delta  \xrightarrow{\sim} (\mbZ/p^\N \mbZ)/\{ \pm 1 \}.
\end{align}
Under 
this bijection,
$\Xi_{2, \N}$ corresponds to the subset $(\mbZ/p^\N \mbZ)^{\times} / \{ \pm 1 \}$ of $(\mbZ/p^\N \mbZ)/ \{ \pm 1 \}$, i.e.,  the image of 
$(\mbZ/p^\N \mbZ)^\times$ under the natural projection  $(\mbZ/p^\N \mbZ) \twoheadrightarrow (\mbZ/p^\N \mbZ)/\{ \pm 1 \}$.
\end{exa}
\SSP

Let $\vec{d} := (d_1, \cdots, d_n)$ be an element of $(\mbZ/p^\N \mbZ)^{n\ma}$.
For each $j = 0, \cdots, \N -1$, denote by $D(\vec{d})_{[j]}$ the diagonal $n \times n$ matrix 
defined as 
\begin{align}
D(\vec{d})_{[j]}  := \begin{pmatrix}
 (-d_1)_{[j]} & 0 & \cdots & 0 & 0   \\
0 &(-d_2)_{[j]} & \cdots & 0 & 0  \\
\vdots & \vdots & \ddots & \vdots & \vdots  \\
0 & 0 & \cdots &  (-d_{n-1})_{[j]} & 0 \\ 
0 & 0 & \cdots & 0 &   (-d_n)_{[j]}
\end{pmatrix}
\end{align}
Also, denote by $\overline{D}(\vec{d})$
the element of $\mfc_\N (S)$ defined to be the image of $(D(\vec{d})_{[0]}, \cdots, D(\vec{d})_{[\N -1]})$ via the quotient $\mfg^{\N\ma} \twoheadrightarrow \mfc_\N$.
The map $(\mbZ/p^\N \mbZ)^{n\ma} \rightarrow \mfc_\N$ given by $\vec{d} \mapsto \overline{D}(\vec{d})$ factors through
the quotient $(\mbZ/p^\N \mbZ)^{n\ma} \twoheadrightarrow \mfS_n \backslash (\mbZ/p^\N \mbZ)^{n\ma}/\Delta$.
The resulting morphism 
\begin{align} \label{Er4}
\mfS_n \backslash (\mbZ/p^\N \mbZ)^{n\ma}/\Delta \rightarrow \mfc_\N (S)
\end{align}
is injective.
Using this injection, we may  regard  $\mfS_n \backslash (\mbZ/p^\N \mbZ)^{n\ma}/\Delta$ (and hence $\Xi_{n, \N}$) as a subset of $\mfc_\N (S)$.

The following assertion for dormant opers  was already established in ~\cite[Proposition 6.14]{Wak20}.
Indeed, the exponent  of  a dormant $(\mr{GL}_n^{(\\N)}, \vartheta)$-oper at each marked point is mapped to its radius  via $\mfS_n \backslash (\mbZ/p^\N \mbZ)^{n\ma} \twoheadrightarrow \mfS_n \backslash (\mbZ/p^\N \mbZ)^{n\ma}/\Delta$.

\SSP
\bpr \label{Prop94}
Let $\nabla^\diamondsuit$ (resp., $\msE^\spadesuit$) be a $p^\N$-nilpotent $(\mr{GL}_n^{(\N)}, \vartheta)$-oper (resp.,  a $p^\N$-nilpotent $\mr{PGL}_n^{(\N)}$-oper) on $\msX$ with  $\N > 1$.
Then, the radius of  $\nabla^\diamondsuit$ (resp., $\msE^\spadesuit$) at $\sigma_i$ lies in $\Xi_{n, \N}$.
\epr
\begin{proof}
It suffices to prove the non-resp'd assertion because of the definition of $p^\N$-nilpotency for $\mr{PGL}_n^{(\N)}$-opers.
To this end,
we may assume, without loss of generality, that $S = \mr{Spec}(R)$ for a local $k$-algebra  $(R, \mfm)$  with $R/\mfm = k$. 

Denote by 
$\vartheta'$ the $n^{(\N -1)}$-theta characteristic obtained from $\vartheta$ via level reduction, and by 
 $\nabla'$ the $(\mr{GL}_n^{(\N -1)}, \vartheta')$-oper obtained from $\nabla^\diamondsuit$.
Since $\nabla'$ is dormant, 
we may associate to it  the exponent  $[d_1, \cdots, d_n]$ at $\sigma_i$, where $d_1, \cdots, d_n \in \mbZ/p^{\N -1}\mbZ$.
Fix  a suitable trivialization $\sigma_i^* (\mcF_\varTheta) \cong \mcO_S^{\oplus n}$, under which we may arrange that 
\begin{align}
\breve{\mu}_i (\nabla^\diamondsuit)^{\langle p^a \rangle}  \left( =  \breve{\mu}_i (\nabla')^{\langle p^a \rangle}  \right)
= \begin{pmatrix}
 (-d_1)_{[a]} & 0 & \cdots & 0 & 0   \\
0 &(-d_2)_{[a]} & \cdots & 0 & 0  \\
\vdots & \vdots & \ddots & \vdots & \vdots  \\
0 & 0 & \cdots &  (-d_{n-1})_{[a]} & 0 \\ 
0 & 0 & \cdots & 0 &   (-d_n)_{[a]}
\end{pmatrix}
\end{align}
for each $a =0, \cdots, \N -2$.
Moreover, the dormancy condition  on $\nabla'$ implies that
the diagonal entries $(-d_1)_{[0]}, \cdots, (-d_{n})_{[0]}$ of 
$\breve{\mu}_i (\nabla^\diamondsuit)^{\langle 1 \rangle}$  are pairwise distinct (cf. ~\cite[Proposition 6.3, (ii)]{Wak20}).
As   $\breve{\mu}_i (\nabla^\diamondsuit)^{\langle p^{\N -1} \rangle}$ commutes with $\breve{\mu}_i (\nabla^\diamondsuit)^{\langle 1 \rangle}$,
 the matrix $\breve{\mu}_i (\nabla^\diamondsuit)^{\langle p^{\N -1} \rangle}$ must be diagonal, i.e., expressed as 
 \begin{align}
 \breve{\mu}_i (\nabla^\diamondsuit)^{\langle p^{\N -1} \rangle}
 =  \begin{pmatrix}
 h_1 & 0 & \cdots & 0 & 0   \\
0 &h_2 & \cdots & 0 & 0  \\
\vdots & \vdots & \ddots & \vdots & \vdots  \\
0 & 0 & \cdots &  h_{n-1}& 0 \\ 
0 & 0 & \cdots & 0 &   h_n
\end{pmatrix}
 \end{align}
 for some $h_1, \cdots, h_n \in H^0 (S, \mcO_S)$.

   The assumption that $\nabla^\diamondsuit$ is $p^\N$-nilpotent  implies
   \begin{align}
   ( \breve{\mu}_i (\nabla^\diamondsuit)^{\langle p^{\N -1} \rangle})^p -  \breve{\mu}_i (\nabla^\diamondsuit)^{\langle p^{\N -1} \rangle} = 0.
   \end{align}
   It follows that $h_j^p - h_j = 0$ for every $j=1, \cdots, n$, i.e., the elements $h_1, \cdots, h_n$ lie in $\mbF_p$.
Therefore, for each $j =1, \cdots, n$, there exists a lifting $d^+_j \in \mbZ/p^\N \mbZ$ of $d_j$ satisfying 
   $(-d^+_j)_{[\N -1]} = h_j$.
The multiset  $[d^+_1, \cdots, d^+_n]$ belongs to $\Xi_{n, \N}$, which coincides with the radius of $\nabla^\diamondsuit$ at $\sigma_i$.
This completes the proof of the assertion.
\end{proof}

\LSP
\subsection{The moduli spaces of dormant/$p^\N$-nilpotent $\mr{PGL}_n^{(\N)}$-opers} \label{SS83}

Let  $\overline{\mcM}_{g, r}$ denote 
the moduli stack classifying $r$-pointed stable curves of genus $g$ over $k$.
Also, let us fix an $r$-tuple  $\rho := (\rho_i)_{i=1}^r \in \mfc_\N (k)^{r\ma}$, where we set $\rho := \emptyset$ if $r =0$.
We define 
\begin{align}\label{Eq114}
\mcO p_{\N,  g, r} \ \left(\text{resp.,} \ \mcO p_{\N, \rho, g, r}; \text{resp.,} \ 
\mcO p_{\N, \rho, g, r}^{^\mr{Zzz...}}; \text{resp.,} \ \mcO p_{\N, \rho, g, r}^{^\mr{nilp}} \right)
\end{align}
to be the category that classifies pairs $(\msX, \msE^\spadesuit)$, where $\msX$ denotes  a pointed curve $\msX$ in $\mr{ob}(\overline{\mcM}_{g, r})$ and $\msE^\spadesuit$ denotes   a $\mr{PGL}_n^{(\N)}$-oper (resp., a $\mr{PGL}_n^{(\N)}$-oper of radii $\rho$; resp.,  a dormant $\mr{PGL}_n^{(\N)}$-oper of radii $\rho$; resp., a $p^\N$-nilpotent $\mr{PGL}_n^{(\N)}$-oper of radii $\rho$) on $\msX$.

Forgetting the data of $\mr{PGL}_n^{(N)}$-opers yields a projection
\begin{align} \label{EQ202}
\Pi_{g, r} : \mcO p_{\N,  \rho, g, r} \rightarrow \overline{\mcM}_{g, r}
\ \Big(\text{resp.,} \ \Pi_{\rho, g, r} : \mcO p_{\N,  \rho, g, r} \rightarrow \overline{\mcM}_{g, r};  \hspace{3mm}\\
 \text{resp.,} \  \Pi'_{\rho, g, r} : \mcO p^{^\mr{Zzz...}}_{\N,  \rho, g, r} \rightarrow \overline{\mcM}_{g, r};
\text{resp.,} \ \Pi''_{\rho, g, r} : \mcO p^{^\mr{nilp}}_{\N,  \rho, g, r} \rightarrow \overline{\mcM}_{g, r}\Big).  \notag
\end{align}
Since any $\mr{PGL}_n^{(\N)}$-oper does not have nontrivial automorphisms (cf. ~\cite[Proposition 5.5]{Wak20}),
the stack $\mcO p_{\N, g, r}$ (resp.,  $\mcO p_{\N, \rho, g, r}$; resp., $\mcO p_{\N, \rho, g, r}^{^\mr{Zzz...}}$; resp., $\mcO p_{\N, \rho, g, r}^{^\mr{nilp}}$) is  fibered in equivalence relations over $\overline{\mcM}_{g, r}$ (in the sense of ~\cite[Definition 3.39]{FGA}), i.e.,
specifies a set-valued \'{e}tale sheaf on  $\overline{\mcM}_{g, r}$.

The following assertion was previously 
established in  ~\cite{Moc2} for $(n, \N) =(2, 1)$, ~\cite{Wak5} for general $n$ with $\N =1$, and ~\cite{Wak20} for arbitrary $(n, \N)$.

\SSP
\bt \label{Th8}
\begin{itemize}
\item[(i)]
The category   $\mcO p_{\N, g, r}$ (resp.,  $\mcO p_{\N, \rho, g, r}$) can be   represented by  
a nonempty  (resp., a possibly empty)
 Deligne-Mumford stacks of finite type over $k$, and the  projection $\Pi_{g, r}$ (resp., 
$\Pi_{\rho, g, r}$) is represented by schemes and is affine.
Moreover, when $\N =1$, the stack $\mcO p_{\N, \rho, g, r}$  is  nonempty.
\item[(ii)]
The category $\mcO p^{^\mr{Zzz...}}_{\N, \rho,  g, r}$ can  be represented by a (possibly empty) proper Deligne-Mumford stack  over $k$,   and the projection $\Pi'_{\rho, g, r}$ is finite.
Moreover, $\mcO p^{^\mr{Zzz...}}_{\N, \rho, g, r}$ is empty unless $\rho$ lies in $\Xi_{n, \N}^{r\ma}$.
\item[(iii)]
Suppose that  $n=2$ and $\mcO p^{^\mr{Zzz...}}_{\N, \rho, g, r} \neq \emptyset$.
Then,  $\mcO p^{^\mr{Zzz...}}_{\N, \rho, g, r}$ is smooth over $k$ and  of dimension $3g-3 +r$.
Moreover,  $\Pi'_{\rho, g, r}$ is 
faithfully flat and generically \'{e}tale.
\end{itemize}
\et
\begin{proof}
The first assertion of  (i) follows from ~\cite[Theorem 5.22, Corollary 6.28, (i)]{Wak20}.
The second assertion follows from ~\cite[Theorem A]{Wak5}.
Assertion (ii) follows from  ~\cite[Proposition 6.14, Theorem 6.17]{Wak20}.
Also, assertion (iii) follows from ~\cite[Corollary 8.17, Theorem 8.29]{Wak20}.
\end{proof}
\SSP

Regarding the moduli stack of $p^\N$-nilpotent $\mr{PGL}_n^{(\N)}$-opers, 
we obtain  the following assertion, which  follows  immediately from 
the definition of $p^\N$-nilpotency.

\SSP
\bpr \label{PRop67}
The category $\mcO p_{\N, \rho, g, r}^{^\mr{nilp}}$ can be represented by a (possibly empty) closed substack  of $\mcO p_{\N, \rho, g, r}$, and the inclusion $\mcO p_{\N, \rho, g, r}^{^\mr{Zzz...}} \hookrightarrow \mcO p_{\N, \rho, g, r}^{^\mr{nilp}}$ is a closed immersion.
\epr
\SSP

We remark  that the nonemptiness of  $\mcO p_{\N, \rho, g, r}^{^\mr{Zzz...}}$ (as well as $\mcO p_{\N, \rho, g, r}^{^\mr{nilp}}$) depends on the choice of $\rho$.
In the case $n=2$,
a combinatorial analysis of $\mcO p_{\N, \rho, g, r}^{^\mr{Zzz...}}$ was carried out in 
 ~\cite{Wak20}  (or ~\cite{Moc2} for $\N =1$), 
 which provides a criterion for when this stack is nonempty.
 Applying this result, we obtain  the following assertion, which will be applied in the proof of Proposition \ref{Prop111}.

\SSP
\bpr \label{Prop77}
Suppose that $n=2$, and 
let $\rho := (\rho_i)_{i=1}^r$ be an element of 
$\Xi_{2, \N}^{r\ma}$
 with $\mcO p^{^\mr{Zzz...}}_{\N, \rho, g, r} \neq \emptyset$.
Then, there exists an $r$-tuple  $\widetilde{\rho} \in \mfc_{\N +1} (k)^{r\ma}$ with $\widetilde{\rho}^{(\N)} = \rho$ (cf. \eqref{Ey78}) and $\mcO p^{^\mr{Zzz...}}_{\N +1, \widetilde{\rho}, g, r} \neq \emptyset$.
\epr
\begin{proof}
Let us take a totally degenerate $r$-pointed stable curve $\msX'$ of genus $g$ over $k$ (cf. ~\cite[Definition 6.25]{Wak20} for the definition of a totally degenerate curve).
Since $\mcO p_{\N, \rho, g, r}^{^\mr{Zzz...}}$ is nonempty, 
the various properties of $\mcO p_{\N, \rho, g, r}^{^\mr{Zzz...}}$ recalled  in Theorem \ref{Th8}, together with the irreducibility of $\overline{\mcM}_{g, r}$,  imply that  the projection $\Pi'_{\rho, g, r} : \mcO p_{\N, \rho, g, r}^{^\mr{Zzz...}} \rightarrow \overline{\mcM}_{g, r}$ is surjective.
The fiber of this projection over the point classifying $\msX'$ is nonempty, so 
 there exists a dormant $\mr{PGL}_2^{(\N)}$-oper $\msE^\spadesuit$ on $\msX'$ of radii $\rho$.

Next, let $G$ denote  the  dual semi-graph of $\msX'$, which is a connected  trivalent  semi-graph of type $(g, r)$, in the sense of ~\cite[Definition 6.6, (ii) and (iii)]{Wak20}.
For a positive integer $\N_0$, recall from ~\cite[Definition 10.17]{Wak20}
the notion of a {\it balance $(p, \N_0)$-edge numbering of prescribed radii} on $G$ (where $G$ is equipped with an ordering on the branches incident to each vertex).
Roughly speaking, such a numbering  is  a certain collection of nonnegative integers $(a_e)_{e \in E}$  indexed by the edge set $E$ of $G$.
It follows from  ~\cite[Proposition 10.18]{Wak20} that, for each  $r$-tuple  $\rho_0 \in \Xi_{2, \N_0}^{r\ma}$,
 there exists a bijective correspondence 
 \begin{align} \label{Eq801h}
 \left(\begin{matrix} \text{the set of isomorphism classes} \\ \text{of  dormant $\mr{PGL}_2^{(\N_0)}$-opers} \\ \text{on $\msX'$ of radii $\rho_0$}\end{matrix} \right) \xrightarrow{\sim}
 \left( \begin{matrix} \text{the set of} \\ \text{balanced $(p, \N_0)$-edge} \\ \text{numberings on $G$ of radii  $\rho_0$}\end{matrix}\right).
 \end{align}
 In particular, the dormant $\mr{PGL}_2^{(\N)}$-oper $\msE^\spadesuit$ corresponds to a balanced $(p, \N)$-edge numbering $(a_e)_{e \in E}$ on $G$ of radii $\rho$.
 If $e_i$ ($i=1, \cdots, r$) denotes the open  edge corresponding to the $i$-th marked point $\sigma_i$,
 then the element 
 $ \overline{\left( \frac{2a_{e_i}+1}{2}\right)}$  in  $(\mbZ/p^\N \mbZ)^\times / \{ \pm 1 \}$ coincides with $\rho_i$ under the identification 
 $\Xi_{2, \N} = (\mbZ/p^\N \mbZ)^{\times}/\{ \pm 1 \}$ given by  \eqref{Eq311}.

 Now, define $\widetilde{\rho}_i$ to  be the  class $\overline{\left( \frac{2a_{e_i}+1}{2}\right)}$, regarded as an  element of $\Xi_{2, \N +1} \left(= (\mbZ/p^{\N +1} \mbZ)^\times / \{ \pm 1 \} \right)$.
The collection $\widetilde{\rho} := (\widetilde{\rho}_i)_{i=1}^r$ satisfies
$\widetilde{\rho}^{(\N)} = \rho$.    
Then, the collection $(a_e)_{e \in E}$  can be regarded  as a balanced $(p, \N +1)$-edge numbering on  $G$, and its radii are given by  $\widetilde{\rho}$.
Applying the  bijection \eqref{Eq801h} to the case $\N_0= \N +1$,
we obtain a dormant $\mr{PGL}_2^{(\N +1)}$-oper on $\msX'$ of radii $\widetilde{\rho}$ corresponding to $(a_e)_e$.
This shows that  $\mcO p_{\N +1, \widetilde{\rho}, g, r}^{^\mr{Zzz...}} \neq \emptyset$,  
 completing  the proof of this proposition.
\end{proof}

\LSP
\subsection{The adjoint bundle of  a dormant $\mr{PGL}_n^{(\N)}$-oper} \label{SS432}

Let $\msX$ be as before and  $\msE^\spadesuit := (\mcE_B, \STR)$ 
 a dormant $\mr{PGL}_n^{(\N)}$-oper on $\msX$.
Denote by $\pi : \mcE \rightarrow X$ (the structure morphism of) the $\mr{PGL}_n$-bundle induced from  $\mcE_B$.
Since $\mr{PGL}_n$ is affine, $\mcE$ corresponds to  an  $\mcO_X$-algebra $\mcO_\mcE$.
Pulling-back the log structure of $X^\mr{log}$ along  $\pi$ yields 
 a log structure on $\mcE$, and  we denote the resulting log scheme by $\mcE^\mr{log}$.
The $\mr{PGL}_n$-action on $\mcE$ induces a $\mr{PGL}_n$-action on the direct image $\pi_* (\mcT_{\mcE^\mr{log}/S^\mr{log}})$ of $\mcT_{\mcE^\mr{log}/S^\mr{log}}$, so we obtain 
\begin{align} \label{YY671}
\widetilde{\mcT}_{\mcE^\mr{log}/S^\mr{log}} :=  \pi_* (\mcT_{\mcE^\mr{log}/S^\mr{log}})^{\mr{PGL}_n},
\end{align}
i.e.,  $\widetilde{\mcT}_{\mcE^\mr{log}/S^\mr{log}}$ is the subsheaf of $\mr{PGL}_n$-invariant sections of $\pi_* (\mcT_{\mcE^\mr{log}/S^\mr{log}})$.

Denote by $\mfg_\mcE$ (resp., $\mfb_{\mcE_B}$) the adjoint vector bundle of $\mcE$ (resp., $\mcE_B$), i.e., the vector bundle on $X$ associated to $\mcE$ (resp., $\mcE_B$) via change of structure group by the adjoint representation $\mr{PGL}_n \rightarrow \mr{GL}(\mfg)$ (resp., $B \rightarrow \mr{GL}(\mfb)$).
We regard  $\mfg_{\mcE}$  as an $\mcO_X$-submodule of $\mcE nd_{\mcO_S} (\mcO_\mcE)$ under the natural identification $\mfg_\mcE = \pi_*(\mcT_{\mcE/X})^{\mr{PGL}_2} \left(= \pi_* (\mcT_{\mcE^\mr{log}/X^\mr{log}})^{\mr{PGL}_2} \right)$.
Differentiating $\pi$ gives rise to  a short exact sequence of $\mcO_X$-modules 
\begin{align} \label{YY662}
0 \rightarrow \mfg_\mcE \xrightarrow{\mr{inclusion}} \widetilde{\mcT}_{\mcE^\mr{log}/S^\mr{log}} \xrightarrow{d_\mcE} \mcT \rightarrow 0
\end{align}
(cf. ~\cite[Section 1.2.5]{Wak5}).

Note that $\STR$ induces a $\mcD^{(\N -1)}$-module structure
 \begin{align} \label{YY788}
 \DMO^\mr{ad}_\STR : {^L}\mcD^{(\N -1)} \rightarrow  \mcE nd_{\mcO_S}(\mfg_\mcE)
 \end{align}
on $\mfg_\mcE$.
To be precise,  if $\STR^{\natural \natural}$ denotes the $\mcD^{(\N -1)}$-module structure on $\mcO_\mcE$ corresponding to $\STR$ (cf. ~\cite[Remark 2.6]{Wak20}), then   $\DMO^\mr{ad}_\STR$ is 
   given by $\DMO^\mr{ad}_\STR (D) (v) = [\STR^{\natural \natural} (D), v]$ for any local sections $D \in \mcD^{(\N -1)}$ and $v \in \mfg_\mcE$.
In particular, we obtain a $\mcD^{(\N-1)}$-module $(\mfg_\mcE, \DMO_\STR^\mr{ad})$,  which has vanishing    $p^\N$-curvature.

The following assertion generalizes ~\cite[Proposition 8.2]{Wak20} to the  higher-rank setting.

\SSP
 \bpr \label{Pee4}
Let us fix 
 $a \in \{ 0, 1, \cdots, \N -1\}$.
 Recall  that $(\mfg_\mcE, \DMO_\STR^{\mr{ad}})$ induces  a flat module 
 $(\mfg_\mcE^{[a]}, \DMO_\STR^{\mr{ad}[a]} : \mfg_\mcE^{[a]} \rightarrow \Omega^{(a)} \otimes \mfg_\mcE^{[a]})$ on $X^{(a)\mr{log}}/S^\mr{log}$ with vanishing $p$-curvature 
 (cf. \eqref{J16}).
  \begin{itemize}
 \item[(i)]
 The
following equalities of $\mcO_S$-modules  hold: 
  \begin{align} \label{Er224}
 \mbR^2 f_* (\mcK^\bullet 
 [\DMO^{\mr{ad}[a]}_\phi])
 =f_*(\mr{Ker}(\DMO_\phi^{\mr{ad} [a]})) = \mbR^1 f_* (\mr{Coker}(\DMO^{\mr{ad}[a]}_\phi)) = 0.
 \end{align}
 \item[(ii)]
 Denote by $D^{(a)}$  
 the reduced effective divisor  on $X^{(a)}$ determined by the union of the marked points of $\msX^{(a)}$.
 Also, for each flat module $(\mcV, \nabla_\mcV)$ on $X^{(a)\mr{log}}/S^\mr{log}$, we write ${^c}\mcV := \mcV (-D^{(a)})$ and write ${^c}\nabla_\mcV$ for the $S^\mr{log}$-connection on ${^c}\mcV \left(\subseteq \mcV \right)$ obtained by restricting $\nabla_\mcV$.
 Then,  we have
  \begin{align} \label{Er225}
 \mbR^2 f_* (\mcK^\bullet 
 [{^c}\DMO^{\mr{ad}[a]}_\phi])
 =f_*(\mr{Ker}({^c}\DMO_\phi^{\mr{ad} [a]})) = \mbR^1 f_* (\mr{Coker}({^c}\DMO^{\mr{ad}[a]}_\phi)) = 0.
 \end{align}
 \item[(iii)]
 $\mbR^1 f_* (\mcK^\bullet [\DMO^{\mr{ad}[a]}_\phi])$
 is a vector bundle  on $S$ of rank 
 $(n^2-1)(2g-2 +r)$.
 Also,  
  $\mbR^1f_*(\mr{Ker}(\DMO_\STR^{\mr{ad} [a]}))$ and $f_* (\mr{Coker}(\DMO^{\mr{ad}[a]}_\STR))$) are  
 vector bundles on $S$   of 
 rank  
  $(n^2-1)(g-1) +   \frac{1}{2}(n^2 -n) r$
  and 
 $(n^2-1)(g-1) + \frac{1}{2} (n^2 + n-2) r$,
 respectively.
   \end{itemize}
 \epr
 \begin{proof}
 We shall prove assertions (i)-(iii)  simultaneously  by induction on $a$.
 For simplicity, we write 
 $\mcF^{[a]} := \mfg_\mcE^{[a]}$ and  $\nabla^{[a]} := \DMO_\STR^{\mr{ad}[a]}$.

Let us begin with   
the base step, i.e., 
 the case  $a = 0$.
 Assertions (i) and (iii) in this case
  were  already established  in ~\cite[Propositions 6.5, (ii), and  6.18]{Wak5}.
  To prove assertion (ii),
we  recall from ~\cite[Section 6.5.1]{Wak5} that 
the $1$-st monodromy operator 
\begin{align}
\mu_i (\nabla_\phi^{\mr{ad}})^{\langle 1 \rangle} : \sigma^*_i (\mcF^{[0]}) \rightarrow \sigma^*_i (\mcF^{[0]})
\end{align}
(for every $i=1, \cdots, r$) of  $\nabla_\phi^{\mr{ad}}$ at $\sigma_i$ can be identified with the adjoint operator 
 on $\sigma^*_i (\mcF^{[0]})$ determined by the monodromy  
 (in the sense of ~\cite[Definition 1.46]{Wak5}) of the $S^\mr{log}$-connection induced  from  $\phi$.
This endomorphism  fits into the following short exact sequence of complexes (concentrated at degrees $0$ and $1$)  under the natural identification $\sigma_i^* (\mcF^{[0]}) = \sigma^*_i (\Omega \otimes \mcF^{[0]})$ arising  from the residue isomorphism $\sigma_i^* (\Omega) \cong \mcO_S$:
  \begin{align} \label{Er500}
\vcenter{\xymatrix@C=46pt@R=36pt{
 0 \ar[r] &{^c}\mcF^{[0]} \ar[r] \ar[d]^-{{^c}\nabla^{[0]}} & \mcF^{[0]}  \ar[r] \ar[d]^-{\nabla^{[0]}} & \bigoplus_{i=1}^r \sigma_{i*} (\sigma^*_i (\mcF^{[0]})) \ar[r] \ar[d]^-{\bigoplus_i \mu_i (\nabla_\phi^{\mr{ad}})^{\langle 1 \rangle}} & 0
 \\
0 \ar[r]  &\Omega \otimes {^c}\mcF^{[0]}  \ar[r]& \Omega \otimes \mcF^{[0]} \ar[r]&  \bigoplus_{i=1}^r \sigma_{i*} (\sigma^*_i (\mcF^{[0]})) \ar[r] & 0.
 }}
\end{align}
The equality   $f_* (\mr{Ker}(\nabla^{[0]})) = 0$ holds by assertion (i), so 
this short exact sequence gives rise to
a long exact sequence of hypercohomology sheaves 
\begin{align}
0 &\rightarrow \bigoplus_{i=1}^r \mr{Ker} (\mu_i (\nabla_\phi^{\mr{ad}})^{\langle 1 \rangle}) 
\rightarrow \mbR^1 f_* (\mcK^\bullet [{^c}\nabla^{[0]}])
\rightarrow \mbR^1 f_* (\mcK^\bullet [\nabla^{[0]}]) \\
&\rightarrow 
\bigoplus_{i=1}^r \mr{Coker} (\mu_i (\nabla_\phi^{\mr{ad}})^{\langle 1 \rangle}) 
\rightarrow \mbR^2 f_* (\mcK^\bullet [{^c}\nabla^{[0]}])
\rightarrow \mbR^2 f_* (\mcK^\bullet [\nabla^{[0]}]) \rightarrow 0, \notag
\end{align}
which is functorial with respect to $S$.
According to  ~\cite[Proposition 6.17]{Wak5},
$\mr{Ker} (\mu_i (\nabla_\phi^{\mr{ad}})^{\langle 1 \rangle})$ and 
$\mr{Coker} (\mu_i (\nabla_\phi^{\mr{ad}})^{\langle 1 \rangle})$ are vector bundles on $S$ of the same rank.

On the other hand,
since the Killing form on $\mfg$ yields an isomorphism of flat modules  $(\mcF^{[0]}, \nabla^{[0]}) \xrightarrow{\sim} (\mcF^{[0]\vee}, \nabla^{[0]\vee})$,
it follows from ~\cite[Corollary 6.16]{Wak5}
that $\mbR^1 f_* (\mcK^\bullet [{^c}\nabla^{[0]}])$
and  $\mbR^1 f_* (\mcK^\bullet [\nabla^{[0]}])$ are vector bundles on $S$ of the same rank.
Hence, the fourth arrow 
$\mbR^1 f_* (\mcK^\bullet [\nabla^{[0]}]) \rightarrow \bigoplus_{i=1}^r \mr{Coker} (\mu_i (\nabla_\phi^{\mr{ad}})^{\langle 1 \rangle})$
in the above exact sequence must be surjective, and the morphism 
$\mbR^2 f_* (\mcK^\bullet [{^c}\nabla^{[0]}])
\rightarrow \mbR^2 f_* (\mcK^\bullet [\nabla^{[0]}])$ induced by  the inclusion $\mcK^\bullet [{^c}\nabla^{[0]}] \hookrightarrow \mcK^\bullet [\nabla^{[0]}]$ turns out to be injective.
This fact together with assertion (i) in the base step shows $\mbR^2 f_* (\mcK^\bullet [{^c}\nabla^{[0]}]) = 0$.

Moreover,  consider 
the conjugate spectral sequence
\begin{align} \label{dE2001}
''E_2^{q_1, q_2} := \mbR^{q_1} f_*(\mcH^{q_2}(\mcK^\bullet [{^c}\DMO^{[0]}])) \Rightarrow \mbR^{q_1 + q_2} f_* (\mcK^\bullet [{^c}\DMO^{[0]}])
\end{align}
associated to $\mcK^\bullet [{^c}\DMO^{[0]}]$ (cf. ~\cite[Eq.\,(757)]{Wak5}).
Since the equality $\mbR^j f_* (\mr{Ker}({^c}\DMO^{[0]})) = 0$ holds for every $j \geq 2$ due to   $\mr{dim}(X/S) =1$,  
the morphism
$\mbR^2 f_* (\mcK^\bullet [{^c}\DMO^{[0]}]) \rightarrow \mbR^1 f_* (\mr{Coker}({^c}\DMO^{[0]}))$
induced by \eqref{dE2001} is surjective, which implies  $\mbR^1 f_* (\mr{Coker}({^c}\DMO^{[0]})) = 0$.
The equality  $f_* (\mr{Ker}({^c}\nabla^{[0]})) = 0$ follows from $\left(f_* (\mr{Ker}({^c}\nabla^{[0]})) \subseteq \right) f_* (\mr{Ker}(\nabla^{[0]})) = 0$ obtained above, 
so the proof of the base step is completed.

Next, to carry out   the induction step, we assume  that  assertions (i)-(iii) with $a$ replaced by $a-1$ (where $a \geq 1$) have been proved.
We begin by considering assertion (i).
Since $\mr{Ker} (\DMO^{[a]}) \subseteq \mr{Ker}(\DMO^{[a-1]})$,
the equality 
\begin{align} \label{YY900}
f_* (\mr{Ker}(\DMO^{[a]})) = 0
\end{align}
follows immediately from  the induction hypothesis $f_* (\mr{Ker}(\DMO^{[a -1]})) = 0$.
Next,   we consider the Hodge-to-de Rham spectral sequence
\begin{align} \label{dE180}
'E_{1}^{q_1, q_2} := \mbR^{q_2} f_* (\mcK^{q_1} [\DMO^{[a]}]) \Rightarrow \mbR^{q_1 + q_2} f_* (\mcK^\bullet [\DMO^{[a]}])
\end{align}
associated to the complex $\mcK^\bullet [\DMO^{[a]}]$ (cf. ~\cite[Eq.\,(755)]{Wak8}).
Since  $\mr{dim}(X/S) =1$, the equality  $\mbR^2 f_* (\mcF^{[a]}) = 0$  holds, and 
the sequence
\begin{align} \label{UU1}
\mbR^1f_* (\mcF^{[a]}) \xrightarrow{\mbR^1 f_* (\nabla^{[a]})}\mbR^1 f_* (\Omega^{(a)} \otimes \mcF^{[a]}) \rightarrow \mbR^2 f_* (\mcK^\bullet [\DMO^{[a]}]) \rightarrow 0
\end{align}
 induced by 
\eqref{dE180} is exact.
According to 
~\cite[Theorem 3.1.1]{Ogu2} (or the discussion
  following ~\cite[Proposition 1.2.4]{Ogu1}),
$\mbR^1 f_* (\Omega^{(a)} \otimes \mcF^{[a]})$ is isomorphic to
$\mbR^1 f_* (\mr{Coker}(\DMO^{[a-1]}))$.
By the induction hypothesis, 
 the equality  $\mbR^1 f_* (\mr{Coker}(\DMO^{[a-1]})) =0$ holds, so the exactness of  \eqref{UU1} implies
\begin{align} \label{dE189}
\mbR^2 f_* (\mcK^\bullet [\DMO^{[a]}]) = 0.
\end{align}

Moreover, 
the conjugate spectral sequence
\begin{align} \label{dE200}
''E_2^{q_1, q_2} := \mbR^{q_1} f_*(\mcH^{q_2}(\mcK^\bullet [\DMO^{[a]}])) \Rightarrow \mbR^{q_1 + q_2} f_* (\mcK^\bullet [\DMO^{[a]}])
\end{align}
associated to $\mcK^\bullet [\DMO^{[a]}]$  gives rise to 
the surjection
$\mbR^2 f_* (\mcK^\bullet [\DMO^{[a]}]) \rightarrow \mbR^1 f_* (\mr{Coker}(\DMO^{[a]}))$.
Hence, it follows from \eqref{dE189} that 
\begin{align} \label{dE193}
\mbR^1 f_* (\mr{Coker}(\DMO^{[a]})) = 0.
\end{align}
By \eqref{YY900}, \eqref{dE189}, and \eqref{dE193},
the proof of assertion (i) is completed.

Assertion (ii) follows immediately from an argument similar to the corresponding assertion  in the base step, using  the assertion (i) just proved  
and the conjugate spectral sequence $''E_2^{q_1, q_2} := \mbR^{q_1} f_*(\mcH^{q_2}(\mcK^\bullet [{^c}\DMO^{[a]}]))$
$\Rightarrow \mbR^{q_1 + q_2} f_* (\mcK^\bullet [{^c}\DMO^{[a]}])$ instead of  \eqref{dE2001}.

Finally, we consider assertion (iii).
It follows from ~\cite[Proposition 2.15, (i)]{Wak20} that  $\mr{Ker}(\nabla^{[a]})$ is a relatively torsion-free sheaf on $X^{(a+1)}$.
Also, since
 $\mr{Coker}(\nabla^{[a]})$ is isomorphic to $\Omega^{(a+1)} \otimes \mr{Ker}(\nabla^{[a]})$ 
 (cf. ~\cite[Theorem 3.1.1]{Ogu2} or the comment following ~\cite[Proposition 1.2.4]{Ogu1}), 
 $\mr{Coker}(\nabla^{[a]})$ is a relatively torsion-free sheaf on $X^{(a+1)}$.
By the equality $\mbR^2 f_* (\mr{Ker}(\nabla^{[a]})) = 0$ and  \eqref{dE193},
both $\mbR^1 f_*(\mr{Ker}(\nabla^{[a]}))$ and $f_*(\mr{Coker}(\nabla^{[a]}))$
turn out to be vector bundles (cf. ~\cite[Chapter III, Theorem 12.11, (b)]{Har}).
Moreover, the exactness of  
the sequence 
\begin{align} \label{UU2}
0 \rightarrow  f_* (\Omega^{(a)} \otimes \mcF^{[a]}) \left(= f_* (\mr{Coker}(\DMO^{[a-1]})) \right) &\rightarrow \mbR^1 f_* (\mcK^\bullet [\DMO^{[a]}])\\
 &\rightarrow \mbR^1 f_* (\mcF^{[a]}) \left(= \mbR^1 f_* (\mr{Ker}(\DMO^{[a-1]})) \right)\rightarrow 0  \notag
\end{align}
induced by \eqref{dE180} and  the induction hypothesis together  imply that
$\mbR^1 f_* (\mcK^\bullet [\nabla^{[a]}])$ is a vector bundle.

 In what follows, we  compute the ranks of the vector bundles under consideration.
 To this end, we may assume that $S = \mr{Spec}(k)$.
 Fix a dormant $n^{(\N)}$-theta characteristic $\vartheta := (\varTheta, \nabla_\vartheta)$ of $\msX$.
 The dormant $\mr{PGL}_n^{(\N)}$-oper $\msE^\spadesuit$ corresponds to a dormant $(\mr{GL}_n^{(\N)}, \vartheta)$-oper $\nabla^\diamondsuit$.
 For each $i=1, \cdots, r$,  let $[d_{i, 1}, \cdots, d_{i, n}]$  be the exponent of $\nabla^\diamondsuit$ at the marked point $\sigma_i$.
 After possibly reordering the elements in this multiset, we may assume that  $0 \leq \widetilde{d}_{i, 1} < \widetilde{d}_{i, 2} < \cdots < \widetilde{d}_{i, n} < p^\N$. 

By  ~\cite[Corollary 6.16]{Wak5} in the case where ``$(\mcF, \nabla)$" is taken to be $(\mcF^{[a]}, \DMO^{[a]})$,
there exists a canonical isomorphism of $k$-vector spaces
 \begin{align}
H^1 (X^{(a)},  \mr{Ker} ({^c}(\DMO^{[a]\vee}))) \xrightarrow{\sim} H^0 (X^{(a)}, \mr{Coker}(\DMO^{[a]}))^\vee.
 \end{align}
 In particular,
 we have
 \begin{align} \label{ddE1}
 h^1 (\mr{Ker} ({^c}(\DMO^{ [a]\vee}))) = h^0 (\mr{Coker}(\DMO^{[a]})).
 \end{align}
 Next,   
 the isomorphism of $\mcD^{(\N -1)}$-modules
 $(\mfg_\mcE, \DMO_\phi^{\mr{ad}}) \xrightarrow{\sim} (\mfg_\mcE^\vee, (\DMO_\phi^\mr{ad})^\vee)$
 induced by the Killing form on $\mfg$ 
 yields  an isomorphism of flat modules
$({^c}\mcF^{[a]}, {^c}\DMO^{[a]}) \xrightarrow{\sim}
({^c}(\mcF^\vee)^{[a]}, {^c}(\DMO^\vee)^{[a]})$.
This restricts to   an isomorphism  $\mr{Ker}({^c}\DMO^{[a]}) \xrightarrow{\sim} \mr{Ker}({^c}(\DMO^\vee)^{[a]})$, which implies
\begin{align} \label{dE155}
h^1 (\mr{Ker}({^c}\DMO^{[a]})) = h^1 (\mr{Ker}({^c}(\DMO^\vee)^{[a]})).
\end{align}
Hence,
the following chain  of equalities holds:
\begin{align} \label{dE149}
h^0 (\mr{Coker}(\DMO^{[a]})) & \stackrel{\eqref{ddE1}}{=} h^1 (\mr{Ker}({^c}(\DMO^{[a]\vee}))) \\
& = h^1 (\mr{Ker}({^c}(\DMO^\vee)^{[a]})) -R_\IN \notag \\
&\stackrel{\eqref{dE155}}{=} h^1 (\mr{Ker}({^c}\DMO^{[a]})) -R_\IN  \notag \\
& = \left(h^1 (\mr{Ker}(\DMO^{[a]})) + (n-1)r + R_\IN \right) -R_\IN \notag \\
& = h^1 (\mr{Ker}(\DMO^{[a]})) +(n-1)r \notag 
\end{align}
(cf. \eqref{dE140} for the definition of $R_\IN$), 
where the second and  fourth equalities follow from Lemma \ref{L010} stated  below.

On the other hand,  
there exists a short exact sequence
\begin{align} \label{ddE6}
0 \rightarrow H^1 (X^{(a)}, \mr{Ker}(\DMO^{[a]})) \rightarrow \mbH^1 (X^{(a)}, \mcK^\bullet [\DMO^{[a]}]) \rightarrow H^0 (X^{(a)}, \mr{Coker}(\DMO^{[a]}))
\rightarrow 0 
\end{align}
induced by  \eqref{dE200}.
 This implies 
\begin{align} \label{eq+}
& \ \  \ \ h^1 (\mr{Ker}(\nabla^{[a]}))  + h^0 (\mr{Coker}(\nabla^{[a]}))  \\
 & = \mr{dim}(\mbH^1 (X^{(a)}, \mcK^\bullet [\nabla^{[a]}]))  \notag \\
&= h^1 (\mr{Ker}(\nabla^{[a-1]})) + h^0 (\mr{Coker}(\nabla^{[a-1]})) \notag  \\
& = 
\left((n^2 -1) (g-1) + \frac{1}{2} (n^2 -n) r  \right) + \left( (n^2-1) (g-1)+ \frac{1}{2} (n^2 + n-2)r\right)
 \notag \\
& =  (n^2-1)(2g-2+r), \notag
\end{align} 
where the second equality follows from the exactness of \eqref{UU2} and the third equality follows from the induction hypothesis.
Consequently,     it follows from 
\eqref{dE149} and \eqref{eq+} that 
\begin{align}
h^1(\mr{Ker}(\nabla^\mr{ad})) & = (n^2 -1) (g-1) + \frac{1}{2} (n^2 -n) r, \\
 h^0(\mr{Coker}(\nabla^\mr{ad})) &= (n^2-1) (g-1)+ \frac{1}{2} (n^2 + n-2)r. \notag
\end{align}
This completes the proof of the proposition.
 \end{proof}
\SSP

\bco \label{Co22}
 $\mbR^1 f_* (\mcS ol (\DMO_\STR^{\mr{ad}}))$ (resp.,  $f_* (\Omega^{(\N)}\otimes \mcS ol (\DMO_\STR^{\mr{ad}}))$) is a  vector bundle on $S$  of rank $(n^2-1)(g-1) +   \frac{1}{2}(n^2 -n) r$ (resp., $(n^2-1)(g-1) +   \frac{1}{2}(n^2 +n -2) r$). 
\eco
\begin{proof}
By the definition of $\nabla_\phi^{\mr{ad}[\N -1]}$, 
 the equality
$\mcS ol (\DMO_\STR^{\mr{ad}}) =  \mr{Ker}(\nabla_\phi^{\mr{ad}[\N -1]})$
holds.
On the other hand, it follows from  ~\cite[Theorem 3.1.1]{Ogu2} (or the discussion following ~\cite[Chapter I, Proposition 1.2.4]{Ogu1}) that the Cartier operator associated to $(\mfg_\mcE^{[\N -1]}, \nabla_\phi^{\mr{ad}[\N -1]})$ induces an isomorphism  $ \mr{Coker} (\nabla_\phi^{\mr{ad}[\N -1]})\xrightarrow{\sim}  \Omega^{(\N)}\otimes \mr{Ker} (\nabla_\phi^{\mr{ad}[\N -1]})$.
These facts together imply
$\mbR^1 f_* (\mcS ol (\DMO_\STR^{\mr{ad}})) \cong \mbR^1 f_* (\mr{Ker}(\nabla_\phi^{\mr{ad}[\N -1]}))$
and $f_* (\Omega^{(\N)}\otimes \mcS ol (\DMO_\STR^{\mr{ad}})) \cong f_* ( \mr{Coker}(\nabla_\phi^{\mr{ad}[\N -1]}))$.
Hence, the assertion 
 follows from Proposition \ref{Pee4}, (ii).
\end{proof}
\SSP

The following lemma was applied in the proof of Proposition \ref{Pee4} above.
See ~\cite[Lemma 8.3]{Wak20} for the case of $n=2$.

\SSP
\ble \label{L010}
We retain  the notation in the proof of Proposition  \ref{Pee4}.
Let $\IN \in \{1, \cdots, \N -1\}$, and 
suppose that $S = \mr{Spec}(k)$.
Also, 
 we  set 
\begin{align} \label{dE140}
R_{\IN} := \sharp \left\{ (i, j, j')  \, \Big| \, 
1 \leq i \leq r,   \  1 \leq j' < j \leq n, \ (d_{i, j} -d_{i, j'})_{[\IN]} \in \{0, p-1 \} \left(\subseteq \mbF_p \right)
\right\}.
   \end{align}
Then, 
the following two equalities hold:
\begin{align} \label{dE141}
h^1 (\mr{Ker}({^c}\DMO^{[\IN]})  &= h^1 (\mr{Ker}(\DMO^{[\IN]})) + (n-1) r + R_\IN,   \\
h^1 (\mr{Ker}({^c}(\DMO^\vee)^{[\IN]})) &= h^1 (\mr{Ker}({^c}(\DMO^{[\IN]\vee})))   +  R_\IN. 
\notag
\end{align}
\ele 
\begin{proof}
To begin with, we set up  some notation.
Let $U_\oslash$ and $\mcD_\oslash^{(\N -1)}$ be as in Section \ref{SS90}.
Let $\msV := (\mcV, \DMO)$ be a  $\mcD^{(\N-1)}$-module (resp., a $\mcD_\oslash^{(\N-1)}$-module) such that $\psi (\DMO) = 0$ and $\mcV$ is a vector bundle of positive rank.
In the non-resp'd case, $({^c}\mcV^{[a]}, {^c}\DMO^{[a]})$ is defined as in the proof of Proposition  \ref{Pee4}.
In the resp'd case, we set ${^c}\mcV^{[a]} := t^{p^a} \cdot \mcV^{[a]} \left(\subseteq \mcV^{[a]} \right)$, which has an $S^\mr{log}$-connection  ${^c}\DMO^{[a]}$ obtained by restricting $\DMO^{[a]}$.
Then, the  inclusion 
${^c}\mcV^{[\IN]}  \hookrightarrow \mcV^{[\IN]}$ restricts to an inclusion
\begin{align} \label{dE123}
\alpha_{\msV}^{[\IN]} : \mr{Ker} ({^c}\DMO^{[\IN]}) \hookrightarrow \mr{Ker} (\DMO^{[\IN]}).
\end{align}
On the other hand, the natural morphism
 ${^c}(\mcV^\vee)^{[a]}
 \rightarrow 
{^c}(\mcV^{[\IN]\vee})$
restricts to a morphism of $\mcO_{X^{(\IN+1)}}$-modules
\begin{align} \label{UU7}
\beta_\msV^{[\IN]} : \mr{Ker}({^c} (\DMO^\vee)^{[\IN]}) \rightarrow \mr{Ker} ({^c}(\DMO^{[\IN]\vee})).
\end{align}

In what follows, we examine  the morphisms $\alpha_\msV^{[\IN]}$ and $\beta_\msV^{[\IN]}$ in the case where   $\msV$ is taken to be the $\mcD_\oslash^{(\N-1)}$-module  $\msO_{\oslash, \EX}^{(\N -1)} := (\mcO_\oslash, \DMO_{\oslash, \EX}^{(\N -1)})$ for 
 $\EX \in (\mbZ/p^\N \mbZ)^\times$ (cf. \eqref{UU3}).
For simplicity, write $\msO_{d} := \msO_{\oslash, \EX}^{(\N -1)}$ and $\DMO_{\EX}:= \DMO_{\oslash, \EX}^{(\N -1)}$.
First, 
 a straightforward calculation  shows 
\begin{align} \label{UU6}
\mr{Ker} ({^c}\DMO^{[\IN]}_{\EX}) = \begin{cases} t^{\widetilde{\EX}_{[0, \IN]}} \cdot \mcO^{(\IN+1)}_{\oslash}& \text{if $\EX_{[\IN]} \neq 0$}; 
\\
 t^{\widetilde{d}_{[0, \IN]}+p^{\IN +1}} \cdot \mcO^{(\IN+1)}_{\oslash}
  & \text{if $\EX_{[\IN]} = 0$}. \end{cases}
\end{align}
Hence, it follows from ~\cite[Proposition 4.4]{Wak20} that $\alpha_{\msO_\EX}^{[\IN]}$ is injective and its cokernel satisfies
\begin{align} \label{dE128}
\mr{length}_{\mcO_\oslash^{(a+1)}}(\mr{Coker}(\alpha^{[\IN]}_{\msO_\EX})) = \begin{cases} 0 & \text{if $\EX_{[\IN]} \neq 0$};\\ 1 & \text{if $\EX_{[\IN]} =0$}.  \end{cases}
\end{align}

Also, by putting
 $c := -\EX$, we see (from $\EX \in (\mbZ/p^\N \mbZ)^\times$) 
that 
\begin{align} \label{UU56}
\mr{Ker} ({^c}(\DMO_{\EX}^{\vee})^{[\IN]}) \left(= \mr{Ker}({^c} \DMO_{c}^{[\IN]}) \right) = \begin{cases} t^{\widetilde{c}_{[0, \IN]}} \cdot \mcO_{\oslash}^{(\IN+1)} &  \text{if
$c_{[\IN]} \neq 0 \left(\Leftrightarrow \EX_{[\IN]} \neq p-1 \right)$};
\\
t^{\widetilde{c}_{[0, \IN]}+p^{\IN+1}} \cdot \mcO_\oslash^{(\IN+1)} &
\text{if
$c_{[\IN]} = 0 \left(\Leftrightarrow \EX_{[\IN]} = p-1 \right)$}.
\end{cases}
\end{align}
and that
\begin{align} \label{UU5}
 \mr{Ker}({^c}(\DMO_{\EX}^{[\IN]})^{\vee}) = 
t^{-\widetilde{d}_{[0, \IN]}+p^{\IN+1}} \cdot \mcO_\oslash^{(\IN+1)} \left(= t^{\widetilde{c}_{[0, \IN]}} \cdot \mcO_\oslash^{(\IN+1)} \right).
\end{align}
Hence, 
$\beta_{\msO_\EX}^{[\IN]}$ is injective and its cokernel satisfies
\begin{align} \label{dE124}
\mr{length}_{\mcO_\oslash^{(a+1)}} (\mr{Coker}(\beta^{[\IN]}_{\msO_\EX})) =  \begin{cases} 1 & \text{if   $\EX_{[\IN]} = p-1$}; \\
0 & \text{if   $\EX_{[\IN]}\neq  p-1$}.
\end{cases}
 \end{align}
 The same  arguments with $d$ replaced by $c$ give
 \begin{align} \label{dE125}
\mr{length}_{\mcO_\oslash^{(a+1)}} (\mr{Coker}(\alpha^{[\IN]}_{\msO_{-\EX}})) &= \begin{cases} 0 & \text{if $\EX_{[\IN]} \neq p-1$};\\ 1 & \text{if $\EX_{[\IN]} =p-1$},  \end{cases}  \\
\mr{length}_{\mcO_\oslash^{(a+1)}} (\mr{Coker}(\beta^{[\IN]}_{\msO_{-\EX}})) &=  \begin{cases} 1 & \text{if   $\EX_{[\IN]} =  0$}; \\
0 & \text{if   $\EX_{[\IN]} \neq  0$}. \notag
\end{cases}
 \end{align}

Now,  let us return to our setting.
We set $(\mcF, \DMO) := (\mfg_\mcE, \DMO_\STR^{\mr{ad}})$.
Both $\alpha^{[a]}_{(\mcF, \DMO)}$ and $\beta^{[a]}_{(\mcF, \DMO)}$  become  isomorphisms
when restricted to $X \setminus \bigcup_{i=1}^r \mr{Im}(\sigma_i)$.
Given each $\EX \in \mbZ/p^\N \mbZ$, we define 
\begin{align}
\mcQ_{\EX}^{[\IN]} := \begin{cases} \mcO_\oslash^{(\IN+1)} /(t^{p^{\IN+1}}) & \text{if $\EX_{[\IN]} \in \{0, p-1 \}$};  \\ 0  & \text{if $\EX_{[a]} \notin \{0, p-1 \}$}.\end{cases}
\end{align}
For each $i =1, \cdots, r$, denote by $U_i$ the formal neighborhood of $\mr{Im}(\sigma_i)$ in $X$, and by 
  $\mr{incl}_{i} : U_i  \rightarrow X$ the natural morphism.
 Note that  $\nabla$ can be identified with the  $\mcD^{(\N -1)}$-module structure on $\overline{\mcE nd} (\mcF_\varTheta)$ induced naturally  from $\nabla^\diamondsuit$.
 Hence, since $\left(\bigoplus_{j=1}^r \msO^{(\N -1)}_{\oslash, d_{i, j}} \right)^\vee \otimes \left(\bigoplus_{j=1}^r \msO^{(\N -1)}_{\oslash, d_{i, j}} \right) \cong  \bigoplus_{1 \leq j, j' \leq n} \msO^{(\N -1)}_{\oslash, d_j -d_j'}$,
the restriction $(\mcF |_{U_i}, \nabla |_{U_i})$ of $(\mcF, \nabla)$ to $U_i$ admits an isomorphism
\begin{align} \label{Et66}
(\mcF |_{U_i}, \nabla |_{U_i}) \xrightarrow{\sim} \msO_0^{\oplus (n-1)} \oplus \bigoplus_{1 \leq j < j' \leq n} \msO_{d_{i, j} - d_{i, j'}} \oplus  \bigoplus_{1 \leq j' < j \leq n} \msO_{d_{i, j} - d_{i, j'}}
\end{align}
under a fixed identification $U_i = U_\oslash$.
By \eqref{dE128} and the first equality  of \eqref{dE125},
 the cokernel  of $\alpha_{(\mcF, \DMO)}^{[\IN]}$ restricted to $U_i$ is isomorphic 
 to
 \begin{align}
 \mcQ_i^{[\IN]} := 
 (\mcO_\oslash^{(\IN+1)} / (t^{p^{\IN+1}}))^{\oplus (n-1)} \oplus \bigoplus_{1 \leq j'  < j \leq n} \mcQ^{[\IN]}_{d_{i, j} -d_{i, j'}}.
 \end{align}
This implies that  $\alpha^{[\IN]}_{(\mcF, \DMO)}$ fits into the short exact sequence 
\begin{align} \label{dE132}
0 \rightarrow \mr{Ker}({^c}\DMO^{[\IN]}) \xrightarrow{\alpha^{[\IN]}_{(\mcF, \DMO)}} \mr{Ker}(\DMO^{[\IN]}) \rightarrow \bigoplus_{i=1}^r \mr{incl}_{i*} ( \mcQ_i^{[\IN]}) \rightarrow 0.
\end{align}

Here, we shall write  $H^j (-) := H^j (X^{(\IN+1)}, -)$ ($j \geq 0$) for simplicity.
Since $H^0 (\mr{Ker}(\DMO^{[\IN]})) 
 =  0$ (cf. \eqref{YY900}),
 the sequence \eqref{dE132} yields  a short exact sequence of $k$-vector spaces
\begin{align}
0 \rightarrow  \bigoplus_{i=1}^r H^0 (\mr{incl}_{i*} (\mcQ_i^{[\IN]}))
\rightarrow  H^1 (\mr{Ker}({^c}\DMO^{[\IN]}))
\rightarrow H^1 (\mr{Ker}(\DMO^{[\IN]}))
\rightarrow 0.
\end{align}
This sequence implies 
\begin{align}
h^1 (\mr{Ker}({^c}\DMO^{[\IN]})) &= h^1 (\mr{Ker}(\DMO^{[\IN]})) + \sum_{i=1}^r\mr{dim} (H^0(\mr{incl}_{i*} (\mcQ_i^{[\IN]}))) \\
& = h^1 (\mr{Ker}(\DMO^{[\IN]})) + \sum_{i=1}^r (n-1 + \sum_{1 \leq j' < j \leq n}\mr{length}_{\mcO_\oslash^{(\IN +1)}} (\mcQ^{[\IN]}_{d_{i, j}-d_{i, j'}})) \notag \\
& =  h^1 (\mr{Ker}(\DMO^{[\IN]}))    +(n-1)r +R_a. \notag
\end{align}
Thus, we have proved the first equality in \eqref{dE141}.

Moreover,  it follows from \eqref{dE124} and the second equality  of \eqref{dE125} that we have the following short exact sequence:
\begin{align} \label{dE130}
0 \rightarrow \mr{Ker}({^c}\DMO^{\vee [\IN]}) \xrightarrow{\beta^{[\IN]}_{(\mcF, \DMO)}} \mr{Ker}({^c}\DMO^{[\IN]\vee})  \rightarrow \bigoplus_{i=1}^r \bigoplus_{1 \leq j' < j \leq n} \mr{incl}_{i*} (\mcQ^{[\IN]}_{d_{i, j}-d_{i, j'}} )
\rightarrow  0.
\end{align}
Similarly  to the above argument, 
 \eqref{dE130} implies  the equality
\begin{align}
h^1 (\mr{Ker}({^c}(\DMO^\vee)^{[\IN]})) = h^1 (\mr{Ker}({^c}(\DMO^{[\IN]\vee}))) + R_{\IN},
\end{align}
which is precisely  the second equality in \eqref{dE141}.
This completes  the proof of this lemma.
\end{proof}

\LSP
\subsection{The level-reduction map} \label{SS241}

Given
 a $\mr{PGL}_n^{(\N +1)}$-oper $\msE^\spadesuit$ on $\msX$, we shall write $\msE^{\spadesuit (\N)}$ for the induced $\mr{PGL}_n^{(\N)}$-oper  via level reduction.
By the definition of $p^\N$-curvature,  this $\mr{PGL}_n^{(\N)}$-oper  is dormant.
For each $\rho \in \mfc_{\N} (k)^{r\ma}$,
we shall write
\begin{align} \label{Er26}
\mcO p_{\N +1, \Rightarrow \rho, g, r}
\end{align}
for the closed substack of $\mcO p_{\N +1, g, r}$ classifying $\mr{PGL}_n^{(\N +1)}$-opers $\msE^\spadesuit$ with $\rho_i (\msE^{\spadesuit (\N)}) = \rho_i$ for every $i =1, \cdots, r$.

The assignment $\msE^\spadesuit \mapsto \msE^{\spadesuit (\N)}$ determines  a morphism of  $\overline{\mcM}_{g, r}$-stacks
\begin{align} \label{eQ2}
\mcO p_{\N +1,\Rightarrow  \rho, g, r} \rightarrow \mcO p^{^\mr{Zzz...}}_{\N, \rho, g, r}.
\end{align}
Using this morphism, we occasionally consider  $\mcO p_{\N +1,\Rightarrow  \rho, g, r}$ as a sheaf on  $\mcO p^{^\mr{Zzz...}}_{\N, \rho, g, r}$ with respect to  the  fpqc topology.

\SSP
\bpr\label{Prop111}
Suppose that  $n=2$ and that $\mcO p_{\N, \rho, g, r}^{^\mr{Zzz...}}$ is nonempty and connected.
Then, 
the sheaf $\mcO p_{\N +1,\Rightarrow  \rho, g, r}$ on $\mcO p^{^\mr{Zzz...}}_{\N, \rho, g, r}$  locally admits  a section;
that is,  each geometric point of $\mcO p^{^\mr{Zzz...}}_{\N, \rho, g, r}$  is contained in an fpqc scheme $U$  over $\mcO p^{^\mr{Zzz...}}_{\N, \rho, g, r}$ such that  the sheaf $\mcO p_{\N +1,\Rightarrow  \rho, g, r}$  admits  a section over $U$.
\epr
\begin{proof}
It follows from Proposition \ref{Prop77} that there exists an $r$-tuple $\widetilde{\rho} \in \Xi_{2, \N +1}^{r\ma}$ with $\widetilde{\rho}^{(\N)} = \rho$ and $\mcO p_{2, \N +1, \widetilde{\rho}, g, r}^{^\mr{Zzz...}} \neq \emptyset$.
Since  both $\mcO p_{\N, \rho, g, r}^{^\mr{Zzz...}}$ and $\mcO p_{\N+1, \widetilde{\rho}, g, r}^{^\mr{Zzz...}}$ are smooth and of  the same dimension (cf. Theorem \ref{Th8}, (ii)),
the composite
\begin{align} \label{Er51}
\mcO p_{\N+1, \widetilde{\rho}, g, r}^{^\mr{Zzz...}} \xrightarrow{\mr{inclusion}} \mcO p_{\N +1,\Rightarrow  \rho, g, r} \xrightarrow{\eqref{eQ2}} \mcO p^{^\mr{Zzz...}}_{\N, \rho, g, r}
\end{align}
 is  flat (cf. ~\cite[Chapter III, Exercise 10.9]{Har}).
Moreover, 
this composite is proper because of the properness of $\mcO p_{\N+1, \widetilde{\rho}, g, r}^{^\mr{Zzz...}}$ over $k$, so 
the connectedness of $\mcO p_{\N, \rho, g, r}^{^\mr{Zzz...}}$ implies that
the composite \eqref{Er51} is surjective.
In particular, when 
regarded as a sheaf on  $\mcO p^{^\mr{Zzz...}}_{\N, \rho, g, r}$ with respect to the fpqc topology,
  $\mcO p_{\N+1, \widetilde{\rho}, g, r}^{^\mr{Zzz...}}$  admits locally a section.
The same holds for 
the sheaf $\mcO p_{\N +1,\Rightarrow  \rho, g, r}$, as it  contains $\mcO p_{\N+1, \widetilde{\rho}, g, r}^{^\mr{Zzz...}} $.
This completes the proof of this assertion.
\end{proof}
\SSP

Given a  dormant $\mr{PGL}_n^{(\N)}$-oper  $\msE^\spadesuit := (\mcE_B, \phi)$ on $\msX$ of radii $\rho$,
we set
\begin{align} \label{Er21}
V_{\msE^\spadesuit} := \mbV (f_* (\Omega^{(\N)} \otimes \mcS ol (\nabla^\mr{ad}_\phi))),
\end{align}
i.e., the total space of  the vector bundle  $f_* (\Omega^{(\N)} \otimes \mcS ol (\nabla^\mr{ad}_\phi))$.
The $S$-schemes  $V_{\msE^\spadesuit}$ defined for the various pairs $(\msX, \msE^\spadesuit)$ can be glued together to form a stack
\begin{align} \label{Er20}
V_\mr{univ}
\end{align}
over $\mcO p^{^\mr{Zzz...}}_{\N, \rho, g, r}$.
That is to say, if $q$ denotes the morphism $S \rightarrow \mcO p^{^\mr{Zzz...}}_{\N, \rho, g, r}$ classifying $\msE^\spadesuit$,
then there exists a functorial isomorphism of $S$-schemes
\begin{align} \label{eQ5}
V_{\msE^\spadesuit} \xrightarrow{\sim} V_\mr{univ} \times_{\mcO p^{^\mr{Zzz...}}_{\N, \rho, g, r}, q} S.
\end{align}
In particular,  it follows from Corollary \ref{Co22} that $V_\mr{univ}$ forms a relative affine space over $\mcO p^{^\mr{Zzz...}}_{\N, \rho, g, r}$ whose relative dimension is  $(n^2 - 1) (g-1) + \frac{1}{2} (n^2 +n-2)r$.

\SSP
\bpr \label{Prop91}
Suppose  that $\mcO p^{^\mr{Zzz...}}_{\N, \rho, g, r} \neq \emptyset$ and that
the fpqc sheaf $\mcO p_{\N +1,\Rightarrow  \rho, g, r}$ on $\mcO p^{^\mr{Zzz...}}_{\N, \rho, g, r}$ locally admits   a section.
Then,
the stack $\mcO p_{\N +1, \Rightarrow  \rho, g, r}$ has a structure of relative affine space over $\mcO p^{^\mr{Zzz...}}_{\N, \rho, g, r}$ modeled on $V_\mr{univ}$.
In particular, $\mcO p_{\N +1, \Rightarrow  \rho, g, r}$ is a nonempty smooth Deligne-Mumford stack over $k$ of dimension $(n^2 +2)(g-1) + \frac{1}{2}(n^2 +n)r$.
\epr
\begin{proof}
Let $S$, $\msX$, $\msE^\spadesuit$, and  $q$ be as above, and choose a dormant $n^{(\N +1)}$-theta characteristic $\widetilde{\vartheta}$ of $\msX$.
Denote by $\theta$ the $n^{(\N)}$-theta characteristic induced by $\widetilde{\vartheta}$ via level reduction, and by $\nabla^\diamondsuit$ the $(\mr{GL}_n^{(\N)}, \vartheta)$-oper corresponding to $\msE^\spadesuit$.

Now, suppose that we are given a $\mr{PGL}_n^{(\N +1)}$-oper $\widetilde{\msE}^\spadesuit$ on $\msX$ inducing $\msE^\spadesuit$ via level reduction.
This corresponds to a dormant $(\mr{GL}_n^{(\N +1)}, \widetilde{\vartheta})$-oper $\widetilde{\nabla}^\diamondsuit$, and 
we obtain the composites
\begin{align}
\alpha (\nabla^\diamondsuit) &: \mcD_{\leq p^\N -1}^{(\N -1)} \otimes \mcF_\varTheta \xrightarrow{\mr{inclusion}} \mcD^{(\N -1)} \otimes \mcF_\varTheta \xrightarrow{\nabla^\diamondsuit} \mcF_\varTheta,  \\
 \alpha (\widetilde{\nabla}^\diamondsuit) &: \mcD^{(\N)}_{\leq p^\N} \otimes \mcF_\varTheta \xrightarrow{\mr{inclusion}} \mcD^{(\N)} \otimes \mcF_\varTheta \xrightarrow{\widetilde{\nabla}^\diamondsuit} \mcF_\varTheta. \notag
\end{align}
The assignment $\widetilde{\msE}^\spadesuit \mapsto \alpha (\widetilde{\nabla}^\diamondsuit)$ determines a map of sets
$\mcO p_{q} (S) \rightarrow \mr{Hom}_{\mcO_X} (\mcD_{\leq p^\N}^{(\N)}\otimes \mcF_\varTheta, \mcF_\varTheta)$,
where $\mcO p_{q} := \mcO p_{\N +1, \Rightarrow \rho, g, r} \times_{\mcO p^{^\mr{Zzz...}}_{\N, \rho, g, r}, q} S$.
The formation of this  map
 is functorial with respect to $S$, so it gives rise to  a morphism of $S$-schemes
\begin{align}
\mcO p_{q} \rightarrow  \mbV (f_* (\mcH om_{\mcO_X} (\mcD^{(\N)}_{\leq p^\N} \otimes \mcF_\varTheta, \mcF_\varTheta))).
\end{align}
Since $\mcD^{(\N)}$ is generated by the sections of $\mcD_{\leq p^\N}^{(\N)}$ as a $\mcD^{(\N)}$-module,
this morphism allows us to view 
$\mcO p_{q}$  as a closed subscheme of $\mbV (f_* (\mcH om_{\mcO_X} (\mcD^{(\N)}_{\leq p^\N} \otimes \mcF_\varTheta, \mcF_\varTheta)))$.

Next, the natural injection $\mcD_{\leq p^\N -1}^{(\N -1)} \rightarrow \mcD_{\leq p^\N}^{(\N)}$ induces a morphism between total spaces of vector bundles
\begin{align}
\gamma :  \mbV (f_* (\mcH om_{\mcO_X} (\mcD^{(\N)}_{\leq p^\N} \otimes \mcF_\varTheta, \mcF_\varTheta))) \rightarrow 
 \mbV (f_* (\mcH om_{\mcO_X} (\mcD^{(\N-1)}_{\leq p^\N -1} \otimes \mcF_\varTheta, \mcF_\varTheta))).
\end{align}
Let  $\sigma$ denote the section $S \rightarrow  \mbV (f_* (\mcH om_{\mcO_X} (\mcD^{(\N-1)}_{\leq p^\N -1} \otimes \mcF_\varTheta, \mcF_\varTheta)))$ determined by $\alpha (\nabla^\diamondsuit)$.
Then, 
  $\mcO p_q$  is contained in  $\gamma^{-1} (\mr{Im} (\sigma))$.
Since $\mcD_{\leq p^\N}^{(\N)}/\mcD_{\leq p^\N -1}^{(\N -1)} \cong \mcT^{\otimes p^\N}$,
it follows from our  assumption  that  $\gamma^{-1} (\mr{Im} (\sigma))$
forms a relative affine space over $S$ modeled  on $\mbV (f_* (\mcH om_{\mcO_X} (\mcT^{\otimes p^\N} \otimes \mcF_\varTheta, \mcF_\varTheta))) \left(= \mbV (f_* (\Omega^{\otimes p^\N} \otimes \mcE nd (\mcF_\varTheta))) \right)$.

Moreover, 
by Lemma \ref{Lem42} proved below, $\mcO p_q$ forms an affine subspace of $\gamma^{-1} (\mr{Im} (\sigma))$, modeled on $\mbV (f_*^{(\N)} (\Omega^{(\N)} \otimes \mcS ol (\nabla^\diamondsuit_{\mcE nd^0})))$,
where $\nabla^\diamondsuit_{\mcE nd^0}$ denotes the $\mcD^{(\N -1)}$-module structure on the sheaf  $\mcE nd^0 (\mcF_\varTheta)$ of $\mcO_{X}$-linear endomorphisms of $\mcF_\varTheta$ with vanishing trace induced from $\nabla^\diamondsuit$.
Since  $(\mfg_\mcE, \nabla^\mr{ad}_\phi)$ can be naturally identified with $(\mcE nd^0(\mcF_\varTheta), \nabla^\diamondsuit_{\mcE nd^0})$,
 we see that  $\mcO p_q$ forms a relative affine space modeled on $V_{\msE^\spadesuit}$.
 Furthermore, it is immediately verified that this affine structure  glues compatibly across $S$, yielding a global structure of  relative affine space  on $\mcO p_{\N +1, \Rightarrow \rho, g, r}^{^\mr{Zzz...}}$ modeled on $V_\mr{univ}$ via \eqref{eQ5}.
 This completes the proof of the assertion.
\end{proof}
\SSP

The following lemma was used in the proof of the above proposition.

\SSP
\ble \label{Lem42}
Let us retain  the notation from  the proof of Proposition \ref{Prop91}.
Let 
$\nu$ be an $\mcO_X$-linear morphism $\mcT^{\otimes p^\N} \rightarrow \mcE nd (\mcF_\varTheta)$ (or equivalently, a global section of $\Omega^{\otimes p^\N} \otimes \mcE nd (\mcF_\varTheta)$), and regard it  as an $\mcO_X$-linear morphism $\mcD_{\leq p^\N}^{(\N)} \otimes \mcF_\varTheta \rightarrow \mcF_\varTheta$ via the natural isomorphism
$\mcD^{(\N)}_{\leq p^\N}/\mcD_{\leq p^\N -1}^{(\N)}\xrightarrow{\sim} \mcT^{\otimes p^\N}$.
Then, $\widetilde{\nabla}^\diamondsuit + \nu = \alpha (\widetilde{\nabla}_1^\diamondsuit)$ for some  $(\mr{GL}_n^{(\N +1)}, \widetilde{\vartheta})$-oper $\widetilde{\nabla}_1^\diamondsuit$ if and only if $\nu \in H^0 (X^{(\N)}, \Omega^{(\N)} \otimes \mcS ol (\nabla^\diamondsuit_{\mcE nd^0}))$.
Moreover, if these equivalent conditions are fulfilled, then the $(\mr{GL}_{n}^{(\N +1)}, \widetilde{\vartheta})$-oper ``$\widetilde{\nabla}_1^\diamondsuit$" is uniquely determined.
\ele
\begin{proof}
Suppose that $\widetilde{\nabla}^\diamondsuit + \nu = \alpha (\widetilde{\nabla}^\diamondsuit_1)$ for some  $(\mr{GL}_n^{(\N +1)}, \widetilde{\vartheta})$-oper $\widetilde{\nabla}_1$.
Then, by the commutativity relation \eqref{eQ10}, $\nu$ must lie  in  $H^0 (X^{(\N)}, \Omega^{(\N)} \otimes \mcS ol (\nabla_{\mcE nd^0}))$.
Moreover, by reversing this discussion, we see that the inverse direction of the desired equivalence can be verified.
\end{proof}

\LSP
\subsection{An affine space structure on $\mcO p_{\N +1, \widetilde{\rho}, g, r}$} \label{SS19e9}

Denote by ${^\dagger}\mfc_{\N+1}$ the  subscheme of $\mfc_{\N+1}$ defined to be the scheme-theoretic image  of the closed  subscheme  
\begin{align}
\left\{ (v_1, \cdots, v_{\N+1}) \in \mfg^{\N+1} \, | \, [v_j, v_{j'}] = 0 \ \text{if} \ j \neq j' \right\}
\end{align}
 under the natural quotient  $\mfg^{\N +1} \twoheadrightarrow \mfc_{\N+1}$.
We shall set ${^\dagger}\Xi_{n, \N+1}$ to be the $k$-scheme fitting into the following Cartesian square diagram:
\begin{align} \label{Er501}
\vcenter{\xymatrix@C=46pt@R=36pt{
{^\dagger}\Xi_{n, \N +1} \ar[r] \ar[d]  &{^\dagger}\mfc_{\N +1} \ar[d] \\
\Xi_{n, \N}\ar[r] & \mfc_\N,  
 }}
\end{align}
where the bottom horizontal arrow arises from \eqref{Er4} and the right-hand vertical arrow is  induced by  the projection onto the first $\N$ factors $\mfg^{\N +1} \twoheadrightarrow \mfg^{\N}$.
We denote the left-hand vertical arrow by 
\begin{align} \label{Et7}
{^\dagger}\pi : {^\dagger}\Xi_{n, \N +1} \twoheadrightarrow \Xi_{n, \N}.
\end{align}
In particular, for each element  $\rho_0 \in \Xi_{n, \N}$ represented by a multiset $\widetilde{\rho}_0 \in \mfS_n \backslash (\mbZ/p^\N \mbZ)^{n\ma}$, the inverse image ${^\dagger}\pi^{-1} (\rho_0)$ can be identified, after choosing  an ordering  of the elements in $\widetilde{\rho}_0$,  with the Lie algebra $\mft$ of the  maximal torus of $\mr{PGL}_n$ consisting of  diagonal matrices.

Let $\rho := (\rho_i)_{i=1}^r$ be an element of $\Xi^{r\ma}_{n, \N}$ satisfying $\mcO p^{^\mr{Zzz...}}_{\N, \rho, g, r} \neq \emptyset$, and write  ${^\dagger}\pi^{-1}(\rho) := ({^\dagger}\pi^{-1}(\rho_i))_{i=1}^r$.
Since the constituents  in  the monodromy operators of a $\mcD^{(\N)}$-module structure commute pairwise,
the radii $\rho_i (\msE^\spadesuit)$ of a  $\mr{PGL}_{n}^{(\N +1)}$-oper  $\msE^\spadesuit$ defines a  point of ${^\dagger}\mfc_\N$.
Moreover, if $\msE^\spadesuit$ is classified by $\mcO p_{\N +1, \Rightarrow \rho, g, r}$,
then
the $r$-tuple  $\rho (\msE^\spadesuit) := (\rho_i (\msE^\spadesuit))_{i=1}^r$ lies in ${^\dagger}\pi^{-1} (\rho)$.
The assignment $\msE^\spadesuit \mapsto \rho (\msE^\spadesuit)$ determines a morphism of $k$-stacks
\begin{align} \label{Et10}
\mr{Rad} : \mcO p_{\N +1, \Rightarrow \rho, g, r} \rightarrow {^\dagger}\pi^{-1} (\rho).
\end{align}
For $\widetilde{\rho} \in {^\dagger}\pi^{-1}(\rho) (k)$,
the inverse image $\mr{Rad}^{-1}(\widetilde{\rho})$ can be identified with $\mcO p_{\N +1, \widetilde{\rho}, g, r}$.

Now, fix  a $k$-rational  point $q \in \mcO p^{^\mr{Zzz...}}_{\N, \rho, g, r} (k)$.
Then, $\mr{Rad}$ restricts to a morphism
\begin{align} \label{Er552}
\mr{Rad} |_{q} : \mcO p_{\N +1, \Rightarrow \rho, g, r} |_q \rightarrow {^\dagger}\pi^{-1} (\rho), 
\end{align}
where $\mcO p_{\N +1, \Rightarrow \rho, g, r} |_q := \mcO p_{\N +1, \Rightarrow \rho, g, r} \times_{ \mcO p^{^\mr{Zzz...}}_{\N, \rho, g, r} , q} \mr{Spec} (k)$.
In what follows, we examine an affine structure on each fiber of $\mr{Rad}|_q$.

Let $S$, $\msX$ be  as before, and
suppose that $S = \mr{Spec}(k)$ and the point $q$ lies over the point of $\overline{\mcM}_{g, r}$ classifying $\msX$.
Also, choose  a dormant $n^{(\N)}$-theta characteristic $\vartheta := (\varTheta, \nabla_\vartheta)$ of $\msX$.
The point $q$ determines  a dormant $\mr{PGL}_n^{(\N)}$-oper $\msE^\spadesuit$ on $\msX$, which corresponds to
a dormant  $(\mr{GL}_n^{(\N)}, \vartheta)$-oper $\nabla^\diamondsuit$.
Let us fix $i \in \{ 1, \cdots, r \}$ and fix an identification of the formal neighborhood $U_i$ of $\mr{Im}(\sigma_i)\subseteq X$ with the formal disk  $U_\oslash$ over $k$.
Then, there exists a representative $[d_{i, 1}, \cdots, d_{i, n}] \in \mfS_n \backslash (\mbZ/p^\N \mbZ)^{n\ma}$ of $\rho_i$ such that 
$0 \leq \widetilde{d}_{i, 1} < \widetilde{d}_{i, 2} < \cdots < \widetilde{d}_{i, n} < p^\N$ and 
we have an isomorphism
 of $\mcD^{(\N -1)}_\oslash$-modules $(\mcF_\varTheta, \nabla^\diamondsuit)|_{U_i} \xrightarrow{\sim} \bigoplus_{j=1}^n \msO^{(\N -1)}_{\oslash, d_j}$.
Similarly to the discussion  surrounding  \eqref{Et66},
the restriction of   $(\mfg_\mcE, \nabla_\phi^\mr{ad})$ to $U_i$ admits an isomorphism
\begin{align}  \label{Ey221}
(\mfg_\mcE, \nabla_\phi^\mr{ad})|_{U_i} 
\xrightarrow{\sim} \mft_{} \otimes_k \msO^{(\N -1)}_{\oslash, 0} \oplus \bigoplus_{1 \leq j <  j' \leq n} \msO^{(\N -1)}_{\oslash, d_{i, j} - d_{i, j'}} \oplus \bigoplus_{1 \leq j' <  j \leq n} \msO^{(\N -1)}_{\oslash, d_{i, j} - d_{i, j'}}, 
\end{align}
since $\nabla^\mr{ad}_\phi$ can be identified with the $\mcD^{(\N -1)}$-module structure on $\overline{\mcE nd} (\mcF_\varTheta)$ induced naturally from 
$\nabla^\diamondsuit$.

We now denote by   $(\mcV, \nabla_\mcV)$ the codomain of this isomorphism.
The inclusion into the first summand   $\mft \otimes \mcO_\oslash \hookrightarrow \mcV$ induces a $k$-linear injection $\iota : \sigma_{i*}(\mft) \hookrightarrow \sigma_{i*} (\sigma^*_i (\mcV))$.
Since $\msO^{(\N -1)}_{\oslash, d_{i, j} - d_{i, j'}} \not\cong \msO^{(\N -1)}_{\oslash, 0}$ for $j \neq j'$,
it follows from Lemma \ref{Lem4490} that the image of
the composite
 $\mcS ol (\nabla_\phi^{\mr{ad}}) \hookrightarrow \mcV \twoheadrightarrow \sigma_{i*} (\sigma^*_i (\mcV))$
  coincides with $\mr{Im}(\iota)$.
It follows that \eqref{Ey221} induces a surjection  $\xi_i : \mcS ol (\nabla_\phi^\mr{ad}) \twoheadrightarrow \left(\sigma_{i*}(\mft) =\right) \sigma_{i*}^{(\N)} (\mft)$, and we obtain a composite
\begin{align}
\xi'_i : \Omega^{(\N)} \otimes\mcS ol (\nabla_\phi^\mr{ad}) \xrightarrow{\mr{id}_{\Omega^{(\N)}}\otimes \xi_i} \Omega^{(\N)} \otimes \sigma^{(\N)}_{i*} (\mft) \xrightarrow{\sim}
\sigma^{(\N)}_{i*} (\sigma^{(\N)*}_i(\Omega^{(\N)}) \otimes_k \mft)
\xrightarrow{\sim} \sigma^{(\N)}_{i*} (\mft),
\end{align}
where the final  arrow arises  from the residue isomorphism $\sigma^{(\N)*}_i (\Omega^{(\N)}) \xrightarrow{\sim} k$.

By setting  $\mcG := \bigcap_{i=1}^r \mr{Ker} (\xi'_i)$, we obtain the following  short exact sequence of $\mcO_{X^{(\N)}}$-modules
\begin{align} \label{Er444}
0 \rightarrow \mcG \rightarrow \Omega^{(\N)} \otimes  \mcS ol (\nabla_\phi^\mr{ad}) \rightarrow \bigoplus_{i=1}^r \sigma_{i*}^{(\N)} (\mft) \rightarrow 0.
\end{align}
Note that both the sheaf $\mcG$ and 
the second arrow  $\mcG \rightarrow \Omega^{(\N)} \otimes  \mcS ol (\nabla_\phi^\mr{ad})$ in this sequence do not depend on the choice of the isomorphism $(\mcF_\varTheta, \nabla^\diamondsuit)|_{U_i} \xrightarrow{\sim} \bigoplus_{j=1}^n \msO^{(\N -1)}_{\oslash, d_j}$.

\SSP
\ble \label{Lem88}
The sequence
\begin{align} \label{Et5}
0 \rightarrow H^0 (X^{(\N)}, \mcG) \rightarrow H^0 (X^{(\N)}, \Omega^{(\N)}\otimes \mcS ol (\nabla_\phi^\mr{ad}))
\rightarrow  \mft^{r} \rightarrow 0
\end{align}
obtained from   \eqref{Er444}  by applying the $\delta$-functor $H^* (X^{(\N)}, -)$ (under the natural identification $H^0 (X^{(\N)}, \sigma^{(\N)}_{i*}(\mft)) = \mft$) is exact.
In particular, we have
\begin{align}
\mr{dim}_k (H^0 (X^{(\N)}, \mcG)) = (n^2 -1)(g-1) + \frac{1}{2} n(n-1) r.
\end{align}
\ele
\begin{proof}
We begin by establishing the exactness of \eqref{Et5}.
Let us use  the notation ``${^c} (-)$" as in Proposition \ref{Pee4}, (ii).
The Cartier operator of $({^c}\mfg_\mcE^{[\N -1]}, {^c}\nabla_\phi^{\mr{ad}[\N -1]})$ induces an isomorphism  
$\mr{Coker}({^c}\nabla_\phi^{\mr{ad}[\N -1]}) \xrightarrow{\sim}\Omega^{(\N)} \otimes \mr{Ker}({^c}\nabla_\phi^{\mr{ad}[\N -1]})$
(cf. ~\cite[Theorem 3.1.1]{Ogu2} or the discussion following ~\cite[Chapter I, Proposition 1.2.4]{Ogu1}).
It follows that 
 \begin{align} \label{Et11}
H^1 (X^{(\N)}, \Omega^{(\N)} \otimes \mr{Ker}({^c}\nabla_\phi^{\mr{ad}[\N -1]})) \cong  H^1 (X^{(\N)}, \mr{Coker}({^c}\nabla_\phi^{\mr{ad}[\N -1]}))  = 0,
\end{align}
 where the second equality follows from
Proposition \ref{Pee4}, (ii).
On the other hand,
Lemma \ref{Lem4490} in the case  $\M = 0$
implies 
$\mr{Ker}({^c}\nabla_\phi^{\mr{ad}[\N -1]}) \subseteq \bigcap_{i=1}^r\mr{Ker} (\xi_i)$.
This inclusion yields
 an injection $\Omega^{(\N)} \otimes \mr{Ker}({^c}\nabla_\phi^{\mr{ad}[\N -1]}) \hookrightarrow \mcG$, which is an isomorphism over $X^{(\N)} \setminus \bigcup_{i=1}^r \mr{Im}(\sigma^{(\N)}_i)$.
Hence,  \eqref{Et11} implies 
 $H^1 (X^{(\N)}, \mcG) = 0$, and hence the sequence  \eqref{Et5} is exact.
 The second assertion follows from the exactness of \eqref{Et5} and  Corollary \ref{Co22}.
\end{proof}
\SSP

\bpr \label{Prop400}
Suppose that $\mcO p^{^\mr{Zzz...}}_{\N, \rho, g, r} \neq \emptyset$ and that the fpqc sheaf $\mcO p_{\N +1, \Rightarrow \rho, g, r}$ on $\mcO p^{^\mr{Zzz...}}_{\N, \rho, g, r}$ admits locally a section.
Let us fix $\widetilde{\rho} \in {^\dagger}\pi^{-1} (\rho)$, and write
$\mcO p_{\N +1,  \widetilde{\rho}, g, r} |_q := \mcO p_{\N+1, \widetilde{\rho}, g, r} \times_{\mcO p^{^\mr{Zzz...}}_{\N, \rho, g, r}, q} \mr{Spec} (k)$.
Then, 
the affine structure on $\mcO p_{\N +1, \Rightarrow \rho, g, r} |_q$ (modeled on the $k$-vector space $H^0 (X^{(\N)}, \Omega^{(\N)}\otimes \mcS ol (\nabla_\phi^\mr{ad}))$) resulting from Proposition \ref{Prop91} 
restricts to  
an affine structure on $\mcO p_{\N +1,  \widetilde{\rho}, g, r} |_q$ modelded on $H^0 (X^{(\N)}, \mcG)$ via the inclusion $H^0 (X^{(\N)}, \mcG) \hookrightarrow H^0 (X^{(\N)}, \Omega^{(\N)}\otimes \mcS ol (\nabla_\phi^\mr{ad}))$. 
\epr
\begin{proof}
The assertion follows from the various definitions involved.
\end{proof}

\SSP
\bpr \label{Prop667}
Let us keep the assumption stated  in Proposition \ref{Prop400}.
Then, 
the morphism 
\begin{align} \label{Eii118}
\mcO p_{\N +1, \Rightarrow \rho, g, r} \rightarrow \mcO p^{^\mr{Zzz...}}_{\N, \rho, g, r} \times_k {^\dagger}\pi^{-1}(\rho)
\end{align}
induced by \eqref{eQ2} and $\mr{Rad}$ is smooth, surjective, and of relative  dimension  $(n^2 -1) (g-1) +\frac{1}{2} n(n-1)r$.
In particular, for each $\widetilde{\rho} \in {^\dagger}\pi^{-1}(\rho)$,
the projection $\mcO p_{\N +1, \widetilde{\rho}, g, r} \rightarrow \mcO p^{^\mr{Zzz...}}_{\N, \rho, g, r}$ is smooth, surjective, and of relative  dimension  $(n^2 -1) (g-1) +\frac{1}{2} n(n-1)r$.
\epr
\begin{proof}
Let $q$ be as before.
Suppose, on the contrary, that 
$\mr{Rad} |_{q}$ is not surjective.
Since ${^\dagger}\pi^{-1} (\rho)$ is irreducible,
the scheme-theoretic image of $\mr{Rad} |_{q}$
must be of strictly smaller dimension.
It follows  that the dimension of any nonempty fiber of $\mr{Rad} |_{q}$ must have dimension strictly  greater than 
\begin{align}
\mr{dim} (\mcO p_{\N +1, \Rightarrow \rho, g, r} |_q) - \mr{dim} ({^\dagger}\pi^{-1} (\rho))
&= (n^2 -1) (g-1) + \frac{1}{2}(n^2 +n-2)r - (n-1)r \\
& = (n^2 -1) (g-1) +\frac{1}{2} n(n-1)r
\end{align}
 (cf. ~\cite[Chapter II, Exercise 3.22]{Har}), where the first equality  follows from Proposition \ref{Prop91}.
 On the other hand,  the  fiber of  $\mr{Rad} |_{q}$ over $\widetilde{\rho} \in {^\dagger}\pi^{-1} (\rho)$ is isomorphic to
 $\mcO p_{\N +1, \widetilde{\rho}, g, r}|_q$, 
 whose dimension is computed as
 \begin{align}
 \mr{dim} (\mcO p_{\N +1, \widetilde{\rho}, g, r}|_q) =\mr{dim} (H^0 (X^{(\N)}, \mcG)) = (n^2-1)(g-1)+ \frac{1}{2} n(n-1)r,
 \end{align}
 where the first equality follows from Proposition \ref{Prop400} and the second equality follows from   Lemma \ref{Lem88}.
 Thus, we obtain a contradiction, which implies the surjectivity of $\mr{Rad} |_{q}$.
 
 By applying  the above discussion  to every geometric point $q \in \mcO p^{^\mr{Zzz...}}_{\N, \rho, g, r}(k)$, we conclude that the morphism  \eqref{Eii118} is surjective and  that all its  fibers  are affine spaces of dimension $(n^2-1)(g-1)+ \frac{1}{2} n(n-1)r$.
 Hence, the required assertion  follows from the smoothness of 
 $\mcO p_{\N +1, \Rightarrow \rho, g, r}$ and ${^\dagger}\pi^{-1}(\rho)$ (cf. Proposition \ref{Prop91} and ~\cite[Chapter III, Exercise 10.9]{Har}).
 \end{proof}
\SSP

\bco \label{Cor559}
Assume the hypothesis of 
 Proposition \ref{Prop400}.
Then, for each $\widetilde{\rho} \in \Xi_{2, \N}^{3\ma}$,
the stack $\mcO p_{\N, \widetilde{\rho}, 0, 3}$ is  isomorphic to $\mr{Spec}(k)$.
\eco
\begin{proof}
The assertion follows from Propositions \ref{Prop400} and \ref{Prop667}, since  the equality $H^0 (X^{(\N)}, \mcG) = 0$ holds when $(g, r, 2) = (0, 3, 2)$.
\end{proof}

\LSP
\subsection{The Hitchin-Mochizuki morphism of higher level} \label{SS199}

Recall  the $S$-scheme  $B_\msX^\flat$ (resp., $B^\flat_{\msX, 0}$)  introduced in  Section \ref{SS5r1}.
Under the situation where  $\M = \N$,  these schemes  defined for the various $\msX$'s can be glued together to form  a stack
\begin{align} \label{Ey34}
B_\mr{univ}^\flat \ \left(\text{resp.,} \ B^\flat_{\mr{univ}, 0}\right)
\end{align}
over $\overline{\mcM}_{g, r}$.
That is to say,  if $\overline{q}$ denotes the morphism $S \rightarrow \overline{\mcM}_{g, r}$ classifying $\msX$,
then we have $B_\mr{univ}^\flat \times_{\overline{\mcM}_{g, r}, \overline{q}} S \cong B^\flat_\msX$ (resp.,  $B^\flat_{\mr{univ}, 0} \times_{\overline{\mcM}_{g, r}, \overline{q}} S \cong B_{\msX, 0}^\flat$).
By 
\eqref{Eu3}  (resp., \eqref{Eu4}),
the stack $B_\mr{univ}^\flat$ (resp., $B_{\mr{univ}, 0}^\flat$)  forms a relative affine space over $\overline{\mcM}_{g, r}$ whose relative dimension equals $(n^2 -1)(g-1) + \frac{1}{2}(n^2 + n-2)r$ (resp., $(n^2 -1)(g-1)+ \frac{1}{2}n(n-1)r$).

For each $\rho \in \mfc_\N (k)^{r\ma}$ (where $\rho = \emptyset$ if $r = 0$),
 the  assignment  $\msE^\spadesuit \mapsto \psi\text{-}\mr{Char}(\msE^\spadesuit)$
(cf. \eqref{Er115}) yields a morphism 
\begin{align} \label{Et61}
\mr{HM}^{(\N)} : \mcO p_{\N +1, \Rightarrow \rho, g, r} \rightarrow B^\flat_\mr{univ},
\end{align}
 which makes  the following square diagram  commute:
\begin{align} \label{EQ145}
\vcenter{\xymatrix@C=46pt@R=36pt{
 \mcO p_{\N +1, \Rightarrow \rho, g, r} \ar[r]^-{\mr{HM}^{(\N)}} \ar[d]_-{\eqref{eQ2}}& B^\flat_\mr{univ} \ar[d]^-{\mr{projection}} \\
 \mcO p^{^\mr{Zzz...}}_{\N, \rho, g, r} \ar[r]_-{\Pi'_{\N, \rho, g, r}} & \overline{\mcM}_{g, r}.
 }}
\end{align}
The following assertion generalizes ~\cite[Theorem A.1]{BeTr},  ~\cite[Theorem 6.1.3]{JoPa}, ~\cite[Chapter II, Theorem 2.3]{Moc1}, and ~\cite[Theorem 3.39]{Wak5}.

\SSP
\bt\label{Cor20}
Suppose that $\mcO p^{^\mr{Zzz...}}_{\N, \rho, g, r} \neq \emptyset$ and that  the fpqc sheaf $\mcO p_{\N + 1, \Rightarrow \rho, g, r}$ on $\mcO p_{\N, \rho, g, r}^{^\mr{Zzz...}}$ admits locally a section.
Then, the following assertions hold:
\begin{itemize}
\item[(i)]
The   morphism
\begin{align} \label{Er207}
\mcO p_{\N +1, \Rightarrow \rho, g, r} \rightarrow \mcO p^{^\mr{Zzz...}}_{\N, \rho, g, r} \times_{\overline{\mcM}_{g, r}} B_\mr{univ}^\flat
\end{align}
induced by  \eqref{EQ145}
 is finite and faithfully flat. 
 \item[(ii)]
 Let $\widetilde{\rho}$ be an element of ${^\dagger}\pi^{-1}(\rho) \cap \Xi_{n, \N+1}^{r\ma}$.
 (By Proposition \ref{Prop667}, the stack $\mcO p_{\N +1, \widetilde{\rho}, g, r}$ is nonempty.)
 Then, $\mr{HM}^{(\N)}$ restricts to a morphism
 \begin{align} \label{Er247} 
\mr{HM}^{(\N)}_{\widetilde{\rho}} : \mcO p_{\N +1, \widetilde{\rho}, g, r} \rightarrow 
B_{\mr{univ}, 0}^\flat,
\end{align}
and the induced morphism
\begin{align} \label{Er247dd}
\mcO p_{\N +1, \widetilde{\rho}, g, r} \rightarrow 
\mcO p^{^\mr{Zzz...}}_{\N, \rho, g, r} \times_{\overline{\mcM}_{g, r}}  B_{\mr{univ}, 0}^\flat
\end{align}
 is finite and faithfully flat.
 \end{itemize}
\et
\begin{proof}
Let us consider assertion (i).
To begin with, we shall prove
the properness of \eqref{Er207}  restricted over each geometric point of $ \mcO p^{^\mr{Zzz...}}_{\N, \rho, g, r}$, by applying the standard   valuative criterion.
Let $\msX$ be as above, and suppose that $S = \mr{Spec} (k)$.
Also, let $R$ be a valuation ring of a field $K$ over $k$ whose residue field $R/\mfm$ is isomorphic to $k$.
Denote by  $\eta$ (resp., $s$) the generic (resp., closed) point of $T := \mr{Spec}(R)$.
For $\Box \in \{ T, \eta, s \}$, we denote by $\msX_\Box := (X_\Box, \{ \sigma_{\Box, i} \}_i)$ the base-change of $\msX$ over $\Box$ (in particular,  $\msX_s = \msX$).
The sheaf of  logarithmic $1$-forms (resp., logarithmic tangent fields; resp., logarithmic differential operators of level $\N$) on the associated log curve $X_\Box^\mr{log}/\Box^\mr{log}$   is denoted by $\Omega_\Box$ (resp., $\mcT_\Box$; resp., $\mcD^{(\N)}_\Box = \bigcup_{j \in \mbZ_{\geq 0}} \mcD_{\leq j, \Box}^{(\N)}$).
According to  ~\cite[Proposition 5.14]{Wak20},
there exists  an $n^{(\N +1)}$-theta characteristic $\vartheta := (\varTheta, \nabla_\vartheta)$ of $X^\mr{log}/\mr{Spec}(k)^\mr{log}$.
In particular, we obtain  the base-change $\vartheta_\Box := (\varTheta_\Box, \nabla_{\vartheta, \Box})$ (resp., $\mcF_{\varTheta, \Box}$) of $\vartheta$ (resp.,  $\mcF_\varTheta$) over $\Box$.

Let us take a $T$-rational point of $B_\msX$, which corresponds, via \eqref{Er99},
to a polynomial 
$Q_T  := t^n + a_2 t^{n-2} + \cdots + a_n$ for some  $a_\ell \in H^0 (X_T^{(\N +1)}, (\Omega_{T}^{(\N+1)})^{\otimes \ell}) \subseteq H^0 (X_T, \Omega_T^{\otimes p^{\N +1} \cdot \ell})$ ($2 \leq \ell \leq n$).
Write $Q_\eta := Q_T |_{X_\eta}$,  $Q_s := Q_T |_{X_s}$, and
suppose that we are given 
a $\mr{PGL}_n^{(\N +1)}$-oper  $\msE^\spadesuit_\eta$
 on $\msX_\eta$ satisfying   $\psi\text{-}\mr{Char} (\msE^\spadesuit_\eta) = Q_\eta$.
 Let $\nabla_\eta^\diamondsuit$ denote the $(\mr{GL}_n^{(\N +1)}, \vartheta_\eta)$-oper 
  corresponding to $\msE_\eta^\spadesuit$.
To prove  the desired properness,  it suffices   to  construct a $\mr{PGL}_n^{(\N +1)}$-oper $\msE_T^\spadesuit$ on $\msX_T$ 
whose
 generic fiber is $\msE^\spadesuit_\eta$ and which satisfies  the equality $\psi\text{-}\mr{Char}(\msE_T^\spadesuit) = Q_T$.

For $\Box \in \{ T, \eta, s \}$,
let  $\langle Q_\Box \rangle$ be the $\mcD^{(\N)}_\Box$-submodule of $\mcD^{(\N)}_\Box\otimes \varTheta_{\Box}$ generated by the sections
\begin{align}
\left( (\partial^{\langle p^{\N +1} \rangle})^n + \langle a_2 |_{X_\Box}, (\partial^{\langle p^{\N +1}\rangle})^{\otimes 2} \rangle \cdot   (\partial^{\langle p^{\N +1}\rangle})^{n-2} 
+ \cdots +  \langle a_n |_{X_\Box}, (\partial^{\langle p^{\N +1}\rangle})^{\otimes n} \rangle
\right) \otimes v 
\end{align}
defined for various local sections $v \in \varTheta_{\Box}$, where $\langle -, - \rangle$'s denote the natural pairings $\Omega_\Box^{\otimes p^{\N +1} \cdot  \ell} \times \mcT^{\otimes p^{\N +1} \cdot  \ell}_\Box \rightarrow \mcO_{X_\Box}$, and $\partial^{\langle p^{\N +1} \rangle}$ denotes the local section of $\mcD^{(\N)}_\Box$ as introduced in Section \ref{SS11}.
The quotient $\mcP_\Box := (\mcD^{(\N)}_\Box \otimes \varTheta_{\Box})/\langle Q_\Box \rangle$ carries a natural  $\mcD^{(\N)}_\Box$-module  structure $\nabla_{\mcP_\Box}$.
For each $j=0, \cdots, p^{\N+1} \cdot n$, we shall set $\mcP_\Box^j$ to be the subbundle of $\mcP_\Box$ defined as 
\begin{align}
\mcP_\Box^j := \mr{Im} \left( \mcD^{(\N)}_{\leq p^{\N +1}\cdot n-j-1, \Box}\otimes \varTheta_\Box \xrightarrow{\mr{inclusion}} \mcD_\Box^{(\N)} \otimes \varTheta_\Box \xrightarrow{\mr{quotient}} \mcP_\Box\right),
\end{align}
where $\mcD^{(\N)}_{\leq -1, \Box} := 0$.
The composite of natural morphisms  $\mcF_{\varTheta, \Box} 
 \hookrightarrow \mcD_\Box^{(\N)} \otimes \varTheta_{\Box} \twoheadrightarrow \mcP_\Box$ is injective 
 and restricts, for each $j=0, \cdots, n$, to an isomorphism $\mcF_{\varTheta, \Box}^j \xrightarrow{\sim} \mcP_{\Box}^{p^{\N +1} \cdot n -n +j}$.
  We regard $\mcF_{\varTheta, \Box}$  as an $\mcO_{X_\Box}$-submodule of $\mcP_\Box$ via the composite injection $\mcF_{\varTheta, \Box} \hookrightarrow \mcP_\Box$.
In the case $\Box = \eta$, 
the composite of  the inclusion $\mcD_\eta^{(\N)} \otimes \varTheta_\eta \hookrightarrow \mcD_\eta^{(\N)} \otimes \mcF_{\varTheta, \eta}$ and 
 the morphism  $\mcD_\eta^{(\N)} \otimes \mcF_{\varTheta, \eta} \rightarrow \mcF_{\varTheta, \eta}$ arising from $\nabla_\eta^\diamondsuit$
  factors through 
the projection  $\mcD_\eta^{(\N)} \otimes \varTheta_{\eta} \twoheadrightarrow \mcP_{\eta}$, as the equality $\mr{Char}(\breve{\psi}(\nabla_\eta^\diamondsuit)) = Q_\eta$ holds.
As a result, we obtain  a surjection between $\mcD^{(\N)}$-modules $\varpi_\eta : (\mcP_\eta, \nabla_{\mcP_\eta}) \twoheadrightarrow (\mcF_{\varTheta, \eta}, \nabla_\eta^\diamondsuit)$.

By the properness of Quot schemes, $\varpi_\eta$ extends to a surjection $\varpi_T: \mcP_T \twoheadrightarrow \mcF'$ for some coherent $\mcO_{X_T}$-module $\mcF'$ that is  flat over $T$ and satisfies $\mr{rank} (\mcF') = \mr{rank}(\mcF_{\varTheta, \eta})$, $\mr{deg} (\mcF') = \mr{rdeg}(\mcF_{\varTheta, \eta})$.
Note that the composite
\begin{align} \label{Eww2}
h : \mcF_{\varTheta, T} \xrightarrow{\mr{inclusion}}  \mcP_T \xrightarrow{\varpi_T} \mcF'
\end{align}
 is injective, since
its generic fiber $h_\eta$  is an isomorphism and $\mcF_{\varTheta, T}$ is locally free.
The kernel  $\mr{Ker}(\varpi_\eta) \left(= \mr{Ker}(\varpi_T) |_{X_\eta} \right)$ of $\varpi_\eta$ is closed under $\nabla_{\mcP_\eta}$,
so
the $T$-flatness of $\mcF'$ implies that $\mr{Ker}(\varpi_T)$ must be closed under $\nabla_{\mcP_T}$.
Hence, there exists a $\mcD^{(\N)}$-module structure on $\mcF'$ that extends $\nabla_\eta^\diamondsuit$ and commutes with $\nabla_{\mcP_T}$ via $\varpi_T$.

By an argument entirely similar to the proof of (the first assertion in) ~\cite[Lemma 5.25]{Wak20}, the coherent sheaf $\mcF'$ is verified to be locally free.
It follows that 
 $\mr{Ker}(\varpi_T)$ defines a subbundle of $\mcP_T$ and  $\mr{Ker}(\varpi_T)|_{X_s}$ coincides, in $\mcP_s$,  with the kernel $ \mr{Ker}(\varpi_s)$ of the special fiber $\varpi_s : \mcP_s \twoheadrightarrow \mcF' |_{X_s}$ of $\varpi_T$.
For each $j=0, \cdots, p^{\N +1} \cdot n -n$,
we set $\mr{Ker} (\varpi_T)^j := \mr{Ker}(\varpi_T) \cap \mcP_{T}^j$.
Since the composite $\mr{Ker} (\varpi_T) |_{X_\eta} \hookrightarrow \mcP_{\eta} \twoheadrightarrow \mcP_\eta/\mcP_\eta^{p^{\N +1} \cdot n - n}$ is an isomorphism,
the identity $\mr{Ker}(\varpi_T)^{p^{\N+1} \cdot n -n} = 0$ holds.
The subquotient $\mr{Ker}(\varpi_T)^j/\mr{Ker}(\varpi_T)^{j+1}$ is contained in $\mcP_T^{j}/\mcP_T^{j+1}$,
so it is $T$-flat.
Hence, the collection $\{ \mr{Ker} (\varpi_T)^j |_{X_s}\}_{j=0}^{p^{\N +1} \cdot n - n}$ forms a decreasing filtration on $\mr{Ker}(\varpi_s) \left( = \mr{Ker}(\varpi_T) |_{X_s}\right)$ and, for each $j=0, \cdots, p^{\N +1} \cdot n - n -1$, we have
\begin{align} \label{Ejkw}
\mr{deg}\left( (\mr{Ker} (\varpi_T)^j |_{X_s})/(\mr{Ker} (\varpi_T)^{j+1} |_{X_s}) \right)
& = \mr{deg} \left((\mr{Ker} (\varpi_T)^j |_{X_\eta})/\mr{Ker} (\varpi_T)^j |_{X_\eta} \right) \\
& = \mr{deg}(\mcP_\eta^{j}/\mcP_\eta^{j+1}) \notag \\
& = \mr{deg} (\mcT^{\otimes (p^{\N +1} \cdot n- 1 -j)}_\eta \otimes \varTheta_\eta) \notag \\
& = \mr{deg} (\varTheta) - (2g-2+r) (p^{\N +1}\cdot n - 1 -j) \notag \\
& < \mr{deg}(\varTheta) - (2g-2+r) (n-1).
\end{align} 

Suppose now that the special fiber  $h_s : \mcF_{\varTheta, s} \rightarrow \mcF' |_{X_s}$  of $h$ is {\it not} injective.
Similarly to  the second assertion in  ~\cite[Lemma 5.25]{Wak20},
there exists  a line subbundle of $\mr{Ker} (h_s)$ whose degree is equal to or greater than the value $\mr{deg}(\mcF_{\varTheta, s}^0/\mcF_{\varTheta, s}^{1}) = \mr{deg}(\varTheta) - (2g-2+r) (n-1)$.
Note that $\mr{Ker}(h_s)$ can be identified with $\mr{Ker} (\varpi_s) \cap \mcP_{s}^{p^{\N +1} \cdot n -n}$ via the isomorphism ``$\mcF_{\varTheta, \Box}^j \xrightarrow{\sim} \mcP_{\Box}^{p^{\N +1} \cdot n -n +j}$" constructed above in the case where $(\Box, j)$ is taken to be $(s, 0)$.
Hence, we obtain a contradiction because of the computation in  the sequence \eqref{Ejkw}.
Thus, $h_s$ must be injective.
By comparing the degrees of $\mcF_{\varTheta, s}$ and $\mcF' |_{X_s}$, we see that $h_s$, and hence $h$, is an isomorphism.

Under the  identification  $\mcF_{\varTheta, T} = \mcF'$  given  by  the isomorphism $h$,
the $\mcD^{(\N)}_T$-module structure on $\mcF'$ constructed above determines 
a $\mcD^{(\N)}$-module structure $\nabla_T^\diamondsuit$ on $\mcF_{\varTheta, T}$.
It is immediately verified that $\nabla_T^\diamondsuit$
 forms a $(\mr{GL}_n^{(\N +1)}, \vartheta_T)$-oper on $\msX_T$.
Since $\varpi_T$ defines a surjection between $\mcD^{(\N)}$-modules  $(\mcP_T, \nabla_{\mcP_T}) \twoheadrightarrow (\mcF_{\varTheta, T}, \nabla_T^\diamondsuit)$,
  the equality  $\mr{Char} (\breve{\psi}(\nabla_T^\diamondsuit)) = Q_T$ holds.
 Hence, the $\mr{PGL}_n^{(\N +1)}$-oper corresponding to $\nabla_T^\diamondsuit$ satisfies the required conditions, and this proves the properness of  \eqref{Er207} over every geometric  point of $\mcO p_{\N, \rho, g, r}^{^\mr{Zzz...}}$.

We now recall that both  $\mcO p_{\N +1, \Rightarrow \rho, g, r}$ and $\mcO p^{^\mr{Zzz...}}_{\N, \rho, g, r} \times_{\overline{\mcM}_{g, r}} B_\mr{univ}^\flat$ are relative affines spaces (in particular, smooth) over $\mcO p^{^\mr{Zzz...}}_{\N, \rho, g, r}$  of the same relative dimension (cf. Proposition \ref{Prop91}).
Combining this fact  with the properness established above, 
we see that \eqref{Er207} is proper, and hence finite.
Also, it follows from  ~\cite[Chapter III, Exercise 10.9]{Har} that 
this morphism is also faithfully flat.
This completes the proof of assertion (i).

Next, we shall consider assertion (ii).
The first assertion follows from ~\cite[Eq.\,(345)]{Wak5} and the observation in the proof of  Proposition \ref{Prop8}.
The second assertion follows from assertion (i) and  the fact that both the domain and codomain of \eqref{Er247dd} are smooth over $\mcO p^{^\mr{Zzz...}}_{\N, \rho, g, r}$ of the same relative dimension (cf. Proposition \ref{Prop667}).
\end{proof}
\SSP

\bde \label{Defgg9}
The morphisms  $\mr{HM}^{(\N)}$ and $\mr{HM}^{(\N)}_{\widetilde{\rho}}$  introduced above are  called   the {\bf level-$\N$ Hitchin-Mochizuki morphism}
 (cf. ~\cite[Appendix A]{BeTr}, ~\cite[Section 3]{LasPa0}, ~\cite[Chapter II, Section 2]{Moc1}, ~\cite[Section 3]{Wak5} for the case $\N =1$).
\ede
\SSP

It follows from  Proposition \ref{Prop94} and  the definition of $p^{\N +1}$-nilpotency that
the inverse image $(\mr{HM}^{(\N)})^{-1} (0_\mr{univ})$ of the zero section $0_\mr{univ} : \overline{\mcM}_{g, r} \rightarrow \left(B_{\mr{univ}, 0}^\flat \subseteq  \right) B_\mr{univ}^\flat$ via $\mr{HM}^{(\N)}$ admits 
 a decomposition
\begin{align} \label{eQ19}
(\mr{HM}^{(\N)})^{-1} (0_\mr{univ}) = 
\coprod_{\widetilde{\rho}}  \mcO p^{^\mr{nilp}}_{\N +1, \widetilde{\rho}, g, r},
\end{align}
where the disjoint union runs over the elements $\widetilde{\rho}$ in $\Xi_{n, \N +1}^{r\ma}$ with $\widetilde{\rho}^{(\N)} = \rho$, i.e., $\widetilde{\rho} \in {^\dagger}\pi^{-1}(\rho) \cap \Xi_{n, \N}^{r\ma}$.
If we fix  such an element $\widetilde{\rho}$,
then 
the restriction of \eqref{eQ2} to $\mcO p_{\N+1, \widetilde{\rho}, g, r}^{^\mr{nilp}}$ defines a morphism 
\begin{align} \label{Et33}
\Lambda_{\widetilde{\rho}, g, r} : \mcO p_{\N+1, \widetilde{\rho}, g, r}^{^\mr{nilp}} \rightarrow \mcO p^{^\mr{Zzz...}}_{\N, \rho, g, r}
\end{align}
over $\overline{\mcM}_{g, r}$.
Since $\mcO p_{\N+1, \widetilde{\rho}, g, r}^{^\mr{nilp}}$ can be identified with $(\mr{HM}^{(\N)}_{\widetilde{\rho}})^{-1}(0_{\mr{univ}})$  (cf. \eqref{Er247}),
   Theorem \ref{Cor20}, (ii),  implies  the following assertion.

\SSP
\bt \label{Th665}
Suppose  that   $\mcO p_{\N, \rho, g, r}^{^\mr{Zzz...}}\neq \emptyset$ and that
the fpqc sheaf $\mcO p_{\N +1, \Rightarrow \rho, g, r}$ on $\mcO p_{\N, \rho, g, r}^{^\mr{Zzz...}}$ admits locally a section.
Then, 
for any $\widetilde{\rho} \in {^\dagger}\pi^{-1}(\rho) \cap \Xi_{n, \N}^{r\ma}$, 
the stack $\mcO p_{\N+1, \widetilde{\rho}, g, r}^{^\mr{nilp}} $ is nonempty, and 
the morphism $\Lambda_{\widetilde{\rho}, g, r}$
  is 
  finite and faithfully flat.
\et
\SSP

\bco \label{Cor15}
Suppose that $n =2$.
Then,  for any $\widetilde{\rho} \in \Xi_{2, \N +1}^{3\ma}$ with $\mcO p_{\N, \widetilde{\rho}^{(\N)}, 0, 3}^{^\mr{Zzz...}}\neq \emptyset$,
the stack $\mcO p_{\N+1,  \rho, 0, 3}^{^\mr{nilp}}$ is isomorphic to $\mr{Spec}(k)$.
\eco
\begin{proof}
It follows from   Theorem \ref{Th665} that $\mcO p_{\N+1,  \rho, 0, 3}^{^\mr{nilp}}$ is nonempty.
Hence, the assertion follows from  
 the fact that $\mcO p_{\N+1,  \rho, 0, 3}^{^\mr{nilp}}$ is a closed subscheme of the $k$-scheme $\mcO p_{\N+1,  \rho, 0, 3}$, which is isomorphic to $\mr{Spec}(k)$ by Corollary \ref{Cor559}.
\end{proof}

\vspace{10mm}
\section{The irreducibility of the moduli space of dormant $\mr{PGL}_2^{(\N)}$-opers} \label{S3}
\LSP

The goal  of this section is  to prove  Theorem \ref{ThA} (cf. Theorem  \ref{Cor49}).
Our argument is basically an extension of the proof given by S. Mochizuki for $\N =1$ (cf. ~\cite[Chapter II, Theorems 1.12 and  2.8]{Moc2}).
We begin by establishing the irreducibility  of $\mcO p_{\N, \rho, g, r}^{^\mr{Zzz...}}$ under the assumption that  $\mr{dim}(\mcO p_{\N, \rho, g, r}^{^\mr{Zzz...}})=1$  (cf. Lemma \ref{T01}).
Then, we proceed by induction, reducing the general case to lower-dimensional cases via an analysis of divisors on the moduli space defined as  the images of various gluing  morphisms (cf. Lemma \ref{Cor11}).

Throughout this section, we assume  that $n=2$.

\LSP
\subsection{Splittings induced from  dormant $\mr{PGL}_2^{(\N)}$-opers} \label{SS39}

Let $(g, r)$ be a pair of nonnegative integers with $2g-2+r > 0$,
 $S$ a scheme (or more generally,   a Deligne-Mumford stack) over $k$,  $\msX := (f: X \rightarrow S, \{ \sigma_i \}_{i=1}^r)$ an $r$-pointed stable curve of genus $g$ over  $S$, and 
   $\vartheta := (\varTheta, \nabla_\vartheta)$  
 a dormant $2^{(\N)}$-theta characteristic   of $\msX$.
 Also, let   us take 
 a dormant $(\mr{GL}_2^{(\N)}, \vartheta)$-oper $\nabla^\diamondsuit$ on $\msX$.
Denote by  $D$ the relative effective divisor on $X$ defined as the union of the $\sigma_i$'s.

Consider the composite
\begin{align} \label{Eq50}
\mcD_{\leq p^{\N -1}}^{(\N -1)} \otimes \varTheta \xrightarrow{\mr{inclusion}} \mcD_{}^{(\N -1)} \otimes \mcF_\varTheta \xrightarrow{\nabla^\diamondsuit} \mcF_\varTheta, 
\end{align}
where the first arrow arises from the inclusions $\mcD_{\leq p^{\N -1}}^{(\N -1)} \hookrightarrow \mcD^{(\N -1)}$ and $\varTheta \hookrightarrow \mcF_\varTheta$.
This composite becomes an isomorphism when restricted to $\mcF_\varTheta = \mcD_{\leq 1}^{(\N -1)} \otimes \varTheta \left(\subseteq  \mcD_{\leq p^{\N -1}}^{(\N -1)} \otimes \varTheta\right)$, 
and hence determines a split surjection 
\begin{align} \label{Eq110}
\delta (\nabla^\diamondsuit) : \mcD_{\leq p^{\N -1}}^{(\N -1)} \otimes \varTheta \twoheadrightarrow \mcF_\varTheta
\end{align}
of the natural short exact sequence
\begin{align}
0 \rightarrow 
\mcF_\varTheta
 \rightarrow \mcD_{\leq p^{\N -1}}^{(\N -1)} \otimes \varTheta \rightarrow (\mcD_{\leq p^{\N -1}}^{(\N -1)}/  \mcD_{\leq 1}^{(\N -1)})\otimes \varTheta \rightarrow 0.
\end{align}

\SSP
\ble \label{Lem49}
Suppose that $S = \mr{Spec} (k)$, and that we are given two dormant $(\mr{GL}_2^{(\N)}, \vartheta)$-opers $\nabla_1^\diamondsuit$, $\nabla_2^\diamondsuit$ on $\msX$.
Then,  $\nabla_1^\diamondsuit = \nabla_2^\diamondsuit$ if and only if $\delta (\nabla_1^\diamondsuit) = \delta (\nabla_2^\diamondsuit)$.
\ele
\begin{proof}
Since the morphism  \eqref{EQ200} for $j =0$ associated to each $(\mr{GL}_2^{(\N)}, \vartheta)$-oper is an isomorphism,
the assertion follows immediately  from the fact that $\mcD^{(\N -1)}$ is generated by sections of $\mcD_{\leq p^{\N -1}}^{(\N -1)}$ as a $\mcD^{(\N -1)}$-module.
\end{proof}
 \SSP

For each $i=1, \cdots, r$ and $m \in \mbZ_{\geq 0}$, 
the isomorphism \eqref{Eq48} restricts to an isomorphism of $\mcO_S$-modules 
 $\sigma^*_{i}({^L}\mcD^{(\N -1)}_{\leq m})\xrightarrow{\sim}\mcB_{S, \leq m}$,
 where $\mcB_{S, \leq m} := \bigoplus_{j=0}^m \mcO_S \cdot \partial_\mcB^{\langle j \rangle} \left(\subseteq \mcB_S \right)$.
 This isomorphism allows us to identify $\sigma_i^* (\mcD_{\leq m}^{(\N -1)} \otimes \varTheta)$ with $\mcB_{S, \leq m} \otimes \sigma_i^* (\varTheta)$.
Hence,
 the pull-back of $\delta (\nabla^\diamondsuit)$ via $\sigma_i$ determines a morphism
 \begin{align} \label{Eq112}
 \kappa_{i}(\nabla^\diamondsuit) := \sigma_i^* (\delta (\nabla^\diamondsuit))  \in \mr{Hom}_{\mcO_S}(\mcB_{S, \leq p^{\N-1}}\otimes\sigma_i^*(\varTheta), \mcB_{S, \leq 1} \otimes \sigma_i^* (\varTheta)).
 \end{align}
 
 \SSP
 \bpr \label{Prop3}
Let us fix $\rho := (\rho_i)_{i=1}^r \in \Xi_{2, \N}^{r\ma}$ ($=((\mbZ/p^\N \mbZ)^{\times}/\{ \pm 1\})^{r\ma}$ via \eqref{Eq311}). 
 Suppose that, for each $j=1, 2$, we are given  
 a dormant $(\mr{GL}_2^{(\N)}, \vartheta)$-oper $\nabla^\diamondsuit_j$ on $\msX$  of radii $\rho$ such that  (the corresponding $\mr{PGL}_2$-oper of) the $(\mr{GL}_2^{(1)}, \vartheta)$-oper obtained from $\nabla_j^\diamondsuit$ by reducing its level is normal, in the sense of ~\cite[Definition 4.53]{Wak5}.
 Then, 
 the following assertions hold:
 \begin{itemize}
 \item[(i)]
 For each $i=1, \cdots, r$,
 the equality $\mu (\nabla_1^\diamondsuit)^{\langle \bullet \rangle} = \mu (\nabla_2^\diamondsuit)^{\langle \bullet \rangle}$ holds.
 \item[(ii)]
 For each $i=1, \cdots, r$, the equality 
  $\kappa_{i}(\nabla_1^\diamondsuit) = \kappa_i (\nabla^\diamondsuit_2)$ holds.
 In particular, the difference $\delta (\nabla_1^\diamondsuit) - \delta (\nabla_2^\diamondsuit)$ specifies an $\mcO_S$-linear morphism
 $(\mcD_{\leq p^{\N -1}}^{(\N -1)}/ \mcD^{(\N -1)}_{\leq 1}) \otimes \varTheta \rightarrow \mcD_{\leq 1}^{(\N -1)}\otimes \varTheta (-D)$.
 \end{itemize}
 \epr
 \begin{proof}
 First, we shall consider assertion (i).
 After possibly replacing $S$ with its covering,
 we can assume that $\mr{Im}(\sigma_i)$ is defined by some function $t$ defined on its open neighborhood; it determines  an identification of $U_\oslash$ with the formal neighborhood of $\mr{Im}(\sigma_i)$.
 
If  $c$ denotes the exponent  of $\nabla_\vartheta$ at $\sigma_i$, then there exists a subset $\{ a, b \}$ of $\mbZ/p^\N \mbZ$ satisfying  $a \neq b$,  $a + b = c$, and $\overline{\left(\frac{a-b}{2}\right)} = \rho_i$ in $(\mbZ/p^\N \mbZ)/\{ \pm 1 \}$.
 For $j=1, 2$,
 there exists an isomorphism of $\mcD_\oslash^{(\N -1)}$-modules 
 \begin{align}
 \alpha_j : (\mcF_\varTheta, \nabla_j^\diamondsuit) |_{U_\oslash} \xrightarrow{\sim} (\mcO_\oslash^{\oplus 2}, \nabla_{\oslash, a}^{(\N -1)} \oplus \nabla_{\oslash, b}^{(\N -1)}).
 \end{align}
 It follows from ~\cite[Proposition 6.3, (i)]{Wak20} that the image of $\varTheta \left(\subseteq \mcF_\varTheta \right)$ under  $\alpha_j$ is generated by an element  $(v_j (t), u_j(t))$ of $H^0 (U_\oslash, (\mcO_\oslash^\times)^{\oplus 2})$, satisfying $v_j (0), u_j (0)\in H^0 (S, \mcO_S^\times)$.
 After possibly composing $\alpha_j$ with an automorphism of $\mcO_\oslash^{\oplus 2}$,
 we can assume that  $v_j (0) = u_j (0) = 1$.
  The restriction of 
  $\alpha_2^{-1}\circ \alpha_1$  to $\sigma_i$  defines  
  an automorphism $\overline{\alpha}_{21}$ of  $\mcB_{S, \leq 1} \otimes \sigma_i^*(\varTheta)$.
  By the above assumption, the
  subbundle $\sigma^*_i (\varTheta) \left(= \mcB_{S, \leq 0} \otimes \sigma^*_i (\varTheta) \right)$ of $\mcB_{S, \leq 1} \otimes \sigma^*_i (\varTheta)$
   is preserved by this automorphism.
  If $P$ denotes the upper-triangular invertible  $2 \times 2$ matrix corresponding to $\overline{\alpha}_{21}$,
 then 
 the equality  $P^{-1}\mu_i (\nabla^\diamondsuit_1)^{\langle p^a\rangle} P= \mu_i (\nabla^\diamondsuit_2)^{\langle p^a\rangle}$ holds for every  $a = 0, \cdots, \N -1$.
    Since we have assumed that both $\nabla_1^\diamondsuit$ and $\nabla_2^\diamondsuit$ are normal,
    $P$ must be a scalar matrix.
    This  implies 
  $\mu_i (\nabla_1^\diamondsuit)^{\langle \bullet \rangle} = \mu_i (\nabla_2^\diamondsuit)^{\langle \bullet \rangle}$, completing 
 the proof of assertion (i).
 
 Finally, assertion (ii) follows from assertion (i) and the definition of $\kappa_i (-)$.
  \end{proof}

\LSP
\subsection{The irreducibility of $\mcO p_{\N, \rho, g, r}^{^\mr{Zzz...}}$ for $(g, r) = (0, 4)$, $(1,1)$} \label{SS55}

As the first step in the proof of Theorem \ref{ThA},
we consider  the case where the moduli space is of dimension $1$.
Let us keep the notation 
from  the previous subsection, 
To begin with, we prove the following assertion.

\SSP
\bpr \label{Lem10}
We suppose that $(g, r) =(0, 4)$ or $(1, 1)$ (or equivalently, $\mr{dim} (\overline{\mcM}_{g, r}) = 1$), and that  the classifying morphism $S \rightarrow \overline{\mcM}_{g, r}$ of $\msX$ is finite and faithfully flat.
Then, for each positive integer $\ell$,
the vector bundle $f_* (\Omega^{\otimes \ell} (-D))$ does not admit any subbundle  isomorphic to $\mcO_S$.
\epr
\begin{proof}
Denote by  $\omega$  the canonical bundle on $X/S$ (hence $\Omega (-D) \cong \omega$).
 Let us  consider a decreasing  filtration $\{ \mcF^j \}_{j =0}^{\ell}$ on $\Omega^{\otimes \ell} (-D)$ defined by 
 \begin{align}
 \mcF^j := \Omega^{\otimes \ell} ((-j-1)D) \left(\cong  \omega^{\otimes \ell} ((\ell - j -1)D)\right)
 \end{align}
  for $j =0, \cdots, \ell -1$ and  $\mcF^{\ell} := 0$.
This induces a decreasing filtration $\{ f_* (\mcF^j)\}_{j=0}^{\ell}$ on  the direct image $f_* (\Omega^{\otimes \ell}(-D))$.
For each $j=0, \cdots, \ell -2$,  we have 
$\mr{deg}(\mcF^{j+1}) = (2g-2) \ell + (\ell - j - 1) r > 2g-2$, which implies  $\mbR^1 f_* (\mcF^{j+1}) = 0$.
It follows that 
 \begin{align}
 f_* (\mcF^j)/f_* (\mcF^{j+1}) 
 &\cong f_* (\mcF^j/\mcF^{j+1}) \notag \\
 & \cong f_* \left(\bigoplus_{i=1}^r \sigma_{i*}( \sigma_i^* (\omega^{\otimes \ell} ((\ell - j -1)D)))\right) \notag \\
& \cong \bigoplus_{i=1}^r \sigma_i^* (\omega^{\otimes \ell} (\ell - j -1) D) \notag \\
& \cong \bigoplus_{i=1}^r \sigma_i^{*}(\omega)^{\otimes (j+1)}, \notag 
 \end{align}
where  the last ``$\cong$" follows from the residue isomorphism $\sigma^*_i (\omega (D)) \left(\cong \sigma_i^* (\Omega) \right)\cong \mcO_S$.

Recall the following sequence of equalities:
  \begin{align}
  c_1 (\sigma_i^* (\omega)^{\otimes (j+1)}) = (j+1) \psi_i = \begin{cases} (j+1) \cdot \mr{deg}(S/\overline{\mcM}_{g, r}) & \text{if $(g, r) = (0, 4)$}; \\ \frac{1}{24} \cdot (j+1)\cdot \mr{deg}(S/\overline{\mcM}_{g, r}) & \text{if $(g, r) = (1, 1)$,}\end{cases}
  \end{align}
  where  $\psi_i$ denotes the  pull-back over $S$ of the $i$-th psi class on $\overline{\mcM}_{g, r}$ (cf. ~\cite[Chapter 17, Lemma 7.4, and  Chapter 19, Eq.\.(7.8)]{ACG}).
 It follows that, for any $j=0, \cdots, \ell -2$,  there is no line subbundle of $f_* (\mcF^{j})/f_* (\mcF^{j+1})$ isomorphic to $\mcO_S$.
 If  $(g, r) = (0, 4)$, then $f_* (\mcF^{\ell -1}) = f_* (\omega^{\otimes \ell}) = 0$.
 On the 
 other hand, if $(g, r) = (1, 1)$, then $f_* (\mcF^{\ell -1})$ is the pull-back of  the $(\ell -1)$-st  tensor power  of  the Hodge  bundle 
  on $\overline{\mcM}_{1, 1}$, which is not isomorphic to $\mcO_S$.
 From these considerations, $f_* (\Omega^{\otimes \ell} (-D))$ does not admit any line subbundle isomorphic to $\mcO_S$.
 This completes the proof of the assertion. 
\end{proof}
\SSP

By applying the above proposition,
we obtain the following assertion.

\SSP
\ble \label{T01}
Suppose that $(g, r) = (1, 1)$ or $(0, 4)$.
Then, for $\rho \in \Xi_{2, \N}^{r\ma}$, the moduli stack $\mcO p^{^\mr{Zzz...}}_{\N, \rho, g, r}$  (in the case $n=2$)
 is  connected (and hence, irreducible), provided  it is nonempty.
\ele
\begin{proof}
With the  notation as above, we  assume further  that $S = \overline{\mcM}_{g, r}$ and $\msX$ is taken to be the universal family of curves.
Also, write $S' := \mcO p^{^\mr{Zzz...}}_{\N, \rho, g, r}$.
For any  stack $Q$ over $S$, we use the  subscript ``$Q$" to denote  the result of base-changing  each object over $S$ to $Q$.
Let us take a dormant $2^{(\N)}$-theta characteristic $\vartheta := (\varTheta, \nabla_{\vartheta})$ of $\msX$.
(It follows from ~\cite[Proposition 5.14]{Wak20} that  such a $2^{(\N)}$-theta characteristic always exists.)
By ~\cite[Proposition 4.55]{Wak5}, 
there exists  a dormant $(\mr{GL}_2^{(\N)}, \vartheta_{S'})$-oper $\nabla^\diamondsuit$ on $\msX_{S'}$ corresponding to  
the universal family of dormant $\mr{PGL}_2^{(\N)}$-opers on  $\msX_{S'}$ such that the induced $\mr{PGL}_2$-oper is normal, in the sense of ~\cite[Definition 4.53]{Wak5}.

Now, suppose that there exist connected components $S_1$, $S_2$ of $S'$ with $S_1 \cap S_2 = \emptyset$.
We can choose a 
smooth proper
 Deligne-Mumford stack $T$ of dimension $1$  that dominates both $S_1$ and $S_2$ as a stack over $\overline{\mcM}_{g, r}$, i.e.,  fits into a commutative diagram  of 
$k$-stacks
\begin{align} \label{Com1}
\vcenter{\xymatrix@C=46pt@R=36pt{
T \ar[r] \ar[d] & S_1\ar[d] \\
S_2 \ar[r] & \overline{\mcM}_{g, r}.
}}
\end{align}
In particular, $T$ is finite and faithfully flat over $\overline{\mcM}_{g, r}$.

For each $i=1,2$,
let $h_i$ be the composite 
  $h_i : T \rightarrow S_i \hookrightarrow S'$.
The pull-back $h^*_i (\delta (\nabla^\diamondsuit))$ of $\delta (\nabla^\diamondsuit)$ (cf. \eqref{Eq110})
determines  a morphism 
$(\mcD_{\leq p^{\N -1}}^{(\N -1)})_T \otimes \varTheta_T \twoheadrightarrow (\mcD_{\leq 1}^{(\N -1)})_T \otimes \varTheta_T$.
The difference $h^*_1 (\delta (\nabla^\diamondsuit)) - h^*_2 (\delta (\nabla^\diamondsuit))$ vanishes when restricted to $(\mcD_{\leq 1}^{(\N -1)})_T \otimes \varTheta_T$, and thus determines a morphism
\begin{align}
\delta' : ( (\mcD_{\leq p^{\N -1}}^{(\N -1)})_T /(\mcD_{\leq 1}^{(\N -1)})_T ) \otimes \varTheta_T \twoheadrightarrow (\mcD_{\leq 1}^{(\N -1)})_T \otimes \varTheta_T (-D_T)
\end{align}
(cf. Proposition \ref{Prop3}, (ii)).
 We abuse notation by writing $\delta'$ for the corresponding 
global section of the vector bundle
\begin{align}
\mcV & := f_{T*} ((\mcD_{\leq 1}^{(\N -1)})_T \otimes \varTheta_T (-D_T) \otimes  (( (\mcD_{\leq p^{\N -1}}^{(\N -1)})_T /(\mcD_{\leq 1}^{(\N -1)})_T ) \otimes \varTheta_T)^\vee) \\
&  \left(\cong f_{T*}(\mcH om (( (\mcD_{\leq p^{\N -1}}^{(\N -1)})_T /(\mcD_{\leq 1}^{(\N -1)})_T ) \otimes \varTheta_T, (\mcD_{\leq 1}^{(\N -1)})_T \otimes \varTheta_T (-D_T))) \right) \notag
\end{align}
on $T$.

Since $S_1 \cap S_2 = \emptyset$,
it follows from Lemma \ref{Lem49} that
the image of the morphism $\mcO_T \rightarrow \mcV$ determined by $\delta'$ specifies a line subbundle of $\mcV$.
On the other hand, one can verify from 
  $\mcD^{(\N -1)}_{\leq j}/\mcD^{(\N -1)}_{\leq j-1} \cong \mcT^{\otimes j}$ that  the vector bundle $\mcV$ admits a filtration  whose graded pieces  are each  isomorphic to $f_{T*}(\Omega^{\otimes \ell} (-D)_T)$ for $\ell > 0$.
This contradicts Proposition \ref{Lem10} applied to  the case where ``$S$" is taken to be $T$.
That is to say,
$\mcO p^{^\mr{Zzz...}}_{\N, \rho, g, r} \left( =S'\right)$ must be connected, which completes the proof of this assertion.
 \end{proof}

\LSP
\subsection{The irreducibility of $\mcO p_{\N, \rho, g, r}^{^\mr{Zzz...}}$ for general $(g, r)$} \label{SS57}

To handle  the general case of pairs  $(g, r)$, we consider  gluing  morphisms on moduli stacks  of  various  types of higher-level $\mr{PGL}_2$-opers.
For the case of dormant opers, we refer to  ~\cite[Example 6.31]{Wak20}.

Let $g_1$, $g_2$, $r_1$, and $r_2$ be nonnegative integers  with $2g_i -1 + r_i > 0$ ($i=1,2$) and $g = g_1 + g_2$, $r = r_1 + r_2$.
These integers determine  the gluing morphism 
\begin{align} \label{EQ300}
\Phi_\mr{tree} : \overline{\mcM}_{g_1, r_1 +1} \times \overline{\mcM}_{g_2, r_2 +1} \rightarrow \overline{\mcM}_{g_1 + g_2, r_1 + r_2}
\end{align}
 obtained by attaching the respective last marked points of curves classified by $\overline{\mcM}_{g_1, r_1 +1}$ and $\overline{\mcM}_{g_2, r_2 +1}$ to form a node.

For $\rho_1\in \Xi_{2, \N}^{r_1\ma}$, $\rho_2\in\Xi_{2, \N}^{r_2\ma}$,   and $\rho_0 \in \Xi_{2, \N}$,
there exists a morphism
\begin{align} \label{EQ301}
\Phi_{\mr{tree}}^{\rho_1, \rho_2, \rho_0} : \mcO p^{^\mr{Zzz...}}_{\N, (\rho_1, \rho_0), g_1, r_1+1} \times \mcO p^{^\mr{Zzz...}}_{\N, (\rho_2, \rho_0), g_2, r_2 +1} \rightarrow \mcO p^{^\mr{Zzz...}}_{\N, (\rho_1, \rho_2), g,r}
\end{align}
obtained by gluing together a pair of  dormant $\mr{PGL}_2^{(\N)}$-opers classified by its domain  along the fibers over the respective last marked points of the underlying curves (cf.  ~\cite[Example 6.31, (i)]{Wak20}).
This morphism  makes   the following square   diagram commute:
\begin{align} \label{EQ305}
\vcenter{\xymatrix@C=46pt@R=36pt{
\mcO p^{^\mr{Zzz...}}_{\N, (\rho_1, \rho_0), g_1, r_1 +1} \times \mcO p^{^\mr{Zzz...}}_{\N, (\rho_2, \rho_0), g_2, r_2 +1} 
\ar[r]^-{\Phi^{\rho_1, \rho_2, \rho_0}_{\mr{tree}}} \ar[d]_-{\Pi'_{(\rho_1, \rho_0), g_1, r_1+1} \times \Pi'_{(\rho_2, \rho_0), g_2, r_2 +1}} 
& \mcO p^{^\mr{Zzz...}}_{\N, (\rho_1, \rho_2), g, r} 
 \ar[d]^-{\Pi'_{(\rho_1, \rho_2), g, r}}  \\
\overline{\mcM}_{g_1, r_1 +1} \times \overline{\mcM}_{g_2, r_2 +1} \ar[r]_-{\Phi_{\mr{tree}}} & \overline{\mcM}_{g, r}.
}}
\end{align}

Moreover,  if we  take $\widetilde{\rho}_1 \in \Xi_{2, \N+1}^{r_1\ma}$, $\widetilde{\rho}_2 \in \Xi_{2, \N +1}^{r_2\ma}$, $\widetilde{\rho}_0 \in \Xi_{2, \N+1}$ with $\widetilde{\rho}_1^{(\N)} = \rho_1$,  $\widetilde{\rho}_2^{(\N)} = \rho_2$, and $\widetilde{\rho}_0^{(\N)} = \rho_0$,
then the same type  of gluing procedure (cf. ~\cite[Section 6.4]{Wak20}) yields a morphism 
\begin{align} \label{Et13}
\Phi_{\mr{tree}}^{\widetilde{\rho}_1, \widetilde{\rho}_2, \widetilde{\rho}_0} : \mcO p^{^\mr{nilp}}_{\N+1, (\widetilde{\rho}_1, \widetilde{\rho}_0), g_1, r_1+1} \times \mcO p^{^\mr{nilp}}_{\N +1, (\widetilde{\rho}_2, \widetilde{\rho}_0), g_2, r_2 +1} \rightarrow \mcO p^{^\mr{nilp}}_{\N +1, (\widetilde{\rho}_1, \widetilde{\rho}_2), g,r}
\end{align}
such that the following diagram is commutative:
\begin{align} \label{Et44}
\vcenter{\xymatrix@C=46pt@R=36pt{
\mcO p^{^\mr{nilp}}_{\N +1, (\widetilde{\rho}_1, \widetilde{\rho}_0), g_1, r_1 +1} \times \mcO p^{^\mr{nilp}}_{\N+1, (\widetilde{\rho}_2, \widetilde{\rho}_0), g_2, r_2 +1} 
\ar[r]^-{\Phi^{\widetilde{\rho}_1, \widetilde{\rho}_2, \widetilde{\rho}_0}_{\mr{tree}}, \rho_0} \ar[d]_-{\Lambda_{(\widetilde{\rho}_1, \widetilde{\rho}_0), g_1, r_1 +1} \times \Lambda_{(\widetilde{\rho}_2, \widetilde{\rho}_0), g_2, r_2 +1}} 
& \mcO p^{^\mr{nilp}}_{\N+1, (\widetilde{\rho}_1, \widetilde{\rho}_2), g, r} 
 \ar[d]^-{\Lambda_{(\widetilde{\rho}_1, \widetilde{\rho}_2), g, r}}  \\
\mcO p^{^\mr{Zzz...}}_{\N, (\rho_1, \rho_0), g_1, r_1 +1} \times \mcO p^{^\mr{Zzz...}}_{\N, (\rho_2, \rho_0), g_2, r_2 +1} 
 \ar[r]_-{\Phi^{\rho_1, \rho_2, \rho_0}_{\mr{tree}}} & \mcO p^{^\mr{Zzz...}}_{\N, (\rho_1, \rho_2), g, r}.
}}
\end{align}

Next, given nonnegative integers $g$, $r$ with $2g + r > 0$, we shall write
\begin{align} \label{EQ302}
\Phi_\mr{loop} : \overline{\mcM}_{g, r+2} \rightarrow \overline{\mcM}_{g+1, r}
\end{align}
 for the gluing morphism obtained by attaching the last two marked points of each curve classified by $\overline{\mcM}_{g, r+2}$ to form a node.

For $\rho \in \Xi_{2, \N}^{r\ma}$ and $\rho_0 \in \Xi_{2, \N}$,
 there exists a morphism
 \begin{align} \label{EQ303}
 \Phi_{\mr{loop}}^{\rho, \rho_0} : \mcO p^{^\mr{Zzz...}}_{\N, (\rho, \rho_0, \rho_0), g, r+2} \rightarrow \mcO p^{^\mr{Zzz...}}_{\N, \rho, g+1, r}
 \end{align}
 obtained by gluing each dormant $\mr{PGL}_n^{(\N)}$-oper classified by its domain   along the fibers over the last two marked points of the underlying curve (cf.  ~\cite[Example 6.31, (ii)]{Wak20}).
 This morphism fits into  the following commutative  diagram:
 \begin{align} \label{EQ304}
\vcenter{\xymatrix@C=46pt@R=36pt{
\mcO p^{^\mr{Zzz...}}_{\N, (\rho, \rho_0, \rho_0), g, r+2}\ar[r]^-{\Phi^{\rho, \rho_0}_{\mr{loop}}} \ar[d]_-{\Pi'_{(\rho, \rho_0, \rho_0), g, r+2}} & \mcO p^{^\mr{Zzz...}}_{\N, \rho, g+1, r} \ar[d]^-{\Pi'_{\rho, g+1, r}} \\
\overline{\mcM}_{g, r +2} \ar[r]_-{\Phi_{\mr{loop}}} & \overline{\mcM}_{g+1, r}.
}}
\end{align}

If  we take $\widetilde{\rho} \in \Xi_{2, \N+1}^{r\ma}$ and $\widetilde{\rho}_0 \in \Xi_{2, \N +1}$ with $\widetilde{\rho}^{(\N)} = \rho$, $\widetilde{\rho}_0^{(\N)} = \rho_0$, then gluing $p^{\N +1}$-nilpotent $\mr{PGL}_2^{(\N +1)}$-opers gives rise to a morphism
 \begin{align} \label{Et55}
\vcenter{\xymatrix@C=46pt@R=36pt{
\mcO p^{^\mr{nilp}}_{\N +1, (\widetilde{\rho}, \widetilde{\rho}_0, \widetilde{\rho}_0), g, r+2}\ar[r]^-{\Phi^{\widetilde{\rho}, \widetilde{\rho}_0}_{\mr{loop}}} \ar[d]_-{\Lambda_{(\widetilde{\rho}, \widetilde{\rho}_0, \widetilde{\rho}_0), g, r+2}} & \mcO p^{^\mr{nilp}}_{\N +1, \widetilde{\rho}, g+1, r} \ar[d]^-{\Lambda_{\widetilde{\rho}, g +1, r}} \\
\mcO p^{^\mr{Zzz...}}_{\N, (\rho, \rho_0, \rho_0), g, r+2}
 \ar[r]_-{\Phi_{\mr{loop}}} & \mcO p^{^\mr{Zzz...}}_{\N, \rho, g+1, r}.
}}
\end{align}

Using the diagrams \eqref{Et44} and \eqref{Et55}, we prove the following assertion.

\SSP
\ble \label{Cor11}
Let  us set $d := 3g-3 +r$.
 Suppose  that $d \geq 1$ and  that
the following conditions are fulfilled: 
\begin{itemize}
\item[(a)]
The stack $\mcO p_{\N, \rho', g', r'}^{^\mr{Zzz...}}$  for any $\rho' \in \Xi_{2, \N}^{r'\ma}$ is either empty or irreducible whenever $3g' -3+r' \leq  d$;
\item[(b)]
The stack $\mcO p_{\N +1, \rho', g', r'}^{^\mr{nilp}}$  for any $\rho' \in \Xi_{2, \N+1}^{r'\ma}$ is either empty or connected whenever $3g' -3+r' <  d$.
\end{itemize}
Also, let $\widetilde{\rho}$ be an element of $\Xi_{2, \N}^{r\ma}$ satisfying $\mcO p_{\N+1, \widetilde{\rho}, g, r}^{^\mr{nilp}} \neq \emptyset$.
Then, the following assertions hold.
\begin{itemize}
\item[(i)]
Any two irreducible components of $\mcO p_{\N+1, \widetilde{\rho}, g, r}^{^\mr{nilp}}$ intersect.
In particular, $\mcO p_{\N+1, \widetilde{\rho}, g, r}^{^\mr{nilp}}$ is connected.
\item[(ii)]
The stack $\mcO p_{\N+1, \widetilde{\rho}, g, r}^{^\mr{Zzz...}}$ is irreducible if it is nonempty.
\end{itemize}
\ele
\begin{proof}
Since $\mcO p_{\N+1, \widetilde{\rho}, g, r}^{^\mr{Zzz...}}$ is smooth over $k$ (cf. Theorem \ref{Th8}, (iii)),
the assertion (ii) follows from assertion (i).
Hence, it suffices to prove assertion (i). 
For simplicity, we set  $\rho := \widetilde{\rho}^{(\N)} \in \Xi_{2, \N}^{r\ma}$.

 First, let us assume that $g \geq 1$.
Since  the commutative square diagram 
  \begin{align} \label{Eq419}
\vcenter{\xymatrix@C=66pt@R=36pt{
\coprod_{\rho_0 \in \Xi_{2, \N}}\mcO p^{^\mr{Zzz...}}_{\N, \rho_0, 1, 1} \times \mcO p^{^\mr{Zzz...}}_{\N, (\rho, \rho_0), g-1, r +1} \ar[r]^-{\coprod_{\rho_0}\Phi^{\emptyset, \rho, \rho_0}_{\mr{tree}}} \ar[d]_-{\coprod_{\rho_0} \Pi'_{\rho_0, 1, 1} \times \Pi'_{(\rho, \rho_0), g-1, r+1}} & \mcO p^{^\mr{Zzz...}}_{\N, \rho, g, r} \ar[d]^-{\Pi'_{\rho, g, r}} \\
\overline{\mcM}_{1, 1} \times \overline{\mcM}_{g -1, r+1}
 \ar[r]_-{\Phi_{\mr{tree}}} &
  \overline{\mcM}_{g, r}
}}
\end{align}
 is Cartesian and the right-hand vertical arrow is   surjective  (cf. Theorems \ref{Th8}, (iii)),
there exists an  element $\rho_0 \in \Xi_{2, \N}$ satisfying  $\mcO p^{^\mr{Zzz...}}_{\N, \rho_0, 1, 1} \times \mcO p^{^\mr{Zzz...}}_{\N, (\rho, \rho_0), g-1, r +1} \neq \emptyset$.
 
 The morphisms $\Phi_{\mr{tree}}^{\emptyset, \widetilde{\rho}, \widetilde{\rho}_0}$ (cf.  \eqref{Et13}) and 
 the diagrams  \eqref{Et44}
  defined  for $(g_1, g_2, r_1, r_2) = (1, g-1, 0, r)$ and $\widetilde{\rho}_0 \in {^\dagger}\pi^{-1}(\rho_0)\cap \Xi_{2, \N +1}$ fit into the following Cartesian diagram: 
 \begin{align} \label{Eq419eee}
\vcenter{\xymatrix@C=66pt@R=36pt{
\coprod_{\widetilde{\rho}_0 \in {^\dagger}\pi^{-1}(\rho_0) \cap \Xi_{2, \N+1}}\mcO p^{^\mr{nilp}}_{\N+1, \widetilde{\rho}_0, 1, 1} \times \mcO p^{^\mr{nilp}}_{\N+1, (\widetilde{\rho}, \widetilde{\rho}_0), g-1, r +1} \ar[r]^-{\coprod_{\widetilde{\rho}_0}\Phi^{\emptyset, \widetilde{\rho}, \widetilde{\rho}_0}_{\mr{tree}}} \ar[d]_-{\coprod_{\widetilde{\rho}_0}\Lambda_{\widetilde{\rho}_0, g, r} \times \Lambda_{(\widetilde{\rho}, \widetilde{\rho}_0), g-1, r +1}} & \mcO p^{^\mr{nilp}}_{\N+1, \widetilde{\rho}, g, r} \ar[d]^-{\Lambda_{\widetilde{\rho}, g, r}} \\
\mcO p_{\N, \rho_0, 1, 1}^{^\mr{Zzz...}} \times \mcO p^{^\mr{Zzz...}}_{\N, (\rho, \rho_0), g-1, r+1}
 \ar[r]_-{\Phi^{\emptyset, \rho, \rho_0}_{\mr{tree}}} &
 \mcO p_{\N, \rho, g, r}^{^\mr{Zzz...}}.
}}
\end{align}
By  the  condition   (b), 
   the image of  $\Phi_{\mr{tree}}^{\emptyset, \widetilde{\rho}, \widetilde{\rho}_0}$  is  connected for any such $\widetilde{\rho}_0$.
 On the other hand, 
since $\mcO p_{\N, \rho, g, r}^{^\mr{Zzz...}}$ is irreducible by the condition (a) and 
 $\Lambda_{\widetilde{\rho}, g, r}$ is finite and faithfully flat (cf. Theorem \ref{Th665}),
the restriction of $\Lambda_{\widetilde{\rho}, g, r}$ to  each irreducible  component of $\mcO p_{\N+1, \widetilde{\rho}, g, r}^{^\mr{nilp}}$ is surjective.
Thus, to complete the proof for $g \geq 1$,
it suffices to show  that the images of the various morphisms  $\Phi_{\mr{tree}}^{\emptyset, \widetilde{\rho}, \widetilde{\rho}_0}$'s (for different  $\widetilde{\rho}_0$'s with $\widetilde{\rho}^{(\N)}_0 = \rho_0$) are all contained in  the same connected component of $\mcO p_{\N +1, \widetilde{\rho}, g, r}^{^\mr{nilp}}$.

Let us take two elements $\widetilde{\rho}_i \in {^\dagger}\pi^{-1}(\rho_0) \cap \Xi_{2, \N+1}$ ($i=1, 2$) with
 $\mcO p^{^\mr{nilp}}_{\N +1, \widetilde{\rho}_i, 1, 1} \times \mcO p^{^\mr{nilp}}_{\N +1, (\widetilde{\rho}, \widetilde{\rho}_i), g-1, r +1} \neq \emptyset$.
Since 
 the square diagram
\begin{align} \label{EQ348}
\vcenter{\xymatrix@C=46pt@R=36pt{
 \coprod_{\zeta \in \Xi_{2, \N}} \mcO p^{^\mr{Zzz...}}_{(\rho_0, \zeta, \zeta), 0, 3} \ar[r]^-{\Phi_\mr{loop}^{\rho_0, \zeta}} \ar[d]_-{\coprod_{\zeta} \Pi'_{(\rho_0, \zeta, \zeta), 0, 3}} & \mcO p_{\N, \rho_0, 1, 1}^{^\mr{Zzz...}} \ar[d]^-{\Pi'_{\rho_0, 1, 1}} \\
 \overline{\mcM}_{0, 3} \ar[r]_-{\Phi_\mr{loop}} & \overline{\mcM}_{1, 1}
 }}
\end{align}
is Cartesian and its right-hand vertical arrow is surjective (cf. Theorem \ref{Th8}, (iii)),
one can find an element $\zeta \in \Xi_{2, \Xi}$ satisfying $\mcO p^{^\mr{Zzz...}}_{(\rho_0, \zeta, \zeta), 0, 3} \neq \emptyset$, or equivalently, $\mcO p^{^\mr{Zzz...}}_{(\zeta, \zeta, \rho_0), 0, 3} \neq \emptyset$.
Let us take $\widetilde{\zeta} \in {^\dagger}\pi^{-1}(\zeta) \cap \Xi_{2, \N +1}$.
Then, it follows from Corollary \ref{Cor15} that both $\mcO p_{\N+ 1, (\widetilde{\zeta}, \widetilde{\zeta}, \widetilde{\rho}_1), 0, 3}^{^\mr{nilp}}$ and 
$\mcO p_{\N+ 1, (\widetilde{\zeta}, \widetilde{\zeta}, \widetilde{\rho}_2), 0, 3}^{^\mr{nilp}}$ are nonempty.

For $i=1, 2$,  the square  diagram
\begin{align} \label{EQ308}
\vcenter{\xymatrix@C=56pt@R=36pt{
\mcO p^{^\mr{nilp}}_{N+1, (\widetilde{\zeta}, \widetilde{\zeta}, \widetilde{\rho}_i), 0, 3} \times 
\mcO p^{^\mr{nilp}}_{N+1, (\widetilde{\rho}, \widetilde{\rho}_i), g-1, r+1} 
\ar[r]^-{\Phi_\mr{loop}^{\widetilde{\rho}_i, \widetilde{\zeta}} \times \mr{id}} 
\ar[d]_-{\Phi_\mr{tree}^{(\widetilde{\zeta}, \widetilde{\zeta}), \widetilde{\rho}, \widetilde{\rho}_i}} &
\mcO p^{^\mr{nilp}}_{N+1, \widetilde{\rho}_i, 1, 1} \times 
\mcO p^{^\mr{nilp}}_{N+1, (\widetilde{\rho}, \widetilde{\rho}_i), g-1, r+1} 
\ar[dd]^-{\Phi_\mr{tree}^{\emptyset, \widetilde{\rho}, \widetilde{\rho}_i}} \\
\mcO p^{^\mr{nilp}}_{N+1, (\widetilde{\zeta}, \widetilde{\zeta}, \widetilde{\rho}), g-1, r+2} 
\ar[d]_-{\varsigma_{g-1, r+2}}
& \\
\mcO p^{^\mr{nilp}}_{N+1, (\widetilde{\rho}, \widetilde{\zeta}, \widetilde{\zeta}), g-1, r+2}  
\ar[r]_-{\Phi_\mr{loop}^{\widetilde{\rho}, \widetilde{\zeta}}} & 
\mcO p^{^\mr{nilp}}_{N+1, \widetilde{\rho}, g, r}
}}
\end{align}
 is commutative, 
 where  $\varsigma_{g-1, r+2}$ denotes the isomorphism
 $\mcO p^{^\mr{nilp}}_{N+1, (\widetilde{\zeta}, \widetilde{\zeta}, \widetilde{\rho}), g-1, r+2} 
 \xrightarrow{\sim} \mcO p^{^\mr{nilp}}_{N+1, (\widetilde{\rho}, \widetilde{\zeta}, \widetilde{\zeta}), g-1, r+2}$
  induced by reordering the marked points of underlying curves.
 Since
its  upper-left corner is nonempty, we have
 $\mr{Im} (\Phi_{\mr{tree}}^{\emptyset, \widetilde{\rho}, \widetilde{\rho}_i}) \cap \mr{Im} (\Phi_\mr{loop}^{\widetilde{\rho}, \widetilde{\zeta}})\neq \emptyset$.
Moreover, by the condition (b), 
 $\mr{Im} (\Phi_\mr{loop}^{\widetilde{\rho}, \widetilde{\zeta}})$ is connected.
 It follows that $\mr{Im} (\Phi_{\mr{tree}}^{\emptyset, \widetilde{\rho}, \widetilde{\rho}_1})$ and $\mr{Im} (\Phi_{\mr{tree}}^{\emptyset, \widetilde{\rho}, \widetilde{\rho}_2})$ are contained in the same connected component of  $\mcO p_{\N+1, \widetilde{\rho}, g, r}^{^\mr{nilp}}$.
 This completes the proof for $g \geq 1$.

Next, assume that $g =0$, which implies $r \geq 4$.
(The following discussion proceeds in essentially the same way as the previous case of $g \geq 1$.)
Let us write $\widetilde{\rho} := (\widetilde{\alpha}, \widetilde{\beta}, \widetilde{\gamma}, \widetilde{\lambda})$ for 
$\widetilde{\alpha}, \widetilde{\beta}, \widetilde{\gamma} \in \Xi_{2, \N+1}$, $\widetilde{\lambda} \in \Xi_{2, \N +1}^{r-3\ma}$, and write 
$\alpha := \widetilde{\alpha}^{(\N)}$, $\beta := \widetilde{\beta}^{(\N)}$,  $\gamma := \widetilde{\gamma}^{(\N)}$, and $\lambda := \widetilde{\lambda}^{(\N)}$ (hence $\rho = (\alpha, \beta, \gamma, \lambda)$).
Since  the commutative square diagram
 \begin{align} \label{Et330}
\vcenter{\xymatrix@C=66pt@R=36pt{
\coprod_{\rho_0 \in \Xi_{2, \N}}  \mcO p^{^\mr{Zzz...}}_{\N, (\alpha, \beta, \gamma, \rho_0), 0, 4} \times \mcO p^{^\mr{Zzz...}}_{\N, (\lambda, \rho_0), 0, r-2} \ar[r]^-{\Phi_\mr{loop}^{(\alpha, \beta, \gamma), \lambda, \rho_0}} \ar[d]_-{\coprod_{\rho_0}\Pi'_{(\alpha, \beta, \gamma, \rho_0), 0, 4} \times \Pi'_{(\lambda, \rho_0), 0, r-2}} & \mcO p_{\N, \rho, 0, r}^{^\mr{Zzz...}} \ar[d]^-{\Pi'_{\rho, 0, r}}
\\
\overline{\mcM}_{0, 4} \times \overline{\mcM}_{0, r-2} \ar[r]_-{\Phi_\mr{loop}} & \overline{\mcM}_{0, r}
}}
\end{align}
is Cartesian and the right-hand vertical arrow is surjective (cf. Theorem \ref{Th8}, (iii)),
there exists 
an element $\rho_0 \in \Xi_{2, \N}$ satisfying 
$ \mcO p^{^\mr{Zzz...}}_{\N, (\alpha, \beta, \gamma, \rho_0), 0, 4} \times \mcO p^{^\mr{Zzz...}}_{\N, (\lambda, \rho_0), 0, r-2} \neq \emptyset$.

The morphisms $\Phi_\mr{loop}^{(\widetilde{\alpha}, \widetilde{\beta}, \widetilde{\gamma}), \widetilde{\lambda}, \widetilde{\rho}_0}$
 defined  for $(g_1, g_2, r_1, r_2) = (0, 0, 3, r-3)$ and $\widetilde{\rho}_0 \in {^\dagger}\pi^{-1}(\rho_0) \cap \Xi_{2, \N +1}$ fit into the following Cartesian diagram:
  \begin{align} \label{Eq419v}
\vcenter{\xymatrix@C=66pt@R=36pt{
\coprod_{\widetilde{\rho}_0 \in {^\dagger}\pi^{-1}(\rho_0)\cap \Xi_{2, \N+1}}\mcO p^{^\mr{nilp}}_{\N+1, (\widetilde{\alpha}, \widetilde{\beta}, \widetilde{\gamma}, \widetilde{\rho}_0), 0, 4} \times \mcO p^{^\mr{nilp}}_{\N +1, (\widetilde{\lambda}, \widetilde{\rho}_0), 0, r -2} \ar[r]^-{\coprod_{\widetilde{\rho}_0}\Phi^{(\widetilde{\alpha}, \widetilde{\beta}, \widetilde{\gamma}), \widetilde{\lambda}, \widetilde{\rho}_0}_{\mr{tree}}} \ar[d]_-{\coprod_{\rho_0}\Lambda_{(\widetilde{\alpha}, \widetilde{\beta}, \widetilde{\gamma}, \widetilde{\rho}_0), 0, 4} \times \Lambda_{(\widetilde{\lambda}, \widetilde{\rho}_0), 0, r -2}} & \mcO p^{^\mr{nilp}}_{\N +1, \widetilde{\rho}, 0, r} \ar[d]^-{\Lambda_{\widetilde{\rho}, 0, r}} \\
\mcO p^{^\mr{Zzz...}}_{\N, (\alpha, \beta, \gamma, \rho_0), 0, 4} \times \mcO p^{^\mr{Zzz...}}_{\N, (\lambda, \rho_0), 0, r-2}
 \ar[r]_-{\Phi_{\mr{tree}}} & \mcO p^{^\mr{Zzz...}}_{\N, \rho, 0, r}.
}}
\end{align}
By the condition (b), the image of $\Phi_\mr{loop}^{(\widetilde{\alpha}, \widetilde{\beta}, \widetilde{\gamma}), \widetilde{\lambda}, \widetilde{\rho}_0}$ is connected for any such $\widetilde{\rho}_0$.
On the other hand, since $\mcO p^{^\mr{Zzz...}}_{\N, \rho, 0, r}$ is irreducible by the condition (a) and $\Lambda_{\widetilde{\rho}, 0, r}$ is finite and faithfully flat (cf. Theorem \ref{Th665}),
the restriction of $\Lambda_{\widetilde{\rho}, 0, r}$ to each irreducible  component of $\mcO p^{^\mr{nilp}}_{\N +1, \widetilde{\rho}, 0, r}$ is surjective.
Thus, to complete the proof for $g =0$,
it suffices to show that the image of the various morphisms  $\Phi_\mr{loop}^{(\widetilde{\alpha}, \widetilde{\beta}, \widetilde{\gamma}), \widetilde{\lambda}, \widetilde{\rho}_0}$ (for $\widetilde{\rho}_0$'s with $\widetilde{\rho}_0^{(\N)} = \rho_0$) lie  in the same connected component of $\mcO p_{\N +1, \widetilde{\rho}, 0, r}^{^\mr{nilp}}$.

Let us take two elements $\widetilde{\rho}_i \in {^\dagger}\pi^{-1}(\rho_0) \cap \Xi_{2, \N+1}$ ($i=1, 2$) with $\mcO p^{^\mr{nilp}}_{\N+1, (\widetilde{\alpha}, \widetilde{\beta}, \widetilde{\gamma}, \widetilde{\rho}_i), 0, 4} \times \mcO p^{^\mr{nilp}}_{\N +1, (\widetilde{\lambda}, \widetilde{\rho}_i), 0, r -2}$ $\neq \emptyset$.
Since  the commutative square diagram
 \begin{align} \label{Et3301}
\vcenter{\xymatrix@C=66pt@R=36pt{
\coprod_{\zeta \in \Xi_{2, \N}} \mcO p_{\N, (\alpha, \beta, \zeta), 0, 3}^{^\mr{Zzz...}} \times \mcO p_{\N, (\gamma, \rho_0, \zeta), 0, 3}^{^\mr{Zzz...}} \ar[r]^-{\Phi_{\mr{tree}^{(\alpha, \beta), (\gamma, \rho_0), \zeta}}} \ar[d]_-{\coprod_\zeta \Pi'_{(\alpha, \beta, \zeta), 0, 3} \times \Pi'_{(\gamma, \rho_0, \zeta), 0, 3}} & \mcO p^{^\mr{Zzz...}}_{\N, (\alpha, \beta, \gamma, \rho_i), 0, 4} \ar[d]^-{\Pi'_{(\alpha, \beta, \gamma, \rho_0), 0, 4}} \\
\overline{\mcM}_{0, 3} \times \overline{\mcM}_{0, 3} \ar[r] & \overline{\mcM}_{0, 4} 
}}
\end{align}
is Cartesian and its right-hand vertical arrow is surjective (cf. Theorem \ref{Th8}, (iii)),
one can find an element $\zeta \in \Xi_{2, \N}$ satisfying $\mcO p_{\N, (\alpha, \beta, \zeta), 0, 3}^{^\mr{Zzz...}} \times \mcO p_{\N, (\gamma, \rho_0, \zeta), 0, 3}^{^\mr{Zzz...}} \neq \emptyset$.
Let us take $\widetilde{\zeta} \in {^\dagger}\pi^{-1}(\zeta) \cap \Xi_{2, \N +1}$.
Then, it follows from Corollary \ref{Cor15} that  both 
$\mcO p_{\N+1, (\widetilde{\alpha}, \widetilde{\beta}, \widetilde{\zeta}), 0, 3}^{^\mr{nilp}} \times \mcO p_{\N+1, (\widetilde{\gamma}, \widetilde{\rho}_1, \widetilde{\zeta}), 0, 3}^{^\mr{nilp}}$ and $\mcO p_{\N+1, (\widetilde{\alpha}, \widetilde{\beta}, \widetilde{\zeta}), 0, 3}^{^\mr{nilp}} \times \mcO p_{\N+1, (\widetilde{\gamma}, \widetilde{\rho}_2, \widetilde{\zeta}), 0, 3}^{^\mr{nilp}}$ are nonempty.

For $i=1, 2$, the  diagram
\begin{align} \label{EQ308v}
\vcenter{\xymatrix@C=56pt@R=36pt{
\mcO p^{^\mr{nilp}}_{\N+1, (\widetilde{\alpha}, \widetilde{\beta}, \widetilde{\zeta}), 0, 3} \times \mcO p^{^\mr{nilp}}_{\N+1, (\widetilde{\gamma}, \widetilde{\rho}_i, \widetilde{\zeta}), 0, 3} \times \mcO p^{^\mr{nilp}}_{\N+1, (\widetilde{\lambda}, \widetilde{\rho}_i), 0, r-2} 
\ar[r]^-{\Phi_\mr{tree}^{(\widetilde{\alpha}, \widetilde{\beta}), (\widetilde{\gamma}, \widetilde{\rho}_i), \widetilde{\zeta}} \times \mr{id}} 
\ar[d]_-{\mr{id}\times \varsigma_{0, 3} \times \mr{id}}
& \mcO p^{^\mr{nilp}}_{\N+1, (\widetilde{\alpha}, \widetilde{\beta}, \widetilde{\gamma}, \widetilde{\rho}_i), 0, 4} \times \mcO p^{^\mr{nilp}}_{\N+1, (\widetilde{\lambda}, \widetilde{\rho}_i), 0, r-2}  \ar[dd]^-{\Phi_\mr{tree}^{(\widetilde{\alpha}, \widetilde{\beta}, \widetilde{\gamma}), \widetilde{\lambda}, \widetilde{\rho}_i}} \\
\mcO p^{^\mr{nilp}}_{\N+1, (\widetilde{\alpha}, \widetilde{\beta}, \widetilde{\zeta}), 0, 3} \times \mcO p^{^\mr{nilp}}_{\N +1, (\widetilde{\gamma}, \widetilde{\zeta}, \widetilde{\rho}_i), 0, 3} \times \mcO p^{^\mr{nilp}}_{\N+1, (\widetilde{\lambda}, \widetilde{\lambda}), 0, r-2} 
\ar[d]_-{\mr{id} \times \Phi_\mr{tree}^{(\widetilde{\gamma}, \widetilde{\zeta}), (\widetilde{\lambda}, \widetilde{\rho}_i), \widetilde{\zeta}}} & \\
\mcO p^{^\mr{nilp}}_{\N+1, (\widetilde{\alpha}, \widetilde{\beta}, \widetilde{\zeta}), 0, 3} \times \mcO p^{^\mr{nilp}}_{\N+1, (\widetilde{\gamma}, \widetilde{\lambda}, \widetilde{\zeta}), 0, r-1} \ar[r]_-{\Phi_\mr{tree}^{(\widetilde{\alpha}, \widetilde{\beta}), (\widetilde{\gamma}, \widetilde{\lambda}), \widetilde{\zeta}}} 
& \mcO p^{^\mr{nilp}}_{\N+1, \widetilde{\rho}, 0, r}
 }}
\end{align}
is commutative, where $\varsigma_{0, 3}$ denotes the isomorphism $\mcO p^{^\mr{nilp}}_{\N+1, (\widetilde{\gamma}, \widetilde{\rho}_i, \widetilde{\zeta}), 0, 3} \xrightarrow{\sim}\mcO p^{^\mr{nilp}}_{\N+1, (\widetilde{\gamma}, \widetilde{\zeta}, \widetilde{\rho}_i), 0, 3}$ induced by reordering the marked points of underlying  curves.
Since its upper-left corner is nonempty,
we have $\mr{Im} (\Phi_\mr{tree}^{(\widetilde{\alpha}, \widetilde{\beta}, \widetilde{\gamma}), \widetilde{\lambda}, \widetilde{\rho}_i}) \cap \mr{Im}(\Phi_\mr{tree}^{(\widetilde{\alpha}, \widetilde{\beta}), (\widetilde{\gamma}, \widetilde{\lambda}), \widetilde{\zeta}}) \neq \emptyset$.
Moreover, by the condition (b), 
$\mr{Im}(\Phi_\mr{tree}^{(\widetilde{\alpha}, \widetilde{\beta}), (\widetilde{\gamma}, \widetilde{\lambda}), \widetilde{\zeta}})$ is connected.
It follows that  $\mr{Im} (\Phi_\mr{tree}^{(\widetilde{\alpha}, \widetilde{\beta}, \widetilde{\gamma}), \widetilde{\lambda}, \widetilde{\rho}_1})$
and $\mr{Im} (\Phi_\mr{tree}^{(\widetilde{\alpha}, \widetilde{\beta}, \widetilde{\gamma}), \widetilde{\lambda}, \widetilde{\rho}_2})$
are contained in the same connected component of  $\mcO p^{^\mr{nilp}}_{\N +1, \widetilde{\rho}, 0, r}$.
This completes the proof for $g =0$.
\end{proof}
\SSP

By applying the above result, we conclude the following theorem.

\SSP
\bt[cf. Theorem \ref{ThA}]\label{Cor49}
(Recall that we have assumed  the equality  $n=2$.)
Let $\rho$ be an element of $\Xi_{2, \N}^{r\ma}$, where we set $\rho := \emptyset$ when $r = 0$.
Then, the following assertions hold.
\begin{itemize}
\item[(i)]
If $\mcO p^{^\mr{nilp}}_{\N, \rho, g, r} \neq \emptyset$, then 
 any two irreducible   components of  $\mcO p^{^\mr{nilp}}_{\N, \rho, g, r}$ intersect.
In particular, $\mcO p^{^\mr{nilp}}_{\N, \rho, g, r}$ is connected (when it is nonempty).
\item[(ii)]
If $\mcO p_{\N, \rho, g, r}^{^\mr{Zzz...}} \neq \emptyset$, then
$\mcO p_{\N, \rho, g, r}^{^\mr{Zzz...}}$ is irreducible.
\end{itemize}
\et
\begin{proof}
We shall  set $d := 3g-3+r$.
Assertions (i), (ii) for  $\N =1$ and  $d\geq 0$ have been  established  in ~\cite[Chapter II, Theorems 1.12 and 2.8]{Moc2}.
Assertion (i) for $\N \geq 1$ and $d = 0$  follows from Corollary \ref{Cor15}.
Moreover, according to Lemma \ref{T01},   assertion (ii) holds when  $\N \geq 1$ and $d =1$.
The remaining cases can be proved by induction on $\N$ and $d$, based on  these base   cases and  Lemma \ref{Cor11}.
\end{proof}

\LSP
\subsection{The projective system of level-reduction maps} \label{SS50}

We here present   an assertion deduced  from the previous theorem, which concerns  the process of  raising the level of each dormant $\mr{PGL}_2^{(\N)}$-oper.

Let $\rho := (\rho_i)_{i=1}^r$ be an $r$-tuple of elements in  $\mbZ_p^{\times}/\{ \pm 1 \}$, where $\mbZ_p := \varprojlim_{\N} \mbZ /p^\N \mbZ$.
For each $\N \in \mbZ_{> 0}$, we denote by $\rho_{i, \N}$ the image of $\rho_i$ via   the natural projection $\mbZ_p^{\times}/\{ \pm 1 \} \rightarrow  (\mbZ/p^\N \mbZ)^{\times}/\{ \pm 1 \}$, and write $\rho_\N := (\rho_{i, \N})_{i=1}^r$.
Under the identification $(\mbZ/p^\N)^\times /\{ \pm 1 \} = \Xi_{2, \N}$ (cf. Example \ref{Exam009}),
the level reductions    of dormant $\mr{PGL}_2^{(\N)}$-opers (for various $\N$'s) yield
 a projective system of $k$-stacks 
\begin{align} \label{Eq12}
\cdots \xrightarrow{\Pi_{\N +1 \Rightarrow \N}} \mcO p^{^\mr{Zzz...}}_{N, \rho_\N, g, r} \xrightarrow{\Pi_{\N \Rightarrow \N -1}} \cdots \xrightarrow{\Pi_{3 \Rightarrow 2}}  \mcO p^{^\mr{Zzz...}}_{2, \rho_2, g, r} \xrightarrow{\Pi_{2 \Rightarrow 1}}  \mcO p^{^\mr{Zzz...}}_{1, \rho_1, g, r}.
\end{align}

\SSP
\bpr \label{Prop59}
(Recall that we have assumed the equality $n=2$.)
Let us retain  the above notation, and suppose that  $\mcO p^{^\mr{Zzz...}}_{\N, \rho_\N, g, r} \neq \emptyset$ for a positive integer  $\N$.
Then, 
the morphism $\Pi_{\N' \Rightarrow \N'-1}$ for every $\N' \leq \N$
 in the projective system \eqref{Eq12} is  finite,  
faithfully flat, and generically \'{e}tale.
\epr
\begin{proof}
The assertion follows from a standard argument using
 Theorem \ref{Cor49}, (ii), and  the various geometric properties of    $\mcO p^{^\mr{Zzz...}}_{\N, \rho, g, r}$  asserted in Theorem \ref{Th8}, (ii) and (iii).
\end{proof}

\SSP
\bco \label{Prop31}
Let  $S$ be a scheme over $k$ and $\msX$ an unpointed  stable curve over $S$ of genus $g >1$.
Also, let $\N$, $\N'$ be positive integers with $\N' \leq \N$ and $\msE^\spadesuit$  a dormant $\mr{PGL}_2^{(\N')}$-oper on $\msX$.
Then, after possibly changing the base over a finite,  faithfully flat, and generically \'{e}tale  morphism $T \rightarrow S$,
there exists a dormant $\mr{PGL}_2^{(\N)}$-oper  whose level reduction to $\N'$  is isomorphic to  $\msE^\spadesuit$.
\eco
\begin{proof}
The assertion follows from Proposition \ref{Prop59} and the fact that $\mcO p^{^\mr{Zzz...}}_{\N, \emptyset, g, 0} \neq \emptyset$ (cf. ~\cite[Corollary 6.27, (ii)]{Wak20}).
\end{proof}

\vspace{10mm}
\section{Applications to certain lifting  properties} \label{S8}
\LSP

In this  section, we present 
 several  applications of 
the surjectivity of the morphisms in \eqref{Eq12}, as established in Proposition \ref{Prop59}.

Denote by $\mcM_g$ the moduli stack classifying smooth proper curves of genus $g>1$ over $k$.
Throughout this section, 
we fix  a  curve  $X$  over $k$ classified by $\mcM_g$.

\LSP
\subsection{Maximally Frobenius-destabilized bundles} \label{SS52}

We begin by addressing the problem of   whether a given vector bundle $\mcV$ admits a  Frobenius descent; that is, whether  there exists a vector bundle $\mcU$ such that  $F_X^*(\mcU) \cong \mcV$, or more generally, $F_X^{\M *}(\mcU)$ for some  $\M \in \mbZ_{> 0}$. 
When  restricted to the  class of stable bundles, this question amounts to knowing the (set-theoretic) image of the {\it generalized Vershiebung   rational map}
\begin{align} \label{Ey1}
\mr{Ver}^n_X : SU^n_{X}  \dashrightarrow  SU^n_{X}
\end{align}
on the moduli space $SU^n_{X}$  of rank $n$ stable bundles on $X$ with trivial determinant (cf. ~\cite[Lemma A]{MeSu}).
To be precise, this rational map is defined by  the assignment $\mcU \mapsto F^*_X (\mcU)$.
Since the Frobenius pull-back of a stable bundle may not remain  stable, $\mr{Ver}^n_X$ is not necessarily defined on all of  $SU^n_{X}$.

Unfortunately, 
our current understanding of this rational map is limited, as few general results are known for arbitrary triples $(p, n, g)$.
Relevant studies in specific cases can
 be found in   ~\cite{LanPa},  ~\cite{LasPa2}, ~\cite{LasPa3}, and ~\cite{Oss1}.

We now  focus on the case $n=2$.
For a rank $2$ vector bundle $\mcV$ of degree $0$,
the stability is equivalent to the condition that every line subbundle of $\mcV$ has negative degree.
On the other hand, if $\mcV$ is indecomposable, 
 the maximal degree of line subbundles is at most $g-1$.
 Taking this into account,
we introduce  the following definition (cf. ~\cite{JoPa}, ~\cite{JRXY}, ~\cite{LanPa},  ~\cite{Li}, ~\cite{Zha} for the case $\N =1$, and ~\cite{KoWa}, ~\cite{Wak22} for general  $\N \in \mbZ_{> 0}$).

\SSP
\bde \label{Def447}
Let $\mcV$ be as above, and suppose further  that $\mcV$ is stable.
Then, we say that $\mcV$ is  {\bf maximally $F^\N$-destabilized} (where $\N \in \mbZ_{> 0}$)
if the pull-back $F^{\N *}_X (\mcV)$ admits a line subbundle of degree $g-1$.
(This definition coincides with those in ~\cite{KoWa} and ~\cite{Wak22} after changing the base of $X$ along $F_{\mr{Spec}(k)}^\N$.)
Moreover, we say that $\mcV$ is  {\bf maximally $F^\infty$-destabilized}
if it is maximally $F^\N$-destabilized for some $\N \in \mbZ_{> 0}$.
\ede
\SSP

Now, let $B_1$ be the set of $k$-rational points of $SU^2_X$ classifying all maximally $F^1$-destabilized bundles.
If the pull-back $\mcU := F^{\N*}_X (\mcV)$ of $\mcV$ is again  stable, 
then $\mcV$ is maximally $F^\infty$-destabilized if and only if $\mcU$ is maximally $F^\infty$-destabilized.
Hence, the subset
\begin{align}
B_\infty := \bigcup_{\N \in \mbZ_{> 0}} ((\mr{Ver}_X^2)^\N)^{-1} (B_1) \subseteq SU_X^2 (k)
\end{align}
coincides with the set of  all maximally $F^\infty$-destabilized bundles in $SU_X^2 (k)$.
Since no point of  $B_\infty$ is periodic under the action of  $\mr{Ver}_X^2$,
there are infinitely many elements of 
$B_\infty$.

As established in ~\cite[Proposition 3.3]{JoXi} (see also ~\cite{JRXY}, ~\cite{JoPa}, ~\cite{KoWa}, and ~\cite{Wak7}),
the Frobenius pull-back functor yields a bijective correspondence between maximally $F$-destabilized bundles  and  dormant opers, and moreover this correspondence can be extended to the case of level $\N$ (cf.  the proof of Proposition \ref{Prop491} stated below).
As an application of this fact together with the result of  the previous section,
we obtain the following statement.

\SSP
\bpr \label{Prop491}
Let $\mcV$ be 
a maximally $F^\infty$-destabilized bundle  classified by $SU_X^2$.
Then, there exists a vector bundle $\mcU$ classified by $SU^2_X$ with $F^*_X (\mcU) \cong \mcV$.
In other words,
the set-theoretic image of $\mr{Ver}^2_X$  contains the set $B_\infty$.
\epr
\begin{proof}
Let us take a maximally $F^\N$-destabilized bundle $\mcV$ classified by $SU_X^2$ (where $\N \in \mbZ_{> 0}$).
By definition, the pull-back $F_X^{\N *} (\mcV)$ admits 
  a unique line subbundle $\mcL$  of degree $g-1$.
Denote by $\nabla_{\mcV, \mr{can}}^{(\N -1)}$  the canonical $\mcD^{(\N-1)}$-module structure discussed in ~\cite[Corollaire 3.3.1]{Mon}.

We claim that   $\mcL$ is not preserved under $\nabla_{\mcV, \mr{can}}^{(\N -1)}$.
Suppose, on the contrary, that  $\nabla_{\mcV, \mr{can}}^{(\N -1)}$ restricts to a $\mcD^{(\N -1)}$-module structure $\nabla_\mcL$ on $\mcL$.
Since $\nabla_{\mcV, \mr{can}}^{(\N -1)}$ has vanishing $p^\N$-curvature,
 there exists a line subbundle $\mcN$ on $\mcV$ with $F^{\N*}_X (\mcN) \cong \mcL$ (cf. ~\cite[Corollary 3.2.4]{LeQu}).
Its  degree satisfies $\mr{deg} (\mcN) = \frac{1}{p^\N} \cdot  \mr{deg} (F^{\N*}_X (\mcN))=\frac{1}{p^\N} \cdot \mr{deg}(\mcL) = \frac{g-1}{p^\N} > 0$, which contradicts the stability assumption on $\mcV$.
This completes the proof of the claim.

By this claim,   the composite
\begin{align} \label{Eq802}
\mcD^{(\N -1)}_{\leq 1} \otimes \mcL \xrightarrow{\mr{inclusion}} \mcD^{(\N -1)} \otimes  F^{\N*}_X (\mcV) \xrightarrow{\nabla_{\mcV, \mr{can}}^{(\N -1)}}
F^{\N*}_X (\mcV)
\end{align}
is injective over the generic point.
By comparing the degrees of the domain and codomain, we see that this composite  is an isomorphism.
Thus,
the collection  $(F^{\N *}_{X}(\mcV), \nabla_{\mcV, \mr{can}}^{(\N -1)}, \mcL)$ forms a dormant $\mr{GL}_2^{(\N)}$-oper, in the sense of ~\cite[Definitions 5.6, 5.7]{Wak20}.

Note that
 \begin{align}
\mcO_X &\cong F^{\N *}_X (\mcO_X) \cong F^{\N *}_X (\mr{det}(\mcV)) \cong  \mr{det}(F^{\N*}_X(\mcV)) \\
 &\cong (\mcD_{\leq 0}^{(\N -1)} \otimes \mcL) \otimes ((\mcD_{\leq 1}^{(\N -1)}/\mcD_{\leq 0}^{(\N -1)}) \otimes \mcL)\cong \mcT \otimes \mcL^{\otimes 2}.\notag
\end{align}
Hence, the trivial $\mcD^{(\N)}$-module structure on $\mcO_X$ can be transposed into  a $\mcD^{(\N)}$-module structure $\nabla_{\widetilde{\vartheta}}$
on $\mcT \otimes \mcL^{\otimes 2}$, and $\widetilde{\vartheta} := (\mcL, \nabla_{\widetilde{\vartheta}})$ defines a dormant  $2^{(\N+1)}$-theta characteristic of $X$.
Let  $\nabla_{\vartheta}$ denote the $2^{(\N)}$-theta characteristic induced by 
 $\nabla_{\widetilde{\vartheta}}$.
Under the identification $F^{\N*}_X (\mcV) = \mcF_\mcL \left(:= \mcD^{(\N -1)}_{\leq 1} \otimes \mcL \right)$ given by \eqref{Eq802},
the triple  $(F^{\N *}_{X}(\mcV), \nabla_{\mcV, \mr{can}}^{(\N -1)}, \mcL)$ defines a dormant $(\mr{GL}_2^{(\N)}, \vartheta)$-oper.
It follows from Corollary \ref{Prop31} and the bijective  correspondence \eqref{Eq107}
that there exists a dormant $(\mr{GL}_2^{(\N +1)}, \widetilde{\vartheta})$-oper  $\nabla^\diamondsuit$ whose
 associated $(\mr{GL}_2^{(\N)}, \vartheta)$-oper is isomorphic to  the above triple. 
Let  $\mcU := \mr{Ker}(\nabla^{\diamondsuit})$, which is endowed   with an $\mcO_X$-module structure via push-forward along the underlying homeomorphism of $F^{\N +1}_X$.
Then, there exists an isomorphism $F^*_X (\mcU) \xrightarrow{\sim} \mcV$.
Since  $\mr{det} (\nabla^\diamondsuit)$ is trivial, $\mr{det}(\mcU)$ is isomorphic to $\mcO_X$.
Also, the stability of $\mcV$ implies that $\mcU$ is stable.
Therefore, the vector bundle $\mcU$ satisfies the desired properties, which completes the proof of the assertion.
\end{proof}
\SSP

We now recall 
 (cf. ~\cite[Section 2.2]{dSa})
that an {\bf $F$-divided sheaf} on $X$ is a collection
\begin{align}
\{ (\mcE_m, \alpha_m) \}_{m \in \mbZ_{\geq 0}},
\end{align}
where each pair $(\mcE_m, \alpha_m)$ consists of a  vector bundle $\mcE_m$ on $X$ and an $\mcO_{X}$-linear isomorphism
$\alpha_m : F_{X}^{*} (\mcE_{m+1}) \xrightarrow{\sim} \mcE_m$.
If we are given an $F$-divided sheaf $\{ (\mcE_m, \alpha_m) \}_{m}$,
then the equalities  $\mr{deg}(\mcE_m) = \mr{deg}(F_{X}^{*} (\mcE_{m+1})) = p \cdot \mr{deg}(\mcE_{m+1})$ for every $m$ imply that 
the degree $\mr{deg}(\mcE_0)$ of $\mcE_0$ must be zero.

It is well-known that 
each $F$-divided sheaf  associates a group representation of a certain pro-algebraic group, often called the {\it stratified fundamental group} of $X$ (cf. ~\cite[Definition 7]{dSa}).
Also, such an object  corresponds to a structure of $\mcD_X^{(\infty)}$-module on $\mcE_0$ (cf. ~\cite[Theorem 1.3]{Gie}), where $\mcD_X^{(\infty)}$ denotes the sheaf of differential operators in the sense of Grothendieck (cf. ~\cite[Section 16.8.1]{Gro}). 
By applying Proposition \ref{Prop491},
we obtain  
the following assertion.

\SSP
\bco \label{Cor4}
Let  $\mcV$ be a (stable and) maximally $F^{(\infty)}$-destabilized  vector bundle on $X$ of rank $2$ and degree $0$.
Then,
there exists an $F$-divided sheaf $\{ (\mcE_m, \alpha_m) \}_{m}$ on $X$ satisfying $\mcE_0 = \mcV$.
In other words, $\mcV$ admits a $\mcD_X^{(\infty)}$-module structure extending its $\mcO_X$-module structure.
\eco
\begin{proof}
Note that we can find a line bundle $\mcL$ on $X$ satisfying $\mcL^{\otimes 2} \cong \mr{det} (\mcV)$ (hence $\mr{deg}(\mcV \otimes \mcL^\vee)$ is trivial).
Since $\mr{deg}(\mcL) = 0$,
the line bundle    $\mcL$ admits a structure of $F$-divided sheaf, i.e., there exists an $F$-divided sheaf  $\{ (\mcL_m, \beta_m) \}_{m \in \mbZ_{\geq 0}}$  with $\mcL_0 = \mcL$.
On the other hand, 
by applying Proposition  \ref{Prop491} inductively to $\mcV \otimes \mcL^\vee$,
we obtain an $F$-divided sheaf  $\{ (\mcF_m, \gamma_m) \}_{m \in \mbZ_{\geq 0}}$ with $\mcF_0 = \mcV \otimes \mcL^\vee$.
The tensor product $\gamma_m \otimes \beta_m$ (for $m \in \mbZ_{\geq 0}$) defines an isomorphism $\left( F^*_X (\mcF_{m+1}) \otimes F^*_X (\mcL_{m+1})\cong \right)F^*_X (\mcF_{m+1} \otimes \mcL_{m+1}) \xrightarrow{\sim} \mcF_m \otimes \mcL_m$.
Under the natural identification $\mcV = \mcF_0 \otimes \mcL$,
the collection $\{ (\mcF_m \otimes \mcL_m, \gamma_m \otimes \beta_m) \}_m$ forms   the desired $F$-divided  sheaf.
\end{proof}

\LSP
\subsection{Deformations to prime-power characteristic} \label{SS51}

In this subsection, we study 
deformations of  flat modules to prime-power characteristic.
Denote by $W_\N$ the ring of Witt vectors of length $\N$ over $k$ (hence $W_1 = k$).
Let us fix  a  smooth deformation $f_\N : X_\N \rightarrow \mr{Spec}(W_\N)$ of $X$ over $W_\N$.
Note that, just as in the case of $X/k$ discussed in Section \ref{SS11},
we can define the notion of a connection on an $\mcO_{X_\N}$-module.
Moreover, the notion of a $\mr{PGL}_n$-oper on $X_\N$ can be formulated in terms of connections,  see ~\cite[Definition 5.2, (i)]{Wak20}.

Let $(\mcF, \nabla)$ be a flat module on $X_\N$, i.e., a pair consisting of an $\mcO_{X_\N}$-module $\mcF$ and a connection $\nabla$ on $\mcF$.
Suppose further  that $\mcF$ is a vector bundle of rank  $n$.

\SSP
\bde[cf. ~\cite{Wak20}, Definition 3.4] \label{Def321}
We say that $(\mcF, \nabla)$ is {\bf dormant} if it is, Zariski locally on $X_\N$, isomorphic to 
$(\mcO_{X_\N}, d)^{\oplus n}$, where $d$ denotes the universal derivation on  $\mcO_{X_\N}$ over $W_\N$.
Moreover, a $\mr{PGL}_n$-oper on $X_\N$ is said to be  {\bf dormant}
if it is induced from  a dormant  flat module of rank $n$  via projectivization.
\ede
\SSP

We recall the diagonal reductions of dormant flat modules  (cf. ~\cite[Section 3]{Wak20})
in the non-logarithmic setting.
Let $(\mcF, \nabla)$ be as above, and denote by $\mcF_1$ the mod $p$ reduction of $\mcF$.
For each open subscheme $U$ of $X \left(= k \times_{W_\N} X_\N \right)$,
we define $\mcV' (U)$ to be the set of sections $v \in \mcF_1 (U)$ admitting a lifting $\widetilde{v}$ in $\mr{Ker}(\nabla) (U)$,
where we regard $U$ as an open subscheme of $X_\N$ via the underlying homeomorphism of the closed immersion $X \hookrightarrow X_\N$.
The Zariski sheaf $\mcV$ on $X$ associated to the presheaf $U \mapsto \mcV' (U)$ is endowed with an $\mcO_{X}$-module structure via push-forward along $F_{X}^\N$.
The morphism $F_X^{\N *} (\mcV) \rightarrow \mcF_1$ corresponding to the inclusion $\mcV \hookrightarrow F_{X*}^\N (\mcF_1)$ via the adjunction relation ``$F^*_X (-) \dashv F_{X*} (-)$"
turns out to be an isomorphism because of the dormancy assumption.
Through this isomorphism, $\nabla_{\mcV, \mr{can}}^{(\N -1)}$ is transferred to a  $\mcD^{(\N -1)}$-module structure ${^{\Diag\!}}\nabla$ on $\mcF_1$ with vanishing $p^\N$-curvature (cf. ~\cite[Proposition 3.11, (ii)]{Wak20}).
Thus, we obtain  a $\mcD^{(\N -1)}$-bundle
\begin{align}
{^{\Diag\!}}\msF := (\mcF_1, {^{\Diag\!}}\nabla),
\end{align}
which is called the {\bf diagonal reduction} of $\msF$ (cf. ~\cite[Definition 3.7]{Wak20}).

This construction yields a map of sets
\begin{align}
\left(\begin{matrix} \text{the set of isomorphism classes} \\ \text{of dormant $\mr{PGL}_n$-opers on $X_\N$}\end{matrix} \right)
\xrightarrow{}
\left(\begin{matrix} \text{the set of isomorphism classes} \\ \text{of dormant $\mr{PGL}_n^{(\N)}$-opers on $X_1$}\end{matrix} \right).
\end{align}
According to ~\cite[Theorem D]{Wak20}, this map becomes bijective 
when $n=2$ and $X$ is sufficiently general in the moduli stack  $\mcM_g$ (cf. ~\cite[Theorem D]{KaWa} for the case of elliptic curves).
Using this bijection and the results established  so far,
we obtain  the following assertion.

\SSP
\bpr \label{Proptt6}
Let $\N'$ be a positive integer with $\N' \leq \N$ and $\msE^\spadesuit$ a dormant $\mr{PGL}_2$-oper on $X_{\N'}$.
Suppose that $X$ is general in $\mcM_g$.
Then, 
there exists a dormant $\mr{PGL}_2$-oper $\msE_\N^\spadesuit$ on $X_{\N}$ 
whose reduction modulo $p^{\N'}$ is isomorphic to  $\msE^\spadesuit$.
\epr
\SSP

In what follows, we reformulate  the above result as a statement about  deformations of differential operators.
For  the precise definitions of the terms involved,  we refer the reader  to ~\cite[Sections 5.7-5.8]{Wak20}.

Let $\N'$ be a positive integer with $\N' \leq \N$ and  $\mcL$ a line bundle on $X_\N$ satisfying  $p^\N \mid 2g-2 + 2 \cdot \mr{deg}(\mcL)$.
 Denote by $\mcL_{\N'}$ 
 the reduction of $\mcL$ modulo $p^{\N'}$.
Suppose that  we are given a second-order linear differential operator $\msD$ on $\mcL_{\N'}$ with unit principal symbol and having a full set of root functions (cf. ~\cite[Definition 5.27, 5.31]{Wak20}).

\SSP
\bpr \label{Proptt5}
Let us keep the above notation, and 
suppose that $X$ is general in $\mcM_g$.
Then, 
there exists 
a second-order  linear differential operator  $\msD_\N$ on $\mcL$ satisfying the following conditions:
\begin{itemize}
\item[(i)]
$\msD_{\N}$  has unit principal symbol and  has  a full set of root functions; 
\item[(ii)]
The mod $p^{\N'}$ reduction of $\msD_{\N}$ coincides with $\msD$.
\end{itemize}
\epr
\begin{proof}
Write $\mcT := \mcT_{X/k}$ and $\mcT_{\N'} := \mcT_{X_{\N'}/W_{\N'}}$ for simplicity.
Since $X$ is assumed  to be general, we may  assume, without loss of generality, that
$X$ is ordinary in the usual sense.
There exists a dormant $2^{(1)}$-theta characteristic $\vartheta := (\varTheta, \nabla_\vartheta)$ of $X_{\N'}$  (cf. ~\cite[Definition 5.12]{Wak20} for the definition of a dormant $2^{(1)}$-theta characteristic in prime-power characteristic) such that $\varTheta^\vee = \mcL_{\N'}$ 
and $\msD$ defines a $(2, \vartheta)$-projective connection (cf. the proof of ~\cite[Theorem 10.25]{Wak20}).

Denote by $(\mcT \otimes \varTheta_1^{\otimes 2}, {^{\Diag\!}}\nabla_\vartheta)$ the diagonal reduction of $(\mcT_{\N'}\otimes \varTheta_{\N'}^{\otimes 2}, \nabla_\vartheta)$.
We  regard  $\mcS ol ({^{\Diag\!}}\nabla_\vartheta)$ as an $\mcO_X$-module via push-forward along the underlying homeomorphism of $F_X^{\N -\N'}$.
Since $F_X^{\N - \N' *} (\mcS ol ({^{\Diag\!}}\nabla_\vartheta)) \cong \mcT \otimes \varTheta_1^{\otimes 2}$ and $p^\N \mid \mr{deg}(\mcT \otimes \varTheta_1^{\otimes 2})$,
  the degree of the line bundle $\mcS ol ({^{\Diag\!}}\nabla_\vartheta)$ 
  is divisible by $p^{\N'}$.
Hence, there exists a line bundle $\mcN$ with $F^{\N' *}_X (\mcN) \cong \mcS ol ({^{\Diag\!}}\nabla_\vartheta)$, 
which yields  an isomorphism  $F^{\N*}_X (\mcN)\xrightarrow{\sim} \mcT \otimes \varTheta_1^{\otimes 2}$.
By this isomorphism, $\nabla_{\mcN, \mr{can}}^{(\N -1)}$ corresponds to   a $\mcD^{(\N -1)}$-module structure $\nabla_{\vartheta, \N}$ on  $ \mcT \otimes \varTheta_1^{\otimes 2}$.
Then, it follows from ~\cite[Proposition 3.19]{Wak20}
that $(\mcT \otimes \varTheta_1^{\otimes 2}, \nabla_{\vartheta, \N})$ admits a unique  diagonal  lifting of the form $(\mcT_\N \otimes (\mcL^\vee)^{\otimes 2}, {^{\Diagg\!}}\nabla_{\vartheta, \N})$ for some connection ${^{\Diagg\!}}\nabla_{\vartheta, \N}$.
The mod $p^{\N'}$ reduction of ${^{\Diagg\!}}\nabla_{\vartheta, \N}$ coincides with
$\nabla_{\vartheta}$.
In particular, the pair $\vartheta_\N := (\mcL^\vee, {^{\Diagg\!}}\nabla_{\vartheta, \N})$ defines a dormant $2^{(1)}$-theta characteristic of $X_\N$.

Next, let $\msE^\spadesuit_{\N'}$ be the dormant $\mr{PGL}_2^{(\N')}$-oper on $X_1$ obtained as the diagonal reduction of the dormant $\mr{PGL}_2$-oper corresponding to the $(2, \vartheta)$-projective connection $\msD$, through   the chain of  bijections established  in ~\cite[Section 5.8]{Wak20}.
According to Corollary \ref{Prop31},
there exists a dormant $\mr{PGL}_2^{\N}$-oper $\msE^\spadesuit_\N$ on $X$ whose level-$\N'$ reduction is isomorphic to $\msE_{\N'}^\spadesuit$.
This also  induces a canonical diagonal lifting over  $X_\N$ (cf. ~\cite[Theorem-Definition 9.4]{Wak20}), which 
corresponds to a $(2, \vartheta_\N)$-projective connection $\msD_\N$.
By the uniqueness of canonical diagonal liftings, the mod $p^{\N'}$ reduction of $\msD_\N$ must coincide  with $\msD_{\N'}$.
Thus, $\msD_\N$ satisfies the  desired properties, completing  the proof of this assertion.
\end{proof}

\vspace{10mm}
\section{Appendix: Generic Miura $\mr{PGL}_2^{(\N)}$-opers and Igusa-structures} \label{S300}
\LSP

In this appendix, we discuss dormant generic Miura $\mr{PGL}_2^{(\N)}$-opers, which are certain variants of dormant $\mr{PGL}_2^{(\N)}$-opers; see ~\cite{KaWa} or ~\cite{Wak98} for the study of such objects.
We show  that dormant generic Miura $\mr{PGL}_2^{(\N)}$-opers on an elliptic curve correspond, in a suitable  sense, to Igusa-structures of level $p^\N$ on that curve (cf. Proposition \ref{Propt201}).
In particular, this correspondence allows us to interpret  the moduli stacks $\mcO p^{^\mr{Zzz...}}_{\N, \rho, g, r}$  (at least, in the case $n=2$) 
as higher-genus generalizations of  the  classical Igusa curves.

Throughout this section, we assume that $p >2$.

\LSP
\subsection{Dormant generic Miura $\mr{PGL}_2^{(\N)}$-opers} \label{SSApp1}

Denote by $\mcM_{\mr{ell}}$ the dense open  substack of $\overline{\mcM}_{1, 1}$
classifying smooth curves.
That is, $\mcM_{\mr{ell}}$ is the moduli stack of elliptic curves (equipped with a distinguished point) over $k$.
Let $\mcM_\mr{ell}^\mr{ord}$ (resp., $\mcM_{\mr{ell}}^\mr{ss}$) denote the open (resp., closed) substack of $\mcM_\mr{ell}$ classifying ordinary  (resp., supersingular) elliptic curves; see ~\cite[Chapter 12]{KaMa} for the study of these moduli spaces.
We also denote by $\msX_\mr{univ} := (X_\mr{univ}, \{ \sigma_\mr{univ} \})$ the universal family of elliptic curves over $\mcM_{\mr{ell}}$

Let $S$ be a scheme over $k$ and  $\msX := (f : X \rightarrow S, \{ \sigma \})$ a pointed curve classified by  an $S$-rational point of 
$\mcM_{\mr{ell}}$.
The family of curves  $X$ carries 
a structure of group scheme over $S$ with   $\sigma$ taken  as the identity section.
To simplify the notation, we set  $\mcD^{(\N -1)} := \mcD^{(\N -1)}_{X/S}$ and $\mcT := \mcT_{X/S}$.

\SSP
\begin{rem}\label{Remark9}
Note that many of the constructions  on pointed stable curves considered   in this paper extend naturally to the case where  the hyperbolic condition $2g-2+r>0$ is not assumed.
For example, one can still define  a dormant $\mr{PGL}_2^{(\N)}$-oper on a proper smooth curve of genus $\leq 1$.
In what follows, when  we refer  to dormant $\mr{PGL}_2^{(\N)}$-opers   on ``$X$" (rather than on the pointed curve  ``$\msX$" as before),
we assume  that the underlying curve is equipped with the trivial log structure,  in  contrast to the framework  discussed in Section \ref{SS11}.

According to ~\cite[Proposition 7.2]{Wak20} (extended to the situation where the condition $2g-2+r >0$ is not assumed),
each dormant $\mr{PGL}_2^{(\N)}$-oper   on $X$ can be identified with a dormant $\mr{PGL}_2^{(\N)}$-oper   on $\msX$ whose radius at the marked point $\sigma$ is represented by the multiset  $[\overline{0}, \overline{1}, \cdots, \overline{n-1}]$. 
\end{rem}
\SSP

Recall that the algebraic group $\mr{PGL}_2$ may  be identified with the automorphism group of the projective line $\mbP^1$.
Under  this identification,
we define  
 $\mr{PGL}_2^\mbA$ to be  the algebraic subgroup of $\mr{PGL}_2$ consisting of all automorphisms $h$ of $\mbP^1$ that fix the point at infinity, i.e.,  $h (\infty) = \infty$.
Each element of $\mr{PGL}_2^\mbA$ may be represented by a matrix  of the form $\begin{pmatrix} a & 0  \\ b & c \end{pmatrix}$ for some $a,  c \in k^\times$ and $b \in k$.
Moreover,  such a representative may be chosen with $a =1$.
We also set $B^\mbA = B \cap \mr{PGL}_2^\mbA$.

Let us consider a pair $\hat{\msE}^\spadesuit := (\mcE_{B^\mbA}, \phi^\mbA)$ consisting of a $B^\mbA$-bundle $\mcE_{B^\mbA}$ on $X$ and an $(\N -1)$-PD stratification $\phi^\mbA$ on the associated $\mr{PGL}_2^\mbA$-bundle $\mcE^\mbA := \mcE_{B^\mbA} \times^{B^\mbA} \mr{PGL}_2^\mbA$ (with respect to the trivial log structure on $X$); see  ~\cite[Definition 2.3]{Wak20} for the definition of a higher-level PD stratification).
By changing the structure group via   the inclusion $B^\mbA \hookrightarrow B$
 (resp.,  $\mr{PGL}_2^\mbA \hookrightarrow \mr{PGL}_2$),
$\mcE^\mbA$ (resp., $\phi^\mbA$) induces a $B$-bundle $\mcE_B$ (resp., an $(\N -1)$-PD stratification $\phi$ on $\mcE := \mcE_B \times^B \mr{PGL}_2$).

\SSP
\bde \label{Def334}
\begin{itemize}
\item[(i)]
We say that $\hat{\msE}^\spadesuit$  is a {\bf  generic Miura $\mr{PGL}_2^{(\N)}$-oper}
(or a {\bf generic Miura $\mr{PGL}_2$-oper of level $\N$}) on $X$ if the associated  pair $(\mcE_{B}, \phi)$ forms a $\mr{PGL}_2^{(\N)}$-oper on $X$.
The notion of an isomorphism between two  generic Miura $\mr{PGL}_2^{(\N)}$-opers can also be defined in a natural way.
\item[(ii)]
A generic Miura $\mr{PGL}_2$-oper $\hat{\msE}^\spadesuit := (\mcE_{B^\mbA}, \phi^\mbA)$ is said to be {\bf dormant}
if the $p^\N$-curvature of $\phi^\mbA$ (in the sense of ~\cite[Definition 2.10, (ii)]{Wak20})
vanishes identically. 
\end{itemize}
\ede
\SSP

\begin{rem} \label{Remark4}
Let us  now describe  (dormant) generic Miura $\mr{PGL}_2^{(\N)}$-opers in terms of $\mcD^{(\N -1)}$-modules (cf. ~\cite[Section 4.3]{KaWa}).
Consider a collection of data
\begin{align} \label{Ew23d}
\hat{\msF}^\heartsuit := (\mcF, \nabla_\mcF, \mcL, \mcN),
\end{align}
where $\mcF$ is  
 a rank $2$ vector bundle on $X$, $\nabla_\mcF$ is a $\mcD^{(\N -1)}$-module structure on $\mcF$, and  $\mcL$, $\mcN$ are line subbundles of $\mcF$ such that  $\mcF = \mcL \oplus \mcN$ and  $\mcN$ is closed under $\nabla_\mcF$.
Suppose that  the composite
\begin{align} \label{Ew24d}
\mcD_{\leq 1}^{(\N -1)} \otimes \mcL \xrightarrow{\mr{inclusion}} \mcD^{(\N -1)} \otimes \mcF \xrightarrow{\nabla_\mcF} \mcF
\end{align}
 is an isomorphism (i.e., $(\mcF, \nabla_\mcF, \mcL)$ defines a $\mr{GL}_2^{(\N)}$-oper).
Then,  we refer to the collection
$\hat{\msF}^\heartsuit$
 as a {\bf generic Miura $\mr{GL}_2^{(\N)}$-oper} on $X$ (cf. ~\cite[Definition 4.8, (i)]{KaWa}).

Let us take a generic Miura $\mr{GL}_2^{(\N)}$-oper.
$\hat{\msF}^\heartsuit$ as above. 
 The vector bundle $\mcF$ together with $\mcN$ induces, via the quotient $\mr{GL}_2 \twoheadrightarrow \mr{PGL}_2$,  a $\mr{PGL}_2^\mbA$-bundle $\mcE^\mbA$ on $X$.
The line subbundle $\mcL$  specifies a $B^\mbA$-reduction $\mcE_{B^\mbA}$ of $\mcE^\mbA$.
Since $\mcN$ is preserved  by  $\nabla_\mcF$, 
we obtain an $(\N -1)$-PD stratification $\phi^\mbA$ on $\mcE^\mbA$ induced  from $\nabla_\mcF$.  
By the assumption that  \eqref{Ew24d} is an isomorphism,
 the resulting pair 
 \begin{align} \label{Ew25d}
 \hat{\msF}^{\heartsuit \Rightarrow \spadesuit} := (\mcE_{B^\mbA}, \phi^\mbA)
 \end{align}
 forms a generic Miura $\mr{PGL}_2^{(\N)}$-oper (cf. ~\cite[Section 4.7]{KaWa}).
If $\nabla_\mcF$ has vanishing $p^\N$-curvature, then  $\hat{\msF}^{\heartsuit \Rightarrow \spadesuit}$ is dormant.
Moreover, it is immediately verified that every  (dormant) generic Miura $\mr{PGL}_2^{(\N)}$-oper arises in this way from some collection
 $\hat{\msF}^\heartsuit$.
\end{rem}
\SSP

Denote  by $\mcM \mcO p^{^\mr{Zzz...}}_{\N}$ the set-valued contravariant functor on $\mcS ch_{/\mcM_\mr{ell}}$ that  assigns, to
each $S$-rational point of $\mcM_\mr{ell}$ classifying $\msX$ as above, the set of isomorphism classes of dormant generic Miura $\mr{PGL}_2^{(\N)}$-opers on $X$.
Applying  the representability of  the stacks $\mcO p^{^\mr{Zzz...}}_{\N, \rho, g, r}$ (cf. ~\cite[Theorem 6.17]{Wak20} and the last comment in Remark \ref{Remark9}),
we see  that $\mcM \mcO p^{^\mr{Zzz...}}_{\N}$ can be represented by a Deligne-Mumford stack of finite type over $k$.

\LSP
\subsection{Igusa-structures on elliptic curves} \label{SSApp2}

Let $V^{(\N)}_{X/S} : X^{(\N)} \rightarrow X$ denote  the dual isogeny  to the $\N$-th relative Frobenius morphism $F^{(\N)}_{X/S}$, which is a cyclic $p^\N$-isogeny (cf. ~\cite[Lemma 12.2.3]{KaMa}).

\SSP
\bde[cf. ~\cite{KaMa}, Definition 12.3.1] \label{Def3821}
An {\bf Igusa-structure of level $p^\N$} on $X$ is an $S$-rational point $P$ of $X^{(\N)}$ that  generates the  
kernel of $V^{(\N)}_{X/S}$.
\ede
\SSP

Denote by $\mcI g_{\N}$ the set-valued 
contravariant functor on $\mcS ch_{/\mcM_\mr{ell}}$ that  assigns, to each $S$-rational point of $\mcM_\mr{ell}$ classifying  $\msX$ as above, the set of  Igusa-structures of level $p^\N$ on $X$.
Then, $\mcI g_{\N}$  can be  represented by a connected  smooth Deligne-Mumford stack of dimension $1$ over $k$ equipped with the finite and faithfully flat morphism  $\mcI g_\N \rightarrow \mcM_\mr{ell}$ given by forgetting the data of Igusa-structures  (cf. ~\cite[Theorem 12.6.1, Corollary 12.6.2]{KaMa}).

In what follows, we construct  dormant generic Miura $\mr{PGL}_2^{(\N)}$-opers by using Igusa-structures of level $p^\N$. 
Let $\msX$ be as before.

To begin with, 
let 
 $\mr{Pic}^0_{X/S}$ denote the relative Jacobian  of $X/S$, which classifies degree-$0$ line bundles on $X$ rigidified along the section $\sigma$.
There exists an isomorphism of group $S$-schemes $\alpha_{X} : X \xrightarrow{\sim} \mr{Pic}^0_{X/S}$
given by $x \mapsto \mcO_X (x -\sigma)$.
The base-change of $\alpha_X$ along $F^\N_{S} : S \rightarrow S$ defines  an isomorphism $\alpha_{X^{(\N)}} : X^{(\N)} \xrightarrow{\sim} \mr{Pic}^0_{X^{(\N)}/S}$.
Under the identifications $X = \mr{Pic}^0_{X/S}$ and  $X^{(\N)} = \mr{Pic}^0_{X^{(\N)}/S}$ given by $\alpha_{X}$ and $\alpha_{X^{(\N)}}$, respectively,
the isogeny  $V^{(\N)}_{X/S} : X^{(\N)} \rightarrow X$ can be identified with   the morphism $F^{(\N)*}_{X/S} (-) : \mr{Pic}^0_{X^{(\N)}/S} \rightarrow \mr{Pic}^0_{X/S}$ given by $\mcL \mapsto F^{(\N)*}_{X/S} (\mcL)$.

Now, let us fix  an Igusa-structure $P$ of level $p^\N$ on $X$.
Under the identification $F^{(\N)*}_{X/S} (-) = V_{X/S}^{(\N)}$ described   above,
$P$ determines a line bundle $\mcL$ on $X^{(\N)}$ admitting an isomorphism $\gamma : F^{(\N)*}_{X/S}(\mcL) \xrightarrow{\sim} \mcO_X$.
Recall that $F^{(\N)*}_{X/S}(\mcL)$ carries a canonical $\mcD^{(\N-1)}$-module structure with vanishing $p^\N$-curvature such that the  subsheaf of horizontal sections coincides with
   $(F^{(\N)}_{X/S})^{-1}(\mcL)$
 (cf. ~\cite[Corollaire 3.3.1]{Mon}).
 This corresponds to  a $\mcD^{(\N -1)}$-module structure 
 \begin{align} \label{Erfh91}
 \nabla_P : {^L}\mcD^{(\N -1)} \rightarrow \mcE nd_{\mcO_S} (\mcO_X)
 \end{align}
  on $\mcO_X$ via $\gamma$,
 and this construction  does not depend on the choice of $\gamma$.

\SSP
\begin{exa} \label{Remark5}
Let us examine   the $\mcD^{(\N -1)}$-module structure  $\nabla_P$ in the case  $\N =1$,   where  the underlying curve is taken to be the universal formal deformation of a supersingular elliptic curve.

Let us take a $k$-rational point $q$ of $\mcM_\mr{ell}$ lying on the supersingular locus $\mcM_{\mr{ell}}^\mr{ss}$.
Denote by  $S_0$  the  formal neighborhood of  $q$  in $\mcM_\mr{ell}$ identified with    $S_0 = \mr{Spec}(k[\![t]\!])$, and by 
 $\msX_0 := (X_0/S_0, \{ \sigma_0 \})$ 
the restriction  of the universal curve $\msX_\mr{univ}$  over $S_0$.
According to ~\cite[Corollary 12.4.4]{KaMa},
there exists a nowhere vanishing invariant $1$-form $\omega_0$ on  $X_0/S_0$ satisfying $C_{X_0/S_0} (\omega_0) = t \cdot \omega_0^{(1)}$, where $C_{(-)}$ denotes the Cartier operator of the relative  smooth curve ``$(-)$" (i.e., the inverse to the operator ``$C^{-1}$'' discussed in ~\cite[Theorem 7.2]{NKat0}), and $\omega^{(1)}_0$ denotes the base-change of $\omega_0$ along the absolute Frobenius morphism $F_{S_0}$.

Now, we set  $S := \mr{Spec} (k[\![u]\!])$, where $u := t^{1/(p-1)}$.
The natural inclusion $k[\![t]\!] \hookrightarrow k[\![u]\!]$ gives rise to a $k$-morphism $s : S \rightarrow S_0$.
Denote by $\msX := (X/S, \{ \sigma \})$ the base-change of $\msX_0$   along $s$.
By putting $\omega := u \cdot s^*(\omega_0)$, we have 
\begin{align}
C_{X/S} (\omega) = C_{X/S} (u \cdot s^*(\omega_0)) = u \cdot C_{X/S} (s^*(\omega_0)) = u \cdot (t \cdot s^*(\omega_0^{(1)})) = u^{p} \cdot  (s^{*}(\omega_0))^{(1)} = \omega^{(1)}.
\end{align}
 It follows from ~\cite[Corollary 7.1.3]{NKa1} that
the connection $\nabla_\omega$  on $\mcO_X$ defined by $\nabla_\omega := d + \omega$  has vanishing $p$-curvature.
Under the identification $X^{(1)} = \mr{Pic}^0_{X^{(1)}/S}$ given by $\alpha_{X^{(1)}}$,
the line bundle $\mcS ol (\nabla_\omega)$ on $X^{(1)}$ corresponds to an Igusa-structure of level $p$ on $X$.
Hence, it determines an $S$-morphism 
\begin{align} \label{Ewi38}
S \rightarrow \mcI g_1 \times_{\mcM_\mr{ell}} S_0.
\end{align}
From the local structure of the projection $\mcI g_1 \rightarrow \mcM_\mr{ell}$ established  in ~\cite[Corollary 12.6.2, (1)]{KaMa},
 we deduce  that this morphism is an isomorphism.
In particular,  if $P$ denotes
the  universal   Igusa-structure  over  the base-change of $\msX_\mr{univ}$ over  $\mcI g_1 \times_{\mcM_\mr{ell}} S_0$,
then the equality 
\begin{align} \label{erk23}
\nabla_\omega = \nabla_P
\end{align}
 holds via the isomorphism \eqref{Ewi38}. 
\end{exa}
\SSP

By using $\nabla_P$ constructed in \eqref{Erfh91},   we obtain a collection of data
\begin{align} \label{Effki}
\hat{\msF}^\heartsuit_P := (\mcO_X^{\oplus 2}, \nabla_\mr{triv} \oplus \nabla_P, \mr{Im}(\Delta), \mr{Im}(\iota_2)),
\end{align}
where $\nabla_\mr{triv}$ denotes the trivial $\mcD^{(\N -1)}$-module structure on $\mcO_X$,  $\Delta$ denotes the diagonal embedding $\mcO_X \hookrightarrow \mcO_X^{\oplus 2}$, and $\iota_2$ denotes the inclusion into the second factor $\mcO_X \hookrightarrow \mcO_X^{\oplus 2}$.

 Suppose further  that  the family of elliptic curves   $X/S$ is ordinary.
 Since $P$ restricts to  a nonzero point  in each fiber of the structure morphism   $f : X \rightarrow S$,
 the pair of  global sections  $(1, 1)$, $(\nabla_\mr{triv} (1), \nabla_P (1))$  forms a basis of  $\mcO_X^{\oplus 2}$.
 It follows that 
 $\hat{\msF}^\heartsuit_P$ gives rise to   a dormant generic Miura $\mr{PGL}_2^{(\N)}$-oper $\hat{\msF}^{\heartsuit \Rightarrow \spadesuit}_P$
 via the construction  in  Remark \ref{Remark4}.
The resulting assignment  $P \mapsto \hat{\msF}^{\heartsuit \Rightarrow \spadesuit}_P$ determines a bijective correspondence
\begin{align} \label{Ey23q}
\begin{pmatrix}\text{The set of Igusa-structures} \\ \text{of level $p^\N$ on $X$} \end{pmatrix}
\xrightarrow{\sim}
\begin{pmatrix}\text{The set of isomorphism classes} \\ \text{of dormant generic Miura} \\ \text{$\mr{PGL}_2^{(\N)}$-opers on $X$} \end{pmatrix}
\end{align}
(cf. ~\cite[Proposition 4.14]{KaWa}).
This bijection is functorial with respect to  $S$, and thus we obtain the following assertion.

\SSP
\bpr \label{Propt201}
\begin{itemize}
\item[(i)]
For every $\N \geq 1$, 
the assignment $(\msX, P) \mapsto (\msX, \hat{\msE}_P^\spadesuit)$ defines 
 an isomorphism of $\mcM_\mr{ell}^\mr{ord}$-stacks 
 \begin{align} \label{Eii811}
 \mcI g_\N  \times_{\mcM_\mr{ell}} \mcM_\mr{ell}^\mr{ord} \xrightarrow{\sim} \mcM \mcO p^{^\mr{Zzz...}}_\N  \times_{\mcM_\mr{ell}} \mcM_\mr{ell}^\mr{ord}.
 \end{align}
 \item[(ii)]
In the case $\N =1$, 
the isomorphism  \eqref{Eii811} extends to an isomorphism of $\mcM_{\mr{ell}}$-stacks
\begin{align}
\mcI g_1 \xrightarrow{\sim} \mcM \mcO p^{^\mr{Zzz...}}_1.
\end{align}
\end{itemize}
\epr
\begin{proof}
Since 
assertion (i) follows from the functorial isomorphism \eqref{Ey23q}, we focus here on proving  assertion (ii).

To  simplify the notation, we set $S := \mcI g_1$,
 and let 
$(\msX, P)$ (where $\msX := (X/S, \{ \sigma \})$) be the universal object defined over $S$.
Write 
 $\nabla_{P, +} := \nabla_\mr{triv} \oplus \nabla_P$.
Also, denote by $\mcG$ the $\mcO_X$-submodule of $\mcO_X^{\oplus 2}$ generated by the two global sections $(1, 1)$, $\nabla_{P, +}((1, 1)) \left(= (0, \nabla_P (1)) \right)$.
The $\mcO_X$-submodule $\mcL$ (resp., $\mcN$) of $\mcG$ generated by  $(1, 1)$ (resp., $(0, \nabla_P (1))$) specifies a line subbundle of $\mcF$.

In what follows, we claim that  $\mcG$ is closed under $\nabla_{P, +}$, and that  the collection $\hat{\msG}^{\heartsuit}_{P} := (\mcG, \nabla_{P, +} |_{\mcG}, \mcL, \mcN)$ forms a generic Miura $\mr{GL}_2^{(1)}$-oper.
Note that  the restriction of $\hat{\msG}^{\heartsuit}_{P}$ to the base-change $\msX_\mr{univ}^\mr{ord}$ of $\msX_\mr{univ}$ over $\mcM_\mr{ell}^\mr{ord}$ coincides with  $\hat{\msF}^\heartsuit_{P^\mr{ord}}$,
  where $P^\mr{ord}$ denotes the restriction of $P$ to $\msX_\mr{univ}^\mr{ord}$.
Hence, it suffices to examine the fiber over  a supersingular curve.
 Let $q$ be a $k$-rational point of $\mcM_\mr{ell}$ lying on $\mcM_\mr{ell}^\mr{ss}$, and $\widetilde{q}$ a unique point of $\mcI g_1$ over $q$.
 Denote by $U$ the formal neighborhood of $\widetilde{q}$ in $\mcI g_1$, and identify $U$ with $\mr{Spec}(k[\![u]\!])$ for some formal parameter $u$.
 Let  $(\msX_U. P_U)$, where $\msX_U := (X_U/U, \{ \sigma_U \})$, be the restriction of $(\msX, P)$ over $U$.
 As discussed in Example  \ref{Remark5},
 there exists, after possibly replacing the parameter $u$ with another,  a nowhere vanishing $1$-form $\omega_0$ on $X_U/U$ satisfying $\nabla_P = d + u \cdot \omega_0$.
Hence,  $\mcG |_{X_U}$ admits  a basis
$\{ (1, 1), (0, u) \}$.
It follows that $\mcG |_{X_U}$ is closed under $\nabla_{P, +} |_{X_U}$, and 
the restriction of $\hat{\msG}^{\heartsuit}_P$ to $X_U$ forms a generic Miura $\mr{GL}_2^{(\N)}$-oper.
This  
completes the proof of the claim.

Consequently, the assignment $P \mapsto \hat{\msG}^{\heartsuit \Rightarrow \spadesuit}_P$ (cf. \eqref{Ew25d}) defines a morphism $\mcI g_1 \rightarrow \mcM \mcO p_1^{^\mr{Zzz...}}$ extending the isomorphism \eqref{Eii811} in the case  $\N =1$.
 Note  that $\mcI g_1$ is smooth over $k$ (cf. ~\cite[Theorem 12.6.1]{KaMa}).
 On the other hand, by  an argument similar to the proof of the smoothness assertion in ~\cite[Theorem 6.3.2]{Wak98},
 $\mcM \mcO p^{^\mr{Zzz...}}_1$ is verified to be smooth over $k$.
Using these facts, we conclude  that the morphism $\mcI g_1 \rightarrow \mcM \mcO p_1^{^\mr{Zzz...}}$ just  constructed  is an isomorphism.
 We finish the proof of this assertion.
  \end{proof}
\SSP

\begin{rem}
As shown in ~\cite[Proposition 5.2.1, (ii)]{Wak8},
there exists  no dormant generic Miura $\mr{PGL}_2^{(\N)}$-oper on any supersingular curve when $\N \geq 2$.
On the other hand, the fiber of the projection $\mcI g_\N \rightarrow \mcM_{\mr{ell}}$ over the  point classifying any such  curve is nonempty  (cf. ~\cite[Corollary 12.6.2]{KaMa}).
Therefore, if $\N \geq 2$,
then the isomorphism \eqref{Eii811} cannot extend to an isomorphism between the entire spaces $\mcI g_\N \xrightarrow{\sim} \mcM \mcO p_\N^{^\mr{Zzz...}}$.
\end{rem}

\LSP
\subsection*{Acknowledgements}
We are grateful for the many constructive conversations we had with
the moduli space of dormant $\mr{PGL}_2^{(\N)}$-opers, who lives in the world of mathematics!
The author was partially supported by 
 JSPS KAKENHI Grant Number 25K06933.

\end{document}